\documentclass[12pt,a4paper]{article}
\usepackage{amscd,amsmath,amssymb,amsfonts,times,color}
\usepackage{mathrsfs}
\usepackage[all]{xy}

   \newtheorem{thm}{Theorem}[section]
   \newtheorem{lem}[thm]{Lemma}
   \newtheorem{prop}[thm]{Proposition}
   \newtheorem{cor}[thm]{Corollary}
   \newtheorem{conj}[thm]{Conjecture}
   \newtheorem{caution}[thm]{Caution}
   \newtheorem{defn}[thm]{Definition}
   
   \newtheorem{exmp}[thm]{Example}
   \newtheorem{rem}[thm]{Remark}
\numberwithin{equation}{subsection}

\begin{document}
\title{\'Etale cohomology of arithmetic schemes \endgraf and zeta values of arithmetic surfaces}
\author{Kanetomo Sato
\footnote{Supported by Grant-in-Aid for Scientific Research (C) 16K05072 \endgraf
}
}
\date{}
\maketitle

\long\def\remind#1{\textcolor[named]{Peach}{}\textcolor[named]{Blue}{{#1}}}

\newcommand{\qed}
{\mbox{}\nolinebreak$\square$}
\newenvironment{pf}{\par\noindent\emph{Proof.}\;}{\hfill\qed\par\smallskip}
\newenvironment{pf*}[1]{\par\noindent\emph{#1.}\;}{\hfill\qed\par\medskip}
\def\smallsymbol#1{%
\mbox{{{\small$#1$}}}}

\newcommand{\bA}{\mathbb A}
\newcommand{\bC}{\mathbb C}
\newcommand{\bF}{\mathbb F}
\newcommand{\bG}{\mathbb G}
\newcommand{\bL}{\mathbb L}
\newcommand{\bN}{\mathbb N}
\newcommand{\bZ}{\mathbb Z}
\newcommand{\bP}{\mathbb P}
\newcommand{\bQ}{\mathbb Q}
\newcommand{\bR}{\mathbb R}
\newcommand{\bW}{\mathbb W}

\newcommand{\cA}{\mathscr A}
\newcommand{\cB}{\mathscr B}
\newcommand{\cC}{\mathscr C}
\newcommand{\cD}{\mathscr D}
\newcommand{\cE}{\mathscr E}
\newcommand{\cF}{\mathscr F}
\newcommand{\cG}{\mathscr G}
\newcommand{\cH}{\mathscr H}
\newcommand{\cK}{\mathcal K}
\newcommand{\cL}{\mathscr L}
\newcommand{\cM}{\!\!{\mathscr M}}
\newcommand{\cN}{\mathscr N}
\newcommand{\cO}{\mathscr O}
\newcommand{\cS}{\mathcal S}
\newcommand{\cT}{\mathscr T}
\newcommand{\cU}{\mathscr U}
\newcommand{\cX}{\mathscr X}
\newcommand{\cY}{\mathscr Y}
\newcommand{\cZ}{\mathscr Z}

\newcommand{\fe}{\mathfrak e}
\newcommand{\fo}{\mathfrak o}
\newcommand{\fO}{\mathfrak O}
\newcommand{\fH}{\mathfrak H}
\newcommand{\fS}{\mathfrak S}
\newcommand{\fT}{\mathfrak T}
\newcommand{\fX}{\mathfrak X}
\newcommand{\fm}{\mathfrak m}
\newcommand{\fM}{\mathfrak M}

\newcommand{\ep}{\epsilon}
\newcommand{\lam}{\lambda}
\newcommand{\Lam}{\Lambda}
\newcommand{\vare}{\varepsilon}
\newcommand{\vG}{\varGamma}
\newcommand{\vL}{\varLambda}

\newcommand{\Ab}{{\mathscr A}\hspace{-1pt}b}
\newcommand{\ac}{\hspace{1.5pt}{\cdot}\hspace{1.5pt}}
\newcommand{\aj}{\text{\rm aj}}
\newcommand{\Alb}{\text{\rm Alb}}
\newcommand{\alb}{\text{\rm alb}}
\newcommand{\alg}{\text{\rm alg}}
\newcommand{\an}{\text{\rm an}}
\newcommand{\Aut}{\text{\rm Aut}}
\newcommand{\Br}{\text{\rm Br}}
\newcommand{\bull}{*}
\newcommand{\cd}{\text{\rm cd}}
\newcommand{\ch}{\text{\rm ch}}
\newcommand{\CH}{\text{\rm CH}}
\newcommand{\Char}{\text{\rm Char}}
\newcommand{\cl}{\text{\rm cl}}
\newcommand{\Cone}{\text{\rm Cone}}
\newcommand{\cont}{\text{\rm cont}}
\newcommand{\Cor}{\text{\rm Cor}}
\newcommand{\corank}{\text{\rm corank}}
\newcommand{\cosp}{\text{\rm cosp}}
\newcommand{\cotor}{\text{\rm cotor}}
\newcommand{\crys}{\text{\rm crys}}
\newcommand{\Cris}{\text{\rm Cris}}
\newcommand{\Cube}{\text{\bf Cube}}
\renewcommand{\det}{\text{\rm det}\hspace{1pt}}
\newcommand{\ECube}{\text{\bf ECube}}
\newcommand{\wCube}{\text{\bf wCube}}
\newcommand{\EwCube}{\text{\bf EwCube}}
\newcommand{\codim}{\text{\rm codim}}
\newcommand{\Coker}{\text{\rm Coker}}
\newcommand{\degn}{\text{\rm degn}}
\renewcommand{\dim}{\text{\rm dim}}
\newcommand{\Div}{\text{\rm Div}}
\renewcommand{\div}{\text{\rm div}}
\newcommand{\divi}{\text{\rm div}}
\newcommand{\dR}{\text{\rm dR}}
\newcommand{\dRW}{\text{\rm dRW}}
\newcommand{\DR}{\text{\rm DR}}
\newcommand{\End}{\text{\rm End}}
\newcommand{\et}{\text{\rm \'et}}
\newcommand{\Et}{\text{\rm \'Et}}
\newcommand{\Ext}{\text{\rm Ext}}
\newcommand{\ez}{\text{\rm \hspace{-1pt}-ez}}
\newcommand{\Fib}{\text{\rm Fib}}
\newcommand{\Fr}{\text{\rm Fr}}
\newcommand{\Frac}{\text{\rm Frac}}
\newcommand{\Gal}{\text{\rm Gal}}
\newcommand{\Gr}{\text{\rm Gr}}
\newcommand{\Gys}{\text{\rm Gys}}
\renewcommand{\H}{\text{\it H}\hspace{1.2pt}}
\newcommand{\h}{\text{\rm h}}
\newcommand{\Hom}{\text{\rm Hom}}
\renewcommand{\hom}{\text{\rm hom}}
\newcommand{\id}{\text{\rm id}}
\newcommand{\Image}{\text{\rm Im}}
\newcommand{\incl}{\text{\rm incl}}
\newcommand{\Jac}{\text{\rm Jac}}
\newcommand{\Ker}{\text{\rm Ker}}
\newcommand{\loc}{\text{\rm loc}}
\newcommand{\logcrys}{\text{\rm log-crys}}
\newcommand{\nd}{\!\!\!\not|\sp}
\newcommand{\Nis}{\text{\rm Nis}}
\newcommand{\Ob}{\text{\rm Ob}}
\newcommand{\op}{\text{\rm op}}
\newcommand{\ord}{\text{\rm ord}}
\newcommand{\ovf}{\!\hspace{-1pt}/\hspace{-.8pt}f}
\newcommand{\Pic}{\text{\rm Pic}}
\newcommand{\pr}{\text{\rm pr}}
\newcommand{\rank}{\text{\rm rank}}
\renewcommand{\Re}{\text{\rm Re}}
\newcommand{\red}{\text{\rm red}}
\newcommand{\reg}{\text{\rm reg}}
\newcommand{\res}{\text{\rm Res}}
\newcommand{\Res}{\text{\rm Res}}
\newcommand{\rk}{\text{\rm rk}}
\newcommand{\Sch}{\text{\bf Sch}}
\newcommand{\Sets}{\text{\bf Sets}}
\newcommand{\sExt}{\cE \! xt}
\newcommand{\sgn}{\text{\rm sgn}}
\newcommand{\sh}{\text{\rm sh}}
\newcommand{\sHom}{\cH \! om}
\newcommand{\sing}{\text{\rm sing}}
\newcommand{\Spec}{\text{\rm Spec}}
\newcommand{\st}{\text{\rm st}}
\newcommand{\supp}{\text{\rm supp}}
\newcommand{\Supp}{\text{\rm Supp}}
\newcommand{\sw}{\text{\rm sw}}
\newcommand{\Tor}{\text{\rm Tor}}
\newcommand{\Tot}{\text{\rm Tot}}
\newcommand{\tors}{\text{\rm tors}}
\newcommand{\tr}{\text{\rm tr}}
\newcommand{\Tr}{\text{\rm Tr}}
\newcommand{\ur}{\text{\rm ur}}
\newcommand{\vol}{\text{\rm vol}}
\newcommand{\vD}{\varDelta}
\newcommand{\Zar}{\text{\rm Zar}}

\newcommand{\tK}{\text{\it K}}
\newcommand{\tL}{\text{\it L}}

\newcommand{\lla}{\longleftarrow}
\newcommand{\lra}{\longrightarrow}
\newcommand{\Lra}{\Longrightarrow}
\newcommand{\Llra}{\Longleftrightarrow}
\newcommand{\hra}{\hookrightarrow}
\newcommand{\Ra}{\Rightarrow}
\newcommand{\sq}{\square}
\newcommand{\qis}{\text{\rm qis}}
\newcommand{\hloc}{\text{\rm -loc}}

\def\bs#1{\boldsymbol{#1}}
\def\lim#1{\us{#1}{\varinjlim}}
\def\spt{\sptilde}
\def\sp{\hspace{.7pt}}
\def\nsp{\hspace{-.7pt}}
\def\ssm{\smallsetminus}
\def\ol#1{\overline{#1}}
\def\os#1#2{\overset{#1}{#2}}
\def\ul#1{\underline{#1}}
\def\us#1#2{\underset{#1}{#2}}
\def\wt#1{\widetilde{#1}}
\def\wh#1{\widehat{#1}}

\def\OK{O_{\hspace{-1.1pt} K}}

\def\rF{\text{\rm F}}
\def\rf{\text{\sl f}\sp}

\def\tC{\wt{C}}
\def\tj{\wt{j}}
\def\tP{\wt{P}}
\def\tV{\wt{V}}
\def\tW{\wt{W}}
\def\tX{\wt{X}}
\def\tY{\wt{Y}}
\def\tom{\wt{\omega}}

\def\sha{
\raisebox{0.059cm}{\underline{\phantom{\hspace{13.3pt}}}}\hspace{-13.3pt}\text{\rm{I\hspace{0.6pt}I\hspace{0.6pt}I}}
\hspace{-13.17pt}\raisebox{0.059cm}{\underline{\phantom{\hspace{13.3pt}}}}\hspace{.7pt}}

\def\ra{\rightarrow}
\def\lra{\longrightarrow}
\def\Lra{\Longrightarrow}
\def\la{\leftarrow}
\def\lla{\longleftarrow}
\def\Lla{\Longleftarrow}
\def\da{\downarrow}
\def\hra{\hookrightarrow}
\def\lmt{\longmapsto}
\def\sm{\setminus}
\def\Gm{{\mathbb G}_{\hspace{-1pt}\mathrm{m}}}
\def\Ga{{\mathbb G}_{\hspace{-1pt}\mathrm{a}}}
\def\L{\bZ/p^n\bZ}
\def\QlZl{\bQ_\ell/\bZ_\ell}
\def\qp{\bQ_p}
\def\zp{\bZ_p}
\def\QpZp{\qp/\zp}
\def\QZ'{\bQ/\bZ'}
\def\Ln{\varLambda_n}
\def\Tn{\fT_n}
\def\Xn{X_{\hspace{-.2pt}n}}
\def\Wn{W_{\hspace{-1pt}n}}
\def\Ws{W\hspace{-1.5pt}}
%
%
%
\def\bigwitt#1#2#3{\bW_{\hspace{-2pt}#2}{\hspace{1pt}}\Omega_{#1}^{#3}}
\def\mwitt#1#2#3{W_{\hspace{-2pt}#2}{\hspace{1pt}}\omega_{#1}^{#3}}
\def\mlogwitt#1#2#3{W_{\hspace{-2pt}#2}{\hspace{1pt}}\omega_{{#1}\hspace{-.7pt},\hspace{.3pt}{\log}}^{#3}}
\def\mwwitt#1#2#3{W_{\hspace{-2pt}#2}{\hspace{1pt}}\tom_{#1}^{#3}}
\def\hwitt#1#2#3{W_{\hspace{-2pt}#2}{\hspace{1pt}}\Xi_{#1}^{#3}}
\def\cwitt#1#2#3{W_{\hspace{-2pt}#2}{\hspace{1pt}}\Lambda_{#1}^{#3}}
\def\witt#1#2#3{W_{\hspace{-2pt}#2}{\hspace{1pt}}\Omega_{\hspace{-.5pt}#1}^{#3}}
\def\logwitt#1#2#3{W_{\hspace{-2pt}#2}{\hspace{1pt}}\Omega_{\hspace{-.5pt}{#1}\hspace{-.7pt},\hspace{.3pt}{\log}}^{#3}}
\def\II{I\hspace{-0.5pt}I}
\def\III{I\hspace{-0.7pt}I\hspace{-0.7pt}I}
\def\IV{I\hspace{-0.7pt}V}
\def\VII{V\hspace{-0.7pt}I\hspace{-0.7pt}I}
\def\XVIII{X\hspace{-0.7pt}V\hspace{-0.7pt}I\hspace{-0.7pt}I\hspace{-0.7pt}I}

\def\pmn{\par\medskip\noindent}
\def\pbn{\par\bigskip\noindent}
\def\vp{v}

\def\sfrac#1#2{\smallsymbol{\dfrac{\sp{#1}\sp}{\sp{#2}\sp}}}
\vspace{-25pt}
\begin{abstract}
In this paper, we deal with the \'etale cohomology of a proper regular arithmetic scheme $X$ with $\zp(r)$ and $\qp(r)$-coefficients, where the coefficients are complexes of \'etale sheaves that the author introduced in \cite{Sa0}.
We will prove that the \'etale cohomology of $X$ with $\qp(r)$-coefficients agrees with the Selmer group of Bloch-Kato for any $r \geqq \text{\rm dim}(X)$.
Using this fundamental result, we further discuss an approach to the study of zeta values (or residue) at $s=r$, via the \'etale cohomology with $\zp(r)$-coefficients, relating Tamagawa number conjecture of Bloch-Kato with a zeta value formula.
As a consequence, we will obtain an unconditional example of an arithmetic surface for which the residue of its zeta function at $s=2$ is computed modulo rational numbers prime to $p$, for infinitely many $p$'s.
\end{abstract}

{\small
\begin{quote}
\end{quote}}


\section{Introduction}\label{sect0}
Let $K$ be an algebraic number field, and let $O_K$ be its integer ring.
Let $X$ be a regular connected scheme which is proper flat over $B:=\Spec(O_K)$, and such that $X_{\!K}=X \otimes_{O_K} K$ is geometrically connected over $K$.
We fix a prime number $p$, and assume that
\begin{quote}
($\star$)\, {\it $X$ has good or log smooth reduction at all places $v$ of $K$ dividing $p$.}
\end{quote}
In this paper, we give a new approach to the values or residue of the zeta function of $X$ at integers $r \geqq \dim (X)$ using \'etale cohomology of $X$ with $\qp(r)$ and $\zp(r)$-coefficients,
 cf.\ \cite{KCT}, \cite{Li2}, \cite{Mo}, \cite{FM}, \cite{FS}.

\subsection{Selmer groups}
Let $\H^1_{\!f}(K,V^i(r))$ be the Selmer group of Bloch-Kato associated with the $p$-adic Galois representation $V^i(r):=\H^i(X_{\ol K},\qp)(r)$.
The first aim of this paper is to relate this group with the \'etale cohomology group $\H^{i+1}(X,\qp(r))$, assuming that $r \geqq d:=\dim(X)$.
Here $\H^*(X,\qp(r))$ is defined as
\[ \H^*(X,\qp(r)):=\qp \otimes_{\zp} \varprojlim_{n \geqq 1} \ \H^*(X,\fT_n(r)) \]
and $\fT_n(r)$ ($n \geqq 1$) denotes the complex of \'etale $\bZ/p^n\bZ$-sheaves on $X$ introduced in \cite{Sa0} under the assumption that $X$ has good or semi-stable reduction at all places $v$ dividing $p$;
we have $\H^*(X,\fT_n(r)) \cong \H^*(X[p^{-1}],\mu_{p^n}^{\otimes r})$ when $r>d$.
See \S\ref{sect1} below for details on this object under the setting of this paper.
The first main result of this paper is the following comparison
(cf.\ \cite{FM} Proposition 5.18, \cite{Sa1} Theorem 9.1, \S10):
\begin{thm}\label{thm1-1}
Assume that $r \geqq d$. Then we have
\[ \H^{i+1}(X,\qp(r)) \cong
\begin{cases}
 \qp & \hbox{{\rm(}when $(i,r)=(2d,d)${\rm)}}, \\
 \H^1_{\!f}(K,V^i(r))
 \quad & \hbox{{\rm(}otherwise{\rm)}}.
\end{cases}
\]
\end{thm}
The key idea of Theorem \ref{thm1-1} is as follows.
By a duality result of Jannsen-Saito-Sato \cite{JSS} and the adjunction between $R\pi_{X\!/\!B !}(=R\pi_{X\!/\!B *})$ and $R\pi_{X\!/\!B}^!$, we have
\begin{equation}\label{eq.intro-1}
 R \pi_{X\!/\!B*}\fT_n(r)_X \cong R\sHom_{B,\,\bZ/p^n\bZ}(R\pi_{X\!/\!B !}\fT_n(d-r)_X,\fT_n(1)_B)[2-2d]
\end{equation}
in $D^+(B_{\et},\bZ/p^n\bZ)$ (see Lemma \ref{lem-direct} below), where the assumption $r \geqq d$ is crucial and $\fT_n(d-r)_X$ is a constructible sheaf placed in degree $0$ by definition.
Using this fact, we introduce the following complexes:
\begin{align*}
\fH^{\geqq i}(X,\fT_n(r)) & := R\sHom_{B,\,\bZ/p^n\bZ}
  (\tau_{\leqq 2d-2-i}R \pi_{X\!/\!B!}\fT_n(d-r)_X,\fT_n(1)_B)[2-2d],\phantom{\big|_{\big|}}\\
\fH^i(X,\fT_n(r)) & :=
  R\sHom_{B,\,\bZ/p^n\bZ}(R^{2d-2-i}\pi_{X\!/\!B!}\fT_n(d-r)_X,\fT_n(1)_B).
\end{align*}
By the proper base change theorem for $R \pi_{X\!/\!B!}$, we have
\[ \fH^i(X,\fT_n(r)) = 0 \quad \hbox{ unless } \;\; 0 \leqq i \leqq 2d-2, \]
and the filtration $\{ \fH^{\geqq i}(X,\fT_n(r)) \}_i$ on the right hand side of \eqref{eq.intro-1} yields a convergent spectral sequence
\begin{equation}\notag
E^{a,i}_2=\H^a(B,\fH^i(X,\fT_n(r))) \Lra \H^{a+i}(X,\fT_n(r)).
\end{equation}
The $E_2$-terms of this spectral sequence are finite (see Proposition \ref{prop1-1} below), and we obtain the following spectral sequence of finite-dimensional $\bQ_p$-vector spaces:
\begin{equation}\label{eq.intro-2}
E^{a,i}_2=\H^a(B,\fH^i(X,\bQ_p(r))) \Lra \H^{a+i}(X,\bQ_p(r)),
\end{equation}
where
\[ \H^a(B,\fH^i(X,\bQ_p(r))):=\qp \otimes_{\zp} \varprojlim_{n \geqq 1} \ \H^a(B,\fH^i(X,\fT_n(r))).\]
Concerning the spectral sequence \eqref{eq.intro-2}, we will prove
\begin{thm}[\S\ref{sect5}]\label{thm1}
Assume $r \geqq d$.
Then the $\bQ_p$-vector space $E_2^{a,i}$ is zero, unless $a=1$ or $(a,i,r)=(3,2d-2,d)$.
Consequently, the spectral sequence \eqref{eq.intro-2} degenerates at $E_2$-terms.
Moreover, we have
\[ E_2^{1,i} \cong \H^1_{\!f}(K,V^i(r)) \]
for any $i$ and $r \geqq d$, which is zero unless $0 \leqq i \leqq 2d-2$.
We have $E_2^{3,2d-2}\cong \qp$, if $r=d$.
\end{thm}
Theorem \ref{thm1-1} is a consequence of this result.
An important point of Theorem \ref{thm1} is the vanishing of $E_2^{2,i}$ for any $i$, which we will prove by
computing the cohomology of all local integer rings with $\fH^i(X,\bQ_p(r))$-coefficients and by a local-global argument using a Hasse principle of Jannsen \cite{J} p.\ 337, Theorem 3\,(c).
As a consequence of the vanishing of $E_2^{2,i}$ (and $E_2^{3,i}$ with $(i,r) \ne (2d-2,d)$), we obtain the following result on Galois cohomology (cf.\ \cite{J} p.\ 317, Conjecture 1, p.\ 349, Question 2, \cite{Fl} \S3, \cite{Ki1} 1.1.7, \cite{Li} 9.1, \cite{So1} Th\'eor\`eme 5):
\begin{cor}[Corollary \ref{cor5-4}\,(2)]\label{cor-main}
Let $S$ be a finite set of places of $K$ including all places which divide $p\ac \infty$ or where $X$ has bad reduction. Assume $r \geqq d$. Then the restriction map
\[ \H^2(G_{\!S},V^i(r)) \lra  \bigoplus_{v \in S} \ \H^2(K_v,V^i(r)) \]
is bijective for any $(i,r) \ne (2d-2,d)$, and injective for $(i,r)=(2d-2,d)$.
In particular, if $r > d$ or $X_{\nsp K}$ has potentially good reduction at all finite places of $K$, then we have
\[ \H^2(G_{\!S},V^i(r))=0 \quad \hbox{ for any $(i,r) \ne (2d-2,d)$}. \]
\end{cor}

\subsection{$\bs{p}$-Tate-Shafarevich groups\quad ($\bs{d=2}$)}
We assume that $X$ is an arithmetic surface, i.e., $d=2$.
Put $T^i:=\H^i(X_{\ol{K}},\zp)$.
In their paper \cite{BK} \S5, Bloch and Kato introduced a homomorphism
\begin{equation}\label{eq-BKmap}
 \alpha^{i,r} : \frac{\H^1(K,T^i \otimes \QpZp(r))}{\sp\H^1_{\!f}(K,T^i \otimes \QpZp(r))\sp}\lra \bigoplus_{v \in P} \ \frac{\H^1(K_v,T^i \otimes \QpZp(r))}{\sp\H^1_{\!f}(K_v,T^i \otimes \QpZp(r))\sp},
\end{equation}
where $P$ denotes the set of all places of $K$, and for each $v \in P$, $K_v$ denotes the local field of $K$ at $v$;
 $\H^1_{\!f}(K,T^i \otimes \QpZp(r))$ (resp.\ $\H^1_{\!f}(K_v,T^i \otimes \QpZp(r))$) denotes the image of $\H^1_{\!f}(K,V^i(r))$ (resp.\ $\H^1_{\!f}(K_v,V^i(r))$).
The cokernel $\Coker(\alpha^{i,r})$ is finite and canonically isomorphic to the Pontryagin dual of $\H^{2-i}(X_{\ol{K}},\qp/\zp(2-r))^{G_{\nsp\nsp K}}$, if $i-2r \leqq -3$.
They also proved that $\Ker(\alpha^{i,r})=:\sha^{(p)}(\H^i(X_K)(r))$, the $p$-Tate-Shafarevich group of the motive $\H^i(X_K)(r)$, is finite for the same $(i,r)$.
The second main result of this paper compares the maps $\alpha^{i,r}$ with $p$-adic Abel-Jacobi mappings
\[ \aj_p^{i,r} : \H^i_{\cM}(X,\bZ(r)) \, \wh{\otimes} \, \zp \lra \H^1_{\!f}(K,T^{i-1}(r))  \]
assuming $r \geqq 2$.
Here $\H^*_{\cM}(X,\bZ(r))$ denotes the motivic cohomology of $X$ (see \S\ref{sect6-1} below), and
for an abelian group $M$, $M \,\wh{\otimes}\, \bZ_p$ denotes its $p$-adic completion
     $\varprojlim{}_n  \ M/p^n$.
We will calculate the above Abel-Jacobi mapping using the Merkur'ev-Suslin theorem \cite{MS} and the Rost-Voevodsky theorem \cite{V1}, \cite{V2}, which together with Theorem \ref{thm1} will play important roles in the following comparison formula:

\begin{thm}[\S\ref{sect6}\,+\,Corollary \ref{cor7-0}]\label{thm2}
Assume $r \geqq 2$, and that $p \geqq 3$ or $B(\bR) = \emptyset$.
Assume further that $\H^3_{\cM}(X,\bZ(r))\{p\}$, the $p$-primary torsion part of $\H^3_{\cM}(X,\bZ(r))$, is finite.
Let $S'$ be the set of the finite places of $K$ which divide $p$ or where $X$ has bad reduction.
Then $\aj_p^{i,r}$ has finite kernel and cokernel for $i=2,3$, and we have
{\allowdisplaybreaks
\begin{align*}
\frac{\chi(\alpha^{1,2})}{\,\chi(\alpha^{0,2})\,}
 & =
  \dfrac{\chi(\aj^{3,2}_p)}{\sp \chi(\aj^{2,2}_p) \sp}
  \ac \dfrac{\sp \#\sp \CH_0(X)\{ p \}\sp}{\sp \#\sp \Pic(\OK)\{ p \} \sp}
  \ac \prod_{v \in S'} \ \dfrac{\, e_v^{2,1,2} \ac e_v^{3,0,2}\,}{e_v^{2,0,2} \ac e_v^{3,1,2}}
  \quad  \phantom{\bigg|_{\big|}} & \hbox{\rm ($r=2$)\,}  \\
\frac{\chi(\alpha^{1,r})}{\sp \chi(\alpha^{0,r}) \ac \chi(\alpha^{2,r})\sp}
 & =
  \dfrac{\,\chi(\aj^{3,r}_p)\,}{\,\chi(\aj^{2,r}_p) \,} \ac \#\sp \H^4_{\cM}(X,\bZ(r))\{p\}
  \ac \prod_{v \in S'} \ \dfrac{\, e_v^{2,1,r} \ac e_v^{3,0,r} \ac e_v^{3,2,r}\,}{e_v^{2,0,r} \ac e_v^{2,2,r} \ac e_v^{3,1,r}}
  \quad & \hbox{\rm ($r \geqq 3$),}
\end{align*}
}where we put $\chi(f):=\#\sp \Coker(f)/\#\sp \Ker(f)$ for a homomorphism $f : M \to N$ of abelian groups with finite kernel and cokernel{\rm;}
for each $v \in S'$ and $a=2,3$, we put
\[ e_v^{a,i,r}:=\# \sp \H^a(B_v,\fH^i(X,\bZ_p(r))), \quad
 \hbox{$B_v:=$ the completion of $B$ at $v$.} \]
See Corollary \ref{cor4-fin}\,{\rm(}2{\rm)} below for the finiteness of $e_v^{a,i,r}$.
\end{thm}
\noindent
The finiteness of $\CH_0(X)$ is due to Bloch \cite{B1}, Kato and Saito \cite{KS}.
By the localization theorem of Levine \cite{Le}, $\H^i_{\cM}(X,\bZ(r))$ is zero for any $i > r+2$ (see Lemma \ref{lem6-2}\,(1) below).
As natural extensions of these facts, we will prove that $\H^4_{\cM}(X,\bZ(r))\{p\}$ is finite for any $r \geqq3$, and that $\H^i_{\cM}(X,\bZ(r))$ is uniquely $p$-divisible for any $i \geqq 5$ and $r \geqq3$, see Propositions \ref{cor6-1} and \ref{cor6-2} below. The formulas in Theorem \ref{thm2} are based on these facts and results.

\subsection{Zeta values modulo rational numbers prime to $\bs{p}$\quad ($\bs{d=2}$)}
Assuming a weak version of $p$-Tamagawa number conjecture (see Conjecture \ref{conj7-1}), we will relate
the formula in Theorem \ref{thm2} with the residue or value at $s=r$ of the zeta function
\[ \zeta(X,s) := \prod_{x \in X_0} \ \frac{1}{\sp 1 - q_x^{-s} \sp } \qquad \hbox{($q_x := \#\kappa(x)$)}, \]
where the product on the right hand side runs through all closed points of $X$ and converges absolutely for any $s$ with $\Re(s) > 2\sp(=d)$. Recall that $\zeta(X,s)$ is meromorphically continued to $\Re(s) > 3/2$ and has a simple pole at $s=2$, see \cite{Se1}.

\begin{thm}[Proposition \ref{prop7-1}]\label{thm3}
Assume $r \geqq 2$ and the following conditions{\rm:}
\begin{itemize}
\item[{\rm (i)}]
$p \geqq r+2$.
\item[{\rm (ii)}]
For any $v \in B_0$ dividing $p$, $v$ is absolutely unramified and $X$ has good reduction at $v$.
\item[{\rm (iii)}]
A weak $p$-Tamagawa number conjecture {\rm(}see Conjecture \ref{conj7-1} below{\rm)} holds for the motives $H^i(X_K)(r)$ with $i=0,1$ {\rm(}resp.\ $i=0,1,2${\rm)}, if $r=2$ {\rm(}resp.\ $r \geqq 3${\rm)}.
\end{itemize}
Then $\H^3_{\cM}(X,\bZ(r))\{p\}$ is finite, and we have
{\allowdisplaybreaks
\begin{align*}
& \displaystyle \us{s=2}{\Res} \ \zeta(X,s) \equiv
\us{s=1}{\Res} \ \zeta_K(s) \ac
\frac{\, \chi(\aj_p^{3,2}) \ac \#\sp \CH_0(X) \ac R^{0,2}_\Phi  \,}{\sp \chi(\aj_p^{2,2}) \ac \#\sp \Pic(\OK) \ac R^{1,2}_\Phi \sp}
\;\; \mod \, \bZ_{(p)}^\times \phantom{\big|_{\Big|}} & \hbox{{\rm (}$r = 2${\rm)}} \\
& \zeta(X,r) \equiv
\frac{\,\chi(\aj_p^{3,r}) \ac \#\sp \H^4_{\cM}(X,\bZ(r))\{p\} \ac R^{0,r}_\Phi \ac R^{2,r}_\Phi \,}{\sp \chi(\aj_p^{2,r}) \ac R^{1,r}_\Phi \sp}
\;\; \mod \, \bZ_{(p)}^\times & \hbox{{\rm (}$r \geqq 3${\rm)}}
\end{align*}
where $\bZ_{(p)}$ denotes the localization of \sp $\bZ$ at $(p)$.
See Conjecture \ref{conj7-1} below for the definition of the number $R^{i,r}_\Phi \in \bR^\times/\bZ_{(p)}^\times$, which is a $p$-adic modification of the Beilinson regulator of the motive $H^i(X_K)(r)$.
}
\end{thm}
\noindent
This result is deduced from Theorem \ref{thm2} and certain comparison results between the alternating products of local terms that appear in Theorem \ref{thm2} with zeta values of the closed fibers of $X \to B$, see Theorems \ref{thm7-1} and \ref{thm7-2} below.
The assumptions (i) and (ii) are essential in this comparison at present, while the reduction at the closed points $v \in B_0$ with $v \nd p$ is arbitrary.

\begin{exmp}
{\rm Let $K$ be an imaginary quadratic field, and let $E$ be an elliptic curve over $K$ with complex multiplication by the integer ring $O_K$ of $K$.
Let $D$ (resp.\ $w$) be the discriminant of $K$ (resp.\ the number of roots of unity contained in $K$).
Let $X$ be a regular model of $E$ which is proper flat over $O_K$.
Let $p$ be a prime number which is prime to $6$ and {\it good for $X$} in the sense that $X$ has good reduction at each place of $K$ lying above $p$.
Then we obtain a formula (without assuming any conjectures)
\[ \us{s=2}\Res \ \zeta(X,s) \equiv 
\frac{\,2\pi \ac \chi(\aj_p^{3,2}) \ac \#\CH_0(X) \ac R_\Phi^{0,2} \,}{\sp w \sqrt{-D\sp} \ac \chi(\aj_p^{2,2}) \ac R_\Phi^{1,2} \sp}
\;\; \mod  \bZ_{(p)}^\times \]
from Corollary \ref{cor-main}, Theorem \ref{thm3} and results of Kings \cite{Ki1} Theorem 1.1.5 and Huber-Kings \cite{HK} Theorem 1.3.1 (see also \cite{Ki2} Theorems 2.1.3 and 2.2.2).
If we assume that $\H^{i+1}_{\cM}(X,\bZ(2))$ is a finitely generated abelian group for $i=0,1,2$, then we have
\[ \rank_{\bZ} \, \H^1_{\cM}(X,\bZ(2)) = 1, \qquad \rank_{\bZ} \, \H^2_{\cM}(X,\bZ(2)) = 2,\qquad
 \# \H^3_{\cM}(X,\bZ(2)) < \infty \]
by Theorem \ref{thm1} (and Proposition \ref{lem6-1}, Corollary \ref{lem6-0}\,(1)\sp--\sp(3) below), and obtain a stronger formula
\[ \us{s=2}\Res \ \zeta(X,s) \equiv
\frac{\, 2\pi \ac R_{\sp\cM}^{0,2} \ac \#\Ker(\reg_{\sp\cD}^{2,2}) \ac \#\CH_0(X) \,}
{\sp \sqrt{-D\sp} \ac \#\Ker(\reg_{\sp\cD}^{1,2}) \ac R_{\sp\cM}^{1,2} \ac \# H^3_{\cM}(X,\bZ(2))\sp}
\;\; \mod  \bZ[T^{-1}]^\times  \]
by Theorem \ref{thm7-3} below, where $T$ denotes the set of all prime numbers which divide $6$ or which are bad for $X$; $\reg_{\sp\cD}^{i+1,2}$ for $i=0,1$ denotes the regulator map to the real Deligne cohomology with $\bZ(2)$-coefficients
\[ \reg_{\sp\cD}^{i+1,2}  : \H^{i+1}_{\cM}(X,\bZ(2)) \lra \H^{i+1}_{\cD}(E_{/\bR},\bZ(2)). \]
For $i=0,1$,
$R_{\sp\cM}^{i,2}$ denotes the volume of $\Coker(\reg_{\sp\cD}^{i+1,2})$ with respect to the same $\bZ$-lattice of $\H^i_{\dR}(E/K)$ as used in the definition of $R^{i,2}_\Phi$.
}
\end{exmp}

\subsection*{Organization of this paper}
In \S\ref{sect1}, we review the definition of the \'etale complexes $\fT_n(r)$ on $X_\et$ and establish their fundamental properties under the setting of this paper.
In \S\ref{etale}--\S\ref{sect-limits}, we further introduce the \'etale complexes $\fH^{\geqq i}(X,\fT_n(r))$ and $\fH^i(X,\fT_n(r))$ on $B_\et$ assuming $r \geqq d$ and prove some preliminary results on those new complexes.
In \S\ref{sect3}--\S\ref{sect5} we will prove Theorems \ref{thm1-1} and \ref{thm1}.
In \S\ref{sect6}, we will compute $p$-adic cycle class maps and $p$-adic Abel-Jacobi mappings assuming $r \geqq d=2$,
 and then prove the formulas in Theorem \ref{thm2}.
In \S\ref{sect7}, we will relate the alternating product of local terms in Theorem \ref{thm2} with zeta values of fibers of $X \to B$.
Finally in \S\ref{sect8}, we will relate the formulas in Theorem \ref{thm2} with zeta values assuming a weak version of $p$-Tamagawa number conjecture.

\subsection{Notation}
Throughout this paper, we fix a prime number $p$, and put $\Ln := \bZ/p^n\bZ$.
\par
If $p$ is {\it invertible} on a scheme $X$, we write $\mu_{p^n}=\mu_{p^n\!\nsp,\sp X}$ ($n \geqq 1$) for the \'etale sheaf of $p^n$-th roots of unity on $X$, and define a $\Ln$-sheaf $\Ln(r)=\Ln(r)_X$ ($r \in \bZ$) on $X_\et$ as
\begin{equation}\label{eqdef1-1}
 \Ln(r) :=
 \begin{cases}
 \mu_{p^n}^{\otimes r} \quad\phantom{\big|_{|}} & (r \geqq 1) \\
 \Ln \quad\phantom{\big|_{|}} & (r=0) \\
 \sHom(\Ln(-r),\Ln) \quad & (r < 0).
 \end{cases}
\end{equation}
This notation will be useful mainly in the case that $r$ is negative.
\par
On the other hand, if $X$ is {\it an $\bF_p$-scheme}, then we write $\logwitt X n r$ ($r \geqq 0$, $n \geqq 1$) for the \'etale subsheaf of the logarithmic part of the Hodge-Witt sheaf $\witt X n r$ (see \cite{il} I (1.12.1)).
If $r < 0$, then we define $\logwitt X n r$ as the zero sheaf.
If $X$ is an equi-dimensional scheme which is of finite type over a field $k$ of characteristic $p$, then we write  $\nu_{X\hspace{-.7pt},\hspace{.3pt}n}^r$ for the sheaf on $X_\et$ defined as the kernel of Kato's boundary map \cite{KCT}
\[  \partial : \bigoplus_{x \in X^0} \ i_{x*} \logwitt x n r \lra  \bigoplus_{x \in X^1} \ i_{x*}\logwitt x n {r-1}, \]
where $i_x : x \to X$ denotes the canonical map for any $x \in X$.
If $X$ is smooth over $k$, then we have $\nu_{X\hspace{-.7pt},\hspace{.3pt}n}^r = \logwitt X n r$ by Gros-Suwa \cite{GS} and Shiho \cite{Sh}.
\par
Unless indicated otherwise, all cohomology groups of schemes are taken over the \'etale topology. 

\section{\'Etale coefficients}\label{sect1}
Let $\fO$ be a Dedekind ring whose fraction field $K$ has characteristic $0$, and let $p$ be a prime number.
We put
\[ B:=\Spec(\fO), \qquad B[p^{-1}]:=\Spec(\fO[p^{-1}]) \quad \hbox{ and } \quad 
 \Sigma:=\Spec\Big(\fO\big/\!\sqrt{(p)}\Big). \]
Let $X$ be a regular connected scheme which is separated, flat of finite type over $B=\Spec(\fO)$.
For a closed point $v \in B$, we put $B_v^\loc:=\Spec(\fO_v^\loc)$ and $Y_v:=X \times_B v$, where $\fO_v^\loc$ denotes the localization of $\fO$ at $v$.
Throughout this paper, we assume
\begin{itemize}
\item[($\star_1$)]
{\it for any $v \in \Sigma$, the reduced part $(Y_v)_\red$ of $Y_v$ has normal crossings on $X$ and the morphism $X \times_B B_v^\loc \to B_v^\loc$ is log smooth with respect to the log structure on $X \times_B B_v^\loc$ associated with $(Y_v)_\red$ and that on $B_v^\loc$ associated with $v$.}
\end{itemize}
See also Remark \ref{rem-log} below for a remark on this assumption.
We write $\pi_{X\!/\!B} : X \to B$ for the structure morphism, and put $d:=\dim(X)$, the absolute dimension of $X$. 
Let $Y$ be the disjoint union of $Y_v$'s for all $v \in \Sigma$.
Let $j$ (resp.\ $\iota$) be the open immersion $X[p^{-1}] \hra X$ (resp.\ closed immersion $Y \hra X$).

In this section, we define a family of complexes of \'etale sheaves
 $\{\Tn(r)\}_{n \geqq 1\nsp,\sp r \in \bZ}$ on $X$ and check several fundamental properties of them using the main results of \cite{SS}, which have been established in \cite{Sa0} and \cite{Sa1} in the case that $X$ has semi-stable reduction at all $v \in \Sigma$.
The coefficients $\{\Tn(r)\}_{n,\sp r}$ play key roles throughout this paper.

\subsection{\'Etale complex $\Tn(r)$}\label{sect1-1}
For $r \geqq 0$, we define a complex $\Tn(r)=\Tn(r)_X \in D^b(X_\et,\Ln)$ by the distinguished triangle
\begin{equation}\label{eq1-1-0}
 \iota_*\nu_{Y\!,\hspace{.3pt}n}^{r-1}[-r-1] \os{g}\lra \Tn(r) \os{t}\lra \tau_{\leqq r}Rj_*\mu_{p^n}^{\otimes r} \os{\sigma}\lra \iota_*\nu_{Y\!,\hspace{.3pt}n}^{r-1}[-r].
\end{equation}
See \cite{Sa0} (3.2.5) and (4.2.1) for the morphism $\sigma$.
By the same arguments as in loc.\ cit.\ 4.2.2, $\Tn(r)$ is concentrated in $[0,r]$, and the pair $(\Tn(r),t)$ is unique up to a unique isomorphism.
For $r<0$, we define $\Tn(r)$ as
\[ \Tn(r) := j_! \Ln(r). \]
See \eqref{eqdef1-1} for the definition of the (locally constant) sheaf $\Ln(r)$ on $(X[p^{-1}])_\et$.
\begin{lem}\label{lem1-0}
\begin{enumerate}
\item[{\rm (1)}]
If $p$ is invertible in $\fO$, then we have $\Tn(r) \cong \Ln(r)$
 for any $r \in \bZ$.
\item[{\rm (2)}]
Assume that
\begin{itemize}
\item[{\rm($\star_2$)}]
 {\it any residue field of $\fO$ of characteristic $p$ is perfect.}
\end{itemize}
Then we have $\Tn(r) \cong Rj_*\Ln(r)=Rj_*\mu_{p^n}^{\otimes r}$ for any $r > d$.
\end{enumerate}
\end{lem}
\begin{pf}
(1) is obvious. We prove (2).
Without loss of generality, we may assume that $\fO$ is local and strict henselian.
Let $k$ be the residue field of $\fO$.
Since $k$ is algebraically closed by assumption, we have $\cd_p(K) \leqq 1$ (\cite{se} Chapter \II, \S3.3).
By this fact and the cohomological dimension of affine varieties \cite{sga4} X.3.2, we have $\tau_{\leqq r}Rj_*\mu_{p^n}^{\otimes r} \cong Rj_*\mu_{p^n}^{\otimes r}$ for any $r \geqq d$.
On the other hand, we have $\nu_{Y\!,\hspace{.3pt}n}^{r-1}=0$ for any $r > d$ again because $k$ is algebraically closed (note that $\dim(Y)=d-1$). The assertion follows from these facts.
\end{pf}

\begin{rem}\label{rem-log}
{\rm
Under ($\star_2$) of Lemma \ref{lem1-0}\,(2),
one does not need the log-smoothness assumption ($\star_1$) to define $\Tn(r)$ for $r >d$, but has only to define $\Tn(r) := Rj_*\mu_{p^n}^{\otimes r}$.
Moreover, one can check that all the results in \S\S\ref{sect1}\sp--\sp\ref{sect5} (resp.\ in \S\ref{sect6}) with `$(r)$' in coefficients (e.g.\ $\Tn(r)$, $\fH^*(X,\Tn(r))$, $\fH^*(X,\zp(r))$, $\fH^*(X,\qp(r))$, $\zp(r)$, $\qp(r)$ and $V^m(r)$, in particular, Corollary \ref{cor5-4})
hold true for $r>d$ (resp.\ $r>2$) without the assumption ($\star_1$), by similar arguments to those in this paper.
We will not get into those details, but leave them to the reader as exercises to simplify the presentation.
}
\end{rem}

\begin{prop}[cf.\ \cite{Sa0} 4.2.8]\label{thm-contra}
Let $\fO'$ be another Dedekind ring which is flat over $\fO$, and
let $X'$ be a scheme which is regular and flat of finite type over $B'$ and satisfies {\rm(}$\star_1${\rm)} over $B'$.
Let $f : X' \to X$ be an arbitrary morphism, and let $g : X'[p^{-1}] \to X[p^{-1}]$ be the induced morphism.
Then for any $n \geqq 1$ and $r \in \bZ$, there exists a unique morphism
\[  f^\sharp : f^*\Tn(r)_X \lra \Tn(r)_{X'} \quad \hbox{ in } \;\; D^b(X'_\et,\Ln) \]
that extends the natural isomorphism $g^*\Ln(r)_{X[p^{-1}]} \cong \Ln(r)_{X'[p^{-1}]}$ on $X'[p^{-1}]$.
\end{prop}
\begin{pf}
The case $r \leqq 0$ is obvious. Assume $r \geqq 1$ and put
$U^1\cO_{\!X}^\times:=\Ker\big(\cO_{\!X}^\times \to \iota_*\cO_{\!Y_\red}^\times\big)$.
We define a filtration
\[ 0 \subset  U^1\nsp R^r\nsp j_*\mu_{p^n}^{\otimes r} \subset F\nsp R^r\nsp j_*\mu_{p^n}^{\otimes r}
 \subset R^r\nsp j_*\mu_{p^n}^{\otimes r} \]
on the sheaf $R^rj_*\mu_{p^n}^{\otimes r}$ as
\begin{align*}
U^1\nsp R^r\nsp j_*\mu_{p^n}^{\otimes r} & := \hbox{the subsheaf generated \'etale locally by symbols of the form} \\
& \qquad \hbox{$\{a,b_1,\dotsc,b_{r-1}\}$ with $a \in U^1\cO_{\!X}^\times$ and $b_j \in j_*\cO_{\!X[p^{-1}]}^\times$,} \phantom{|_{\big|}} \\
F\nsp R^r\nsp j_*\mu_{p^n}^{\otimes r} & := \hbox{the subsheaf generated \'etale locally by $U^1\nsp R^r\nsp j_*\mu_{p^n}^{\otimes r}$ and the symbols} \\
& \qquad \hbox{$\{a_1,a_2,\dotsc,a_r\}$ with $a_j \in \cO_{\!X}^\times$.}
\end{align*}
We have $R^r\nsp j_*\mu_{p^n}^{\otimes r}/F\nsp R^r\nsp j_*\mu_{p^n}^{\otimes r} \cong \iota_*\nu_{Y\!,\hspace{.3pt}n}^{r-1}$ by \cite{SS} 1.1 (see also Remark \ref{rem:kato} below) and the same arguments as in \cite{Sa0} 3.4.2, and hence
\begin{equation}\label{eq1-1-1}
\cH^{r}(\fT_n(r)) \cong F\nsp R^r\nsp j_*\mu_{p^n}^{\otimes r}.
\end{equation}
The assertion follows from this fact and \cite{Sa0} 2.1.2\,(1).
\end{pf}

\begin{rem}\label{rem:kato}
{\rm
The assumption in \cite{SS} 1.1 that the base field $K$ contains a primitive $p$-th root of unity can be removed by the following argument due to K.\ Kato, \cite{KSS}.
Without loss of generality, we may assume that $\fO$ is henselian local and that $X$ is an affine scheme of the from
\[ X=\Spec(\fO[t_0,t_1,\dotsc,t_d]/(t_0^{e_0}t_1^{e_1}\dotsb t_c^{e_c}-\pi)) \]
for some integers $0 \leqq c \leqq d$ and $e_0,e_1,\dotsc,e_c \geqq 1$ and some prime element $\pi \in \fO$.
Put $\varpi:=\sqrt[p-1]{\pi\sp}$ and $\fO'':=$ the valuation ring of $K(\varpi)$.
There is a finite flat extension of $X$
\[ X''=\Spec(\fO''[T_0,T_1,\dotsc,T_d]/(T_0^{e_0}T_1^{e_1} \dotsb T_c^{e_c}-\varpi)) \]
with $T_i:=\sqrt[p-1]{\sp t_i\sp}$, which is quasi-log smooth over $\fO''$ and $K(\varpi)$ contains a primitive $p$-th root of unity.
Hence \cite{SS} 1.1 is applicable for $X''$, and we obtain the same assertion for $X$ by a standard norm argument.}
\end{rem}

\begin{prop}[cf.\ \cite{Sa0} 4.3.1]\label{prop:bock}
For any $r \in \bZ$ and $m,n \geqq 1$, there exists a canonical distinguished triangle of the following form{\rm:}
\[\xymatrix{
\fT_n(r) \ar[r]^{\ul {p^m}} & \fT_{n+m}(r) \ar[r]^-{{\mathscr R}^m} & \fT_{m}(r) \ar[r]^-{\delta_{m,n}} & \fT_n(r)[1] \quad\; \hbox{in \;\; $D^b(X_{\et})$}.
}\]
Here $\ul {p^m}$ {\rm(}resp.\ ${\mathscr R}^m${\rm)} is a unique morphism that extends the natural  inclusion $\Ln(r) \hra \vL_{n+m}(r)$ {\rm(}resp.\ the natural surjection $\vL_{n+m}(r) \twoheadrightarrow \vL_m(r)${\rm)} on $(X[p^{-1}])_\et$ and satisfies
\[  \ul {p^m} \circ {\mathscr R}^m = \text{``}{\times}\sp  p^m\text{"}  : \vL_{n+m}(r) \lra \vL_{n+m}(r) \]
The arrow $\delta_{m,n}$ is a canonical morphism which extends the Bockstein morphism $\vL_m(r) \ra \Ln(r)[1]$ in $D^b((X[p^{-1}])_\et)$ associated with the exact sequence $0 \ra \Ln(r) \to \vL_{n+m}(r) \to \vL_{m}(r) \ra 0$.
\end{prop}

\begin{pf}
On obtains the assertion by repeating the proof of \cite{Sa0} 4.3.1, using \cite{SS} 1.1 in place of \cite{Sa0} 3.3.7\,(1).
\end{pf}

\subsection{Purity and duality}
Let $Z$ be an integral closed subscheme of $Y$, and let $i_{\nsp Z} : Z \hra Y$ and $\iota_{\nsp Z} : Z \hra X$ be the natural closed immersions. Put $c:=\codim_X(Z)$.
We define the Gysin morphism for $\iota_{\nsp Z}$ as the composite
\begin{equation}\label{eq-Gys}
 \Gys_{\iota_{\nsp Z}} : \nu_{Z,n}^{r-c}[-r-c] \os{\Gys_{i_{\nsp Z}}}\lra Ri_{\nsp Z}^!\nu_{Y\!,\hspace{.3pt}n}^{r-1}[-r-1] \os{g}\lra R\iota_{\nsp Z}^!\fT_n(r) \quad \hbox{in }\;\; D^+(Z_\et,\Ln).
\end{equation}
See \eqref{eq1-1-0} for $g$, and \cite{Sa0} 2.2.1 for $\Gys_{i_{\nsp Z}}$ (see also \cite{Sa} 2.4.1).
\begin{prop}\label{thm-purity}
\begin{enumerate}
\item[{\rm (1)}]
$\Gys_{\iota_{\nsp Z}}$ induces an isomorphism
$\nu_{Z,n}^{r-c}[-r-c] \cong \tau_{\leqq r+c}R\iota_{\nsp Z}^!\fT_n(r)$ for any $r \in \bZ$.
\item[{\rm (2)}]
Assume further the condition {\rm($\star_2$)} of Lemma \ref{lem1-0}\,{\rm(}2{\rm)}.
Then the above $\Gys_{\iota_{\nsp Z}}$ is an isomorphism for any $r \geqq d$.
\end{enumerate}
\end{prop}
\begin{pf}
(1)\; We obtain the assertion by repeating the proof of \cite{Sa0} 4.4.7,
 using \cite{SS} 1.1 and 4.5 in place of \cite{Sa0} 3.3.7.
More precisely, our task is to prove that
\[ \tau_{\leqq r+c-1} Ri_{\nsp Z}^!(\tau_{\geqq r+1}\iota^*Rj_*\mu_{p^n}^{\otimes r}) = 0, \]
which is reduced, by a standard argument using \cite{SS} 1.1, to showing the semi-purity of Hagihara in our situation:
\[ R^qi_{\nsp Z}^!(\iota^*R^mj_*\mu_p^{\otimes r}) = 0 \quad \hbox{
 for any $m$ and $q$ with $q \leqq c-2$.} \]
This last vanishing is further reduced to the case that $K$ contains a primitive $p$-th root of unity by the argument in Remark \ref{rem:kato}, and then checked by the arguments in \cite{Sa0} A.2.9 and the fact that
 the sheaf $U^1\nsp R^m\nsp j_*\mu_{p}^{\otimes m}$ introduced in the proof of Proposition \ref{thm-contra} has a finite descending filtration for which each graded quotient is a free $(\cO_T)^p$-modules for some irreducible component $T$ of $Y$, see \cite{SS} 4.5 and the last display in the proof of loc.\ cit.\ 4.4.
\par
(2)\; Under the assumptions, the left arrow in \eqref{eq-Gys} is an isomorphism by \cite{Sa} 1.3.2 and 4.3.2.
The right arrow in \eqref{eq-Gys} is an isomorphism as well by the facts that $\tau_{\leqq r}Rj_*\mu_{p^n}^{\otimes r} \cong Rj_*\mu_{p^n}^{\otimes r}$ for any $r \geqq d$ (see the proof of Lemma \ref{lem1-0}\,(2)) and that $R\iota^!Rj_*=0$.
\end{pf}

\begin{cor}[cf.\ \cite{Sa0} 4.4.9]\label{cor-purity}
For any closed immersion $\iota_{\nsp Z} : Z \hra X$ of codimension $\geqq r+1$ and any $q \leqq 2r+1$,
 we have $R^q\iota_{\nsp Z} ^!\fT_n(r) = 0$.
\end{cor}
\begin{pf}
One obtains the corollary by the same arguments as in the proof loc.\ cit.\ 4.4.9, using Proposition \ref{thm-purity}\,(1) in place of loc.\ cit.\ 4.4.7.
\end{pf}
\pmn

Let $x$ and $y$ be points of $X$ such that $y \in \ol{\{x\}}$ and such that $c:=\codim_X(y)=\codim_X(x)+1$.
To proceed our preliminaries on the complex $\fT_n(r)$, we introduce the following residue diagram:
\begin{equation}\label{eq-RD}
\xymatrix{
\H^{r-c+1}(x,\Ln(r-c+1))
\ar[d]_{\Gys_{\iota_{\nsp x}}} \ar[r]^-{\partial} & \H^{r-c}(y,\Ln(r-c))
 \ar[d]^{\Gys_{\iota_{\nsp y}}} \\
\H^{r+c-1}_x(\Spec(\cO_{X,x}),\fT_n(r)) \ar[r]^-{\delta} & \H^{r+c}_y(\Spec(\cO_{X,y}),\fT_n(r)),
}\end{equation}
where the coefficient $\Ln(s)=\Ln(s)_z$ on a point $z$ denotes the \'etale complex $\logwitt z n s [-s]$
 (resp.\ the \'etale sheaf defined in \eqref{eqdef1-1}) if $\ch(z)=p$ (resp.\ $\ch(z) \ne p$).
If $\ch(z) \ne p$, then the Gysin map $\Gys_{\iota_{\nsp z}}$ for $\iota_{\nsp z} : z \hra \Spec(\cO_{X,z})$ is defined as the cup product with Gabber's cycle class $\cl_{X}(z) \in \H^{2c'}_z(\Spec(\cO_{X,z}),\mu_{p^n}^{\otimes c'})$, where $c':=\codim_X(z)$.
The arrow $\partial$ denotes the boundary map of Galois cohomology {\rm\cite{KCT}}, and $\delta$ denotes the connecting map of a localization long exact sequence of \'etale cohomology.
 
\begin{lem}\label{lem-compati}
The diagram \eqref{eq-RD} is anti-commutative.
\end{lem}
\begin{pf}
See \cite{JSS} Theorem 3.1.1 for the case $\ch(y) \ne p$.
The case $\ch(x)=\ch(y)=p$ follows from the definition of the Gysin morphism in \cite{Sa0} 2.2.1.
We check the case that $\ch(x)=0$ and $\ch(y)=p$, using the results in \cite{Sa0} as follows.
Put $Z:=\ol{\{y\}}$, the Zariski closure of $\{y\}$ in $X$.
We write RD$(X,x,y,r)$ for the diagram \eqref{eq-RD}.
Since the problem is \'etale local on $X$, we may assume that $X$ is affine and that $X$ is a closed subscheme of an affine space $\bA^N_{\fO}=:X'$. Let $\xi$ be the generic point of $X$ and put $c':=\codim_{X'}(X)$.
The diagram RD$(X',\xi,\eta,r+c')$ is anti-commutative for any generic point $\eta$ of $Y$ by \cite{Sa0} 6.1.1. Hence there exists a Gysin morphism for $i : X \hra X'$
\[ \Gys_i : \Tn(r)[-2c'] \lra Ri^!\Tn(r+c') \quad \hbox{ in } \;\; D^+(X_\et,\Ln), \]
which induces an isomorphism $\Tn(r)[-2c'] \cong \tau_{\leqq r+c'}Ri^!\Tn(r+c')$,
by the same arguments as in loc.\ cit.\ \S6.3.
Moreover, one obtains the transitivity assertion in loc.\ cit.\ 6.3.3 for the closed immersions $Z \hra X \hra X'$
 by the same arguments as in the proof of loc.\ cit.\ 6.3.3, where we have again used the fact that
 the diagram RD$(X',\xi,\eta,r+c')$ is anti-commutative for any generic point $\eta$ of $Y$.
Thus the anti-commutativity of RD$(X,x,y,r)$ follows from that of RD$(X',x,y,r+c')$ (loc.\ cit.\ 6.1.1) and
the purity in Proposition \ref{thm-purity}\,(1) for $Z \hra X$ and $Z \hra X'$.
\end{pf}
\par\medskip

The compatibility in Lemma \ref{lem-compati} plays an important role in the following results:

\begin{prop}\label{thm-trace}
\begin{enumerate}
\item[{\rm (1)}]
Let $\fO'$ be another Dedekind ring which is flat over $\fO$, and let $X'$ be a scheme which is regular and separated flat of finite type over $B'$ and satisfies {\rm(}$\star_1${\rm)} over $B'$.
Let $f : X' \to X$ be an arbitrary morphism, and let $\psi : X'[p^{-1}] \to X[p^{-1}]$ be the induced morphism.
Put $c:=\dim(X[p^{-1}])-\dim(X'[p^{-1}])$.
Then for any $n \geqq 1$ and $r \geqq 0$, there exists a unique morphism
\begin{equation*}
\tr_f : Rf_! \Tn(r-c)_{X'}[-2c] \lra \Tn(r)_X \quad \hbox{in} \;\; D^+(X_\et,\Ln)
\end{equation*}
that extends the push-forward map $\tr_\psi : R\psi_!\Ln(r-c)[-2c] \to \Ln(r)$ on $(X[p^{-1}])_\et$.
We will often write $\tr_{X'\!/\!X}$ for $\tr_f$ in what follows.
\item[{\rm (2)}]
Assume further the condition {\rm($\star_2$)} of Lemma \ref{lem1-0}\,{\rm(}2{\rm)}.
Then the adjunction morphism of $\tr_{X\!/\!B}=\tr_{\pi_{X\!/\!B}}$ is an isomorphism for any $r \geqq d${\rm:}
\begin{equation*}
 \Tn(r)_X[2(d-1)] \cong R\pi_{X\!/\!B}^! \Tn(r+1-d)_B \quad \hbox{in} \;\; D^+(X_\et,\Ln).
\end{equation*}
\end{enumerate}
\end{prop}
\begin{pf}
If $f$ is a locally closed immersion, the assertion (1) follows from Lemma \ref{lem-compati}, see \cite{Sa0} 6.3.4\,(2).
One can check (1) in the general case, using  \cite{SS} 1.1 and 4.5 and the arguments in \cite{Sa0} \S\S7.1--7.2;
in the step corresponding to loc.\ cit.\ 7.1.2, it is enough to consider locally free $(\cO_T)^p$-modules $\cF$ for each irreducible component $T$ of $Y$ in place of `locally free $(\cO_Y)^p$-modules $\cF$'
 (and the assumption on the perfectness of $k$ is unnecessary).
\par
As for the assertion (2) with $r=d$, see loc.\ cit.\ 7.3.1, where we have used the absolute purity \cite{FG} and the duality in \cite{JSS} Theorem 4.6.2. The assertion (2) in the case $r>d$ directly follows from the absolute purity, Lemma \ref{lem1-0}\,(2) and the base change isomorphism $R\pi_{X\!/\!B}^!Rj_{U*}=Rj_*R\pi_{X_U\!/\nsp U}^!$ (\cite{sga4} \XVIII.3.1.12.3), where $j_U$ denotes the open immersion $U:=B[p^{-1}] \hra B$.
\end{pf}

\begin{cor}\label{cor-trace}
Let $\beta : B' \to B$ be a flat morphism such that $B'$ is regular of dimension $\leqq 1$ and such that $X':=X \times_B B'$ satisfies {\rm ($\star_1$)} over $B'$.
Let $\alpha : X' \to X$ be the first projection. Then the following diagram commutes in $D^+(B'_\et,\Ln)$ for any $r \geqq d-1${\rm:}
\[ \xymatrix{
R\pi_{X'\!/\!B'!} \Tn(r)_{X'}[2(d-1)] \ar[rr]^-{\tr_{X'\!/\!B'}} && \Tn(r+1-d)_{B'} \\
\beta^*R\pi_{X\!/\!B!} \Tn(r)_X[2(d-1)] \ar[u]^{\alpha^*} \ar[rr]^-{\beta^*\tr_{X\!/\!B}} && \beta^*\Tn(r+1-d)_B. \ar[u]_{\beta^*}} \]
\end{cor}
\begin{pf}
The assertion follows from the uniqueness of the trace morphisms for $\Tn(r)$ and the base change property in \cite{sga4} \XVIII.2.9.
\end{pf}

\begin{cor}\label{cor-duality}
\begin{enumerate}
\item[{\rm (1)}]
Assume that $\fO$ is a strict henselian discrete valuation ring with algebraically closed residue field, and let $v$ be the closed point of $B$. Then there is a trace map
\[\xymatrix{
\tr_{X,Y} : \H^{2d}_c(X,\iota_*R\iota^!\Tn(d)) \ar[r]^-{\tr_{X\!/\!B}} &
 \H^2_v(B,\Tn(1)) & \ar[l]_-{\; \Gys_{\iota_v}}^-{\; \simeq} \Ln,} \]
where $\iota_v: v \hra B$ denotes the closed point of $B$ and the subscript $c$ means the \'etale cohomology with proper support over $B$.
Moreover, for any constructible $\Ln$-sheaf $\cF$ on $X$ and any $i \geqq 0$, the induced pairing
\[ \H^i_c(X,\cF) \times \Ext^{2d-i}_{X\!,\sp\Ln}(\cF\!,\iota_*R\iota^!\Tn(d)) \lra \Ln  \]
is a non-degenerate pairing of finite $\Ln$-modules.
\item[{\rm (2)}]
Assume that $\fO$ is an algebraic integer ring. Then there is a trace map
\[\xymatrix{
\tr_X : \H^{2d+1}_c(X,\Tn(d)) \ar[r]^-{\tr_{X\!/\!B}} &
 \H^3_c(B,\Tn(1)) \ar[r]^-{\tr_B}_-{\simeq} & \Ln,} \]
where the subscript $c$ means the \'etale cohomology with compact support {\rm(}see e.g.\ {\rm\cite{KCT}} \S3{\rm)}.
Moreover, for any constructible $\Ln$-sheaf $\cF$ on $X$ and any $i \geqq 0$, the induced pairing
\[ \H^i_c(X,\cF) \times \Ext^{2d+1-i}_{X\!,\sp\Ln}(\cF\!,\Tn(d)) \lra \Ln  \]
is a non-degenerate pairing of finite $\Ln$-modules.
\end{enumerate}
\end{cor}
\begin{pf}
(1)\;
By Proposition \ref{thm-trace}\,(2) for $r=d$ and the purity in Proposition \ref{thm-purity}\,(2), we have isomorphisms
\begin{align*}
R\iota^!\Tn(d) \cong R\iota^!R\pi_{X\!/\!B}^!\Tn(1)[-2(d-1)]
 = R\pi_{Y\!/\!v}^!R\iota_v^!\Tn(1)[-2(d-1)] \cong R\pi_{Y/v}^!\Ln[-2d].
\end{align*}
The assertion follows from this fact and the isomorphisms compatible with Yoneda pairings
\[ \H^*_c(X,\cF) \cong \H^*_c(Y,\iota^*\!\cF),\qquad
\Ext^*_{X\!,\sp\Ln}(\cF\!,\iota_*R\iota^!\Tn(d)) \cong \Ext^*_{Y\!,\sp\Ln}(\iota^*\!\cF\!,R\iota^!\Tn(d)),  \]
where we have used the proper base change theorem to obtain the left isomorphism.
See e.g.\ \cite{KSc} Chapter \II, Proposition 2.6.4 for the right isomorphism.
\par
(2)\;
The assertion follows from Proposition \ref{thm-trace}\,(2) and \cite{JSS} Proposition 2.4.1\,(3), Corollary 2.5.1.
\end{pf}

\begin{rem}\label{rem-PF}
{\rm
The push-forward morphism $\tr_f$ in Proposition \ref{thm-trace}\,(1) satisfies the projection formula in \cite{Sa0} 7.2.4, by the same arguments as in loc.\ cit. \sp
See also the proof of Proposition \ref{thm-trace}\,(1) as to how we modified loc.\ cit.\ 7.1.2 in our situation.
}
\end{rem}

\subsection{Cycle class morphism}\label{sect1-2}
To construct a cycle class morphism from Bloch's cycle complex (see \eqref{eq-cycle} below),
 we formulate a version of $\fT_n(r)$ with log poles and a purity for this coefficient; see also \cite{Z} for a construction assuming Gersten's conjecture for Bloch's cycle complex.
Let $D$ be a reduced normal crossing divisor on $X$ which is flat over $B$ and such that $D \cup Y_\red$ also has simple normal crossings on $X$ and such that the pair $(X,D)$ is quasi-log smooth over $B$ in the sense of \cite{SS} 5.2. We define $\fT_n(r)_{(X,D)}$ by the following distinguished triangle analogous to \eqref{eq1-1-0}:
\begin{equation}\label{eq1-1+0}
 \iota_*\nu_{(Y,E),n}^{r-1}[-r-1] \os{g}\lra \Tn(r)_{(X,D)} \os{t}\lra \tau_{\leqq r}R\psi_*\mu_{p^n}^{\otimes r} \os{(\star)}\lra \iota_*\nu_{(Y,E),n}^{r-1}[-r],
\end{equation}
where we put $E:=Y_\red \cap D$
and $\nu_{(Y,E),n}^{r-1}:=\phi_*\nu_{Y \ssm E}^{r-1}$ with $\phi : Y \ssm E \hra Y$; $\psi$ denotes the open immersion $X \ssm (Y \cup D) \hra X$. See also \cite{Sa1} 3.5 and 3.6. 
When $D=\emptyset$, we have $\Tn(r)_{(X,\emptyset)}=\Tn(r)_X$.
The following propositions concerning the complex $\fT_n(r)_{(X,D)}$ play fundamental roles in our construction of cycle class maps.

\begin{prop}[cf.\ \cite{Sa1} 6.5]\label{prop-purity}
Let $Z$ be a closed subset of $X$ of codimension $\geqq c$. Then we have
\[ \H^q_Z(X,\Tn(r)_{(X,D)})\cong
\begin{cases}
0 & \hbox{{\rm(}$q<r+c${\rm)}} \\
\H^{r+c}_{Z \ssm D}(X \ssm D,\Tn(r)) \quad & \hbox{{\rm(}$q=r+c${\rm)}.}
\end{cases} \]
In particular, if $Z$ has pure codimension $c$ on $X$, then we have
\[\H^q_Z(X,\Tn(c)_{(X,D)})\cong
\begin{cases}
0 & \hbox{{\rm(}$q < 2c${\rm)}} \\
\Ln[Z^0 \ssm D] \quad & \hbox{{\rm(}$q=2c${\rm)}},
\end{cases} \]
where $\Ln[Z^0 \ssm D]$ means the free $\Ln$-module generated over the set $Z^0 \ssm D$.
\end{prop}
\begin{pf}
One obtains the assertion by repeating the arguments in the proof of loc.\ cit.\ 6.5,
 using \cite{SS} 1.1 and 4.5 (resp.\ Corollary \ref{cor-purity} of the previous subsection) in place of \cite{Sa1} 3.3 (resp.\ \cite{Sa0} 4.4.9).
We do not need to assume the existence of primitive $p$-th roots of unity in $K$ by the argument in Remark \ref{rem:kato}.
\end{pf}

\begin{prop}[cf.\ \cite{Sa1} 4.3]\label{lem-homotopy}
Let $E \to X$ be a vector bundle of rank $a$, and let $f : \bP:=\bP(E \oplus 1) \to X$ be its projective completion.
Let $\bP':=\bP(E)$ the projective bundle associated with $E$, regarded as the infinite hyperplane section of $\bP$.
Then the composite morphism
\[ \Tn(r)_X \lra Rf_*\Tn(r)_{\bP} \lra Rf_*\Tn(r)_{(\bP,\bP')} \]
is an isomorphism in $D^+(X_\et,\Ln)$.
\end{prop}
\begin{pf}
One can extend the Dold-Thom isomorphism (loc.\ cit.\ 4.1) and the distinguished triangle in loc.\ cit.\ 3.12 to the situation of this section, by repeating the same arguments as in the proofs of loc.\ cit.\ 4.1 and 3.12, using \cite{SS} 1.1 and 4.5 (note also Remark \ref{rem:kato} of this section).
The assertion follows from those facts and Remark \ref{rem-PF}.
\end{pf}
\par\medskip

Let $\Et/X$ be the underlying category of $X$-schemes of the \'etale site $X_\et$.
For a scheme $U$ and $r \geqq 0$, let $z^r(U,*)$ be Bloch's cycle complex \cite{B2}.
We define a complex $\bZ(r)$ of presheaves on $\Et/X$ by the assignment
\[ \bZ(r) : U \in \Ob(\Et/X) \; \longmapsto \; z^r(U,*)[-2r], \]
which is in fact a complex of sheaves in the Zariski and the \'etale topologies.
We call $\bZ(r)$ the {\it motivic complex} of $X$ of weight $r$. For a closed subset $C \subset X$ and $U \in \Ob(\Et/X)$, put $C_U:=C \times_X U$ and 
let $z^r_{C_U}(U,q)$ be the subgroup of $z^r(U,q)$ consisting of the cycles on $U \times \vD^q$ of codimension $r$ whose support is contained in $C_U \times \vD^q$ (and which satisfies the face condition).
The collection $\{z^r_{C_U}(U,q)\}_{q \geqq 0}$ forms a subcomplex of $z^r(U,*)$, and
we define a subcomplex $\bZ(r)_{C \subset X} \subset \bZ(r)$ by the assignment
\[ \bZ(r)_{C \subset X} : U \in \Ob(\Et/X) \; \longmapsto \; z^r_{C_U}(U,*)[-2r]. \]
By Propositions \ref{prop-purity} and \ref{lem-homotopy}, Lemma \ref{lem-compati} and the same arguments as in \cite{Sa1} \S7 (see also Remark \ref{rem-Sa} below), one obtains a {\it cycle class morphism}
\begin{equation}\label{eq-cycle}
 \cl_{C \subset X,\Ln} :
 \bZ(r)_{C \subset X} \otimes \Ln \lra R\ul{\vG}_C(X,\fT_n(r)) \quad \hbox{ in } D(X_\et,\Ln)
\end{equation}
for any $r \geqq 0$, which yields the {\it cycle class map} on hypercohomology groups
\[ \cl_{C \subset X,\Ln} : \H^*_C(X_\Zar,\bZ(r)\otimes \Ln) \lra \H^*_C(X,\fT_n(r)). \]
When $C=X$, the group on the left hand side will be denoted by $\H^*_{\cM}(X,\Ln(r))$, and the map $\cl_{X,\Ln}:=\cl_{X \subset X,\Ln}$ will be computed in Lemma \ref{lem6-2}\,(3) below under the assumption that $d=2$.

\begin{rem}\label{rem-Sa}
{\rm
To follow the arguments in \cite{Sa1} \S7, we have used the projection formula in \cite{Sa0} Corollary 7.2.4,
 which has been extended to our situation in Remark \ref{rem-PF}.
We also need to extend the compatibility fact in \cite{Sa0} Corollary 6.3.3 to our situation, where the push-forward morphism in Proposition \ref{thm-trace}\,(1) plays the role of $\Gys_i$ of loc.\ cit.\ 6.3.3.
One can easily check the details by Lemma \ref{lem-compati} and the proof of loc.\ cit.\ 6.3.3.
}
\end{rem}

Let $C_2 \subset C_1$ be closed subsets of $X$, and let $\phi : X':=X \ssm C_2 \hra X$ be the natural open immersion.
Put $C':=C_1 \ssm C_2$. Then the squares in $D(X_\et,\Ln)$
\begin{equation}\label{eq-cd-cycle}
\xymatrix{
\bZ(r)_{C_2 \subset X} \otimes \Ln \ar[r] \ar[d]_{\cl_{C_2 \subset X,\Ln}}
 & \bZ(r)_{C_1 \subset X} \otimes \Ln \ar[r]^-{\phi^\sharp} \ar[d]_{\cl_{C_1 \subset X,\Ln}}
  & R\phi_*\bZ(r)_{C' \subset X'}\otimes \Ln \ar[d]_{\cl_{C' \subset X',\Ln}} \\
R\ul{\vG}_{C_2}(X,\fT_n(r)) \ar[r] & R\ul{\vG}_{C_1}(X,\fT_n(r)) \ar[r] &
 R\phi_* R\ul{\vG}_{C'}(X',\fT_n(r))
}\end{equation}
are commutative by the construction of cycle class morphisms.
From this commutative diagram, one obtains another commutative diagram in $D(X_\et,\Ln)$
\begin{equation}\label{eq-cd-cycle2}
\xymatrix{
R\phi_*\bZ(r)_{C' \subset X'}\otimes \Ln \ar[r]^-{\delta} \ar[d]_{\cl_{C' \subset X',\Ln}} &
R\ul{\vG}_{C_2}(X,\bZ(r)_{C_1 \subset X}\otimes \Ln)[1] \ar[d]_{\cl_{C_2 \subset X,\Ln}} &
 \ar[l]_-\gamma \bZ(r)_{C_2 \subset X} \otimes \Ln [1] \ar[ld]^{\; \cl_{C_2 \subset X,\Ln}} \\
R\phi_* R\ul{\vG}_{C'}(X',\fT_n(r)) \ar[r]^-{\delta} &
R\ul{\vG}_{C_2}(X,\fT_n(r))[1],
}\end{equation}
where the arrows $\delta$ are the connecting morphisms of localization triangles (see \cite{Sa0} 1.9).
\begin{rem}
{\rm
The arrow $\gamma$ of \eqref{eq-cd-cycle2} is {\it not} an isomorphism, or equivalently,
the upper row of \eqref{eq-cd-cycle} does {\it not} fit into any distinguished triangle in $D(X_\et,\Ln)$.
If one considers localization triangles in the Zariski topology, then the morphism corresponding to $\gamma$ of \eqref{eq-cd-cycle2} is an isomorphism by Levine \cite{Le} Theorem 1.7.
}
\end{rem}

\subsection{Cospecialization and a residue diagram}\label{sect1-3}
In this subsection, we consider a residue map and prove its contravariance, which will be useful in \S\ref{sect2-3} below.
We suppose that $\pi_{X\!/\!B} : X \to B$ is {\it proper}, and that $\fO$ is a {\it henselian discrete valuation ring with algebraically closed residue field}.
Put $I_K:=\Gal(\ol K/K)$. By the duality in Corollary \ref{cor-duality}\,(1) and the Poincar\'e duality for $X_{\ol K}$, the cospecialization map
\[ \cosp_X : \H^i(Y,\Ln) \cong \H^i(X,\Ln) \lra \H^i(X_{\ol K},\Ln)^{I_K} \]
induces a canonical homomorphism
\[ \res_X : \H^{i'}(X_{\ol K},\mu_{p^n}^{\otimes d-1})_{I_K} \lra \H^{i'+2}_Y(X,\Tn(d)), \]
where we put $i':=2(d-1)-i$.
\begin{prop}\label{prop-res}
For any $i \geqq 0$, the following diagram is anti-commutative{\rm:}
\[\xymatrix{
\H^1(I_K,\H^i(X_{\ol K},\mu_{p^n}^{\otimes d})) \ar[r]^-{\alpha} \ar[d]_\vare &
\H^i(X_{\ol K},\mu_{p^n}^{\otimes d-1})_{I_K} \ar[d]^{\res_X} \\
\H^{i+1}(X_{\nsp K},\mu_{p^n}^{\otimes d}) \ar[r]^{\delta_X} & \H^{i+2}_Y(X,\Tn(d)),
}\]
where the left vertical arrow is an edge map of a Hochschild-Serre spectral sequence, and the upper horizontal arrow denotes the composite map
\begin{align*}
\H^1(I_K,\H^i(X_{\ol K},\mu_{p^n}^{\otimes d}))
 &\lra \H^1(I_K,\mu_{p^n}) \otimes \H^i(X_{\ol K},\mu_{p^n}^{\otimes d-1})_{I_K} \\
 & \; \cong \; K^\times \otimes \H^i(X_{\ol K},\mu_{p^n}^{\otimes d-1})_{I_K} \os{\ord_K}\lra \bZ \otimes \H^i(X_{\ol K},\mu_{p^n}^{\otimes d-1})_{I_K}.
\end{align*}
The bottom horizontal arrow is the connecting map of a localization long exact sequence.
\end{prop}
\begin{pf}
The following diagram of trace maps and boundary maps are commutative:
\[
\xymatrix{
\H^{2d-1}(X_{\nsp K},\mu_{p^n}^{\otimes d}) \ar[r]^-{\tr_{X\!/\!B}} \ar[d]_{\delta_X} &
\H^1(I_K,\mu_{p^n}) \ar@{=}[r]^-\sim \ar[d]_{\delta_B} & K^\times \otimes \Ln \ar[d]^{-\ord_K} \\
\H^{2d}_Y(X,\Tn(d)) \ar[r]^-{\tr_{X\!/\!B}} & \H^2_v(B,\Tn(1)) \ar@{=}[r]^-{\tr_{B,v}}_-{\sim}& \Ln ,
}\]
where $v$ denotes the closed point of $B$, and $\tr_{B,v}$ means $\tr_{X,Y}$ for $(X,Y)=(B,v)$ (see Corollary \ref{cor-duality}\,(1)).
See Lemma \ref{lem-compati} for the commutativity of the right square.
The assertion follows from this commutativity and the following obvious commutative square:
\[\xymatrix{
\H^0(I_K,\H^{i'}(X_{\ol K},\Ln)) & \ar@{=}[l] \H^{i'}(X_{\ol K},\Ln)^{I_K} \\
\ar[u] \H^{i'}(X_{\nsp K},\Ln) & \ar[l]\ar[u] \H^{i'}(X,\Ln).
}\]
The details are straight-forward and left to the reader.
\end{pf}
\par\medskip

The following consequence of Proposition \ref{prop-res} will be useful later.
Let $\fO'$ be another strict henselian local ring which is flat over $\fO$ and whose residue field $k'$ is algebraically closed. Let $X'$ be a scheme which is regular, proper flat of finite type over $B'$ and satisfies ($\star_1$) over $B'$.
Put $L:=\Frac(\fO')$, $I_L:=\Gal(\ol L/L)$ and $Y':=X' \otimes_{\fO'} k'$.
Assume that \[ \dim(X')=\dim(X)=d \]
(hence that $\dim(X'_L)=\dim(X_K)=d-1$).
Under this setting we obtain:

\begin{cor}\label{cor-residue}
For any morphism $f : X' \to X$ and any $i \geqq 0$, the diagram
\[\xymatrix{
\H^{i}(X'_{\ol L},\mu_{p^n}^{\otimes d-1})_{I_L} \ar[d]_{\res_{X'}} & \ar[l]_-{\;\;f^\sharp}
\H^{i}(X_{\ol K},\mu_{p^n}^{\otimes d-1})_{I_K} \ar[d]^{\res_X} \\
 \H^{i+2}_{Y'}(X',\Tn(d))  & \ar[l]_-{\;\;f^\sharp} \H^{i+2}_Y(X,\Tn(d))
}\]
is commutative, that is, the map $\res_X$ is contravariant in $X$.
\end{cor}

\begin{pf}
In the diagram of Proposition \ref{prop-res},
the composite map $\delta \circ \ep$ is contravariant in $X$ by Proposition \ref{thm-contra}, and the map $\alpha$ is surjective by the fact that $\cd(I_K)=1$. The corollary follows from these facts and Proposition \ref{prop-res}.
\end{pf}

\section{A filtration on the direct image}\label{etale}
Let $\pi_{X\!/\!B} : X \to B=\Spec(\fO)$ be as in the beginning of \S\ref{sect1}.
In this section, we assume
\begin{itemize}
\item[($\star_2$)]
 {\it any residue field of $\fO$ of characteristic $p$ is perfect.}
\end{itemize}
Under this assumption, we introduce objects $\fH^*(X,\Tn(r))$ of $D^+(B_{\et},\Ln)$ for $r \geqq d=\dim(X)$, which play central roles throughout this paper.
The \'etale cohomology of $B$ with coefficients in these new objects will be related to the \'etale cohomology of $X$ with coefficients in $\Tn(r)$ by the spectral sequence \eqref{ss:leray} below.

\subsection{\'Etale complex $\fH^m(X,\Tn(r))$}\label{sect2.1}
\begin{lem}\label{lem-direct}
For any $r \geqq d$, we have 
\begin{equation}\label{selmer:q-isom}
R \pi_{X\!/\!B*}\Tn(r) \cong R\sHom_{B\nsp,\sp\Ln}(R \pi_{X\!/\!B!}\Tn(d-r),\Tn(1))[2-2d]
\end{equation}
in $D^+(B_{\et},\Ln)$.
\end{lem}
\begin{pf}
Since $r \geqq d$ by assumption, there exists a canonical isomorphism
\begin{equation}\label{eq-RHom}
 \Tn(r) \cong R\sHom_{X\!,\sp\Ln}(\Tn(d-r),\Tn(d)) \quad \hbox{ in } \;\; D^+(X_{\et},\Ln),
\end{equation}
which is obvious if $r=d$, and otherwise a consequence of Lemma \ref{lem1-0}\,(2) and the adjunction in \cite{sga4}\,\XVIII.3.1.10 for the open immersion $X[p^{-1}] \hra X$.
Hence we have
\begin{align*}
 R \pi_{X\!/\!B*}\Tn(r)
   & \cong R \pi_{X\!/\!B*} R\sHom_{X\!,\sp\Ln}(\Tn(d-r),R \pi_{X\!/\!B}^!\Tn(1)[2-2d])
   \phantom{\big|_|}\\
   & \cong R\sHom_{B\nsp,\sp\Ln}(R \pi_{X\!/\!B!}\Tn(d-r),\Tn(1))[2-2d]
\end{align*}
in $D^+(B_{\et},\Ln)$,
by Proposition \ref{thm-trace}\,(2) and \cite{sga4}\,\XVIII.3.1.10 for $\pi_{X\!/\!B}$.
\end{pf}

\begin{defn}\label{def:mot}
{\rm
For each $m \in \bZ$, we define
\begin{align*}
\fH^{\leqq m}(X,\Tn(r)) & := R\sHom_{B\nsp,\sp\Ln}(\tau_{\geqq 2(d-1)-m}R \pi_{X\!/\!B!}\Tn(d-r)_X,\Tn(1)_B)[2-2d],
 \phantom{\big|_|}\\
\fH^m(X,\Tn(r)) & := R\sHom_{B\nsp,\sp\Ln}(R^{2(d-1)-m}\pi_{X\!/\!B!}\Tn(d-r)_X,\Tn(1)_B),
\end{align*}
which are objects of $D^+(B_{\et},\Ln)$.
}
\end{defn}
\begin{caution}
{\rm
$\fH^m(X,\Tn(r))$ is NOT the sheaf $R^m\pi_{X\!/\!B*}\Tn(r)$, but a complex of sheaves.
}
\end{caution}

By Lemma \ref{lem-direct} and the proper base change theorem (for $R \pi_{X\!/\!B!}$), we have
\begin{align}
\fH^{\leqq m}(X,\Tn(r)) & \cong
\begin{cases}
    0    &  \hbox{{\rm ($m \leqq -1$)}}\phantom{\big|_{|}} \\
   R \pi_{X\!/\!B*}\Tn(r)_X  \;\;\; &  \hbox{{\rm ($m \geqq 2(d-1)$)}}
\end{cases}
 \phantom{\bigg|_{\Big|}} \label{filter} \\
\fH^m(X,\Tn(r)) &= \; 0 \qquad \hbox{unless \; $0 \leqq m \leqq 2(d-1)$}.  \label{filter2}
\end{align}
For any $m \in \bZ$, we have a natural distinguished triangle of the form
\begin{multline}\label{eq2-1}
\fH^{\leqq m-1}(X,\Tn(r)) \lra \fH^{\leqq m}(X,\Tn(r)) \lra \fH^m(X,\Tn(r))[-m] \phantom{\big|_{|}} \\
   \lra \fH^{\leqq m-1}(X,\Tn(r))[1]. 
\end{multline}
The data $\{ \fH^{\leqq m}(X,\Tn(r)) \}_{m \leqq 2(d-1)}$ form a finite ascending filtration on
 $\fH^{\leqq 2(d-1)}(X,\Tn(r))$ $\cong R \pi_{X\!/\!B*}\Tn(r)_X$, and yield a convergent spectral sequence
\begin{equation}\label{ss:leray}
E^{a,b}_2=\H^a(B,\fH^b(X,\Tn(r))) \Lra \H^{a+b}(X,\Tn(r)).
\end{equation}
\par
To illustrate our complex $\fH^m(X,\Tn(r))$,
 we show here the following proposition assuming that $\pi_{X\!/\!B}$ is proper.
See Proposition \ref{prop:local-str} below for more detailed computations without the properness assumption.
\begin{prop}\label{ex:trace}
Assume that $\pi_{X\!/\!B} : X \to B$ is proper and that $r \geqq d$.
\begin{enumerate}
\item[{\rm(1)}]
Let $U \subset B[p^{-1}]$ be an open subset for which $\pi_{X_U/U} : X_U=X \times_B U \to U$ is
 smooth {\rm(}and proper{\rm)}.
Then $\fH^{m}(X,\Tn(r))|_U$ is the locally constant constructible sheaf placed in degree $0$,
 associated with $\H^m(X_{\ol K},\mu_{p^n}^{\otimes r})$.
\item[{\rm(2)}]
Assume further that the generic fiber $X_{\nsp K}$ is geometrically connected over $K$.
Then the trace map $\tr_{X\!/\!B} : R\pi_{X\!/\!B*}\Tn(r)_X[2(d-1)] \to \Tn(r+1-d)_B$ induces an isomorphism
\begin{equation}\label{isom:trace}
\fH^{2(d-1)}(X,\Tn(r)) \cong \Tn(r+1-d)_B.
\end{equation}
\end{enumerate}
\end{prop}
To prove this proposition, we need the following lemma:
\begin{lem}\label{lem1-1}
Let $Z$ be a scheme and let $\cF$ be a locally constant constructible $\Ln$-sheaf on $Z_\et$.
Then we have
\[ \sHom(\cF,\Ln)_{\ol x} \cong \Hom(\cF_{\ol x},\Ln) \quad \hbox{and} \quad
 \sExt^q_{Z,\sp\Ln}(\cF,\Ln) = 0 \quad \hbox{{\rm(}$q \geqq 1${\rm)}}. \]
\end{lem}
\begin{pf*}{\it Proof of Lemma \ref{lem1-1}}
Since $\cF$ is a pseudo-coherent $\Ln$-module on $Z_\et$ in the sense of \cite{milne}\,p.\ 80, we have
\[ \sExt^q_{Z\nsp,\sp\Ln}(\cF,\Ln)_{\ol x} \cong \Ext^q_{\Ln}(\cF_{\ol x},\Ln) \]
for any $q \geqq 0$ by loc.\ cit.\ \II.3.20. The assertions follow from this fact and the fact that $\Ln$ is an injective $\Ln$-module.
\end{pf*}
\par\medskip
\begin{pf*}{\it Proof of Proposition \ref{ex:trace}}
(1)\;
By definition, we have
\[\fH^{m}(X,\Tn(r))|_U = R\sHom_{U\!,\sp\Ln}(R^{2(d-1)-m}\pi_{X_U\!/U*}\Ln(d-r),\Ln(1)). \]
Since $R^{2(d-1)-m}\pi_{X_U/U*}\Ln(d-r)$ is locally constant and constructible by the proper smooth base change theorem,
the object on the right hand side is isomorphic to the sheaf
\[ \sHom_{U\!,\sp\Ln}(R^{2(d-1)-m}\pi_{X_U\!/\nsp U*}\Ln(d-r),\Ln(1)) \]
placed in degree $0$, by Lemma \ref{lem1-1}. Then the assertion follows from the Poincar\'e duality.
\par
(2)\;
We have
\[ \fH^{2(d-1)}(X,\Tn(r)) = R\sHom_{B\nsp,\sp\Ln}(\pi_{X\!/\!B*}\Tn(d-r)_X,\Tn(1)_B) \phantom{\big|_{\big|}} \]
by definition, and $\pi_{X\!/\!B*}\Tn(d-r)_X \cong \Tn(d-r)_B$ for $r \geqq d$ by the connectedness of the geometric fibers.
The assertion follows from this fact and Lemma \ref{lem1-0}\,(2) for $B$.
\end{pf*}

\subsection{Local computations}
We investigate here the local structure of $\fH^m(X,\Tn(r))$ around the closed points on $B$ {\it without assuming} that $\pi_{X\!/\!B}$ is proper.
For a closed point $v \in B$, we often write $Y_v$ (resp.\ $Y_{\ol v}$, $X_{\ol v}$) for $X \times_B v$ (resp.\ $X \times_B \ol v$, $X \times_B B_{\ol v}^\sh$), where $B_{\ol v}^\sh$ denotes the spectrum of the strict henselization of $\fO_v=\cO_{B,v}$ at its maximal ideal.
\begin{prop}\label{prop:local-str}
Let $v$ be a closed point on $B$, and let $q$ and $m$ be integers.
We write $\iota_v$ for the closed immersion $v \hra B$ and $j_v$ for the open immersion $B \ssm v \hra B$.
Assume $r \geqq d$. Then
\begin{enumerate}
\item[{\rm(1)}]
We have $R^q \iota_v^!\fH^m(X,\Tn(r))=0$ unless $q =2$, and a canonical isomorphism
\[ (R^2\iota_v^!\fH^m(X,\Tn(r)))_{\ol v} \cong \H^{m+2}_{Y_{\ol v}}(X_{\ol v}, \Tn(r)). \]
Moreover, we have $R \iota_v^!\fH^m(X,\Tn(r))=0$, if $\ch(v)=p$ and $r > d$.
\item[{\rm(2)}]
We have
\[ (R^qj_{v*}j_v^*\fH^m(X,\Tn(r)))_{\ol v} \cong \H^q(I_v,\H^m(X_{\ol K},\mu_{p^n}^{\otimes r})), \]
where $I_v$ denotes the inertia subgroup of $G_{\!K}$ at $v$. Consequently, we have
\[ R^qj_{v*}j_v^*\fH^m(X,\Tn(r))=0 \]
unless $q=0$ or $1$, by the fact that $\cd_p(I_v)=1$ {\rm(}see {\rm\cite{se}}\,Chapter \II, \S3.3{\rm)}.
\item[{\rm(3)}]
We have
\[ \cH^q(\fH^m(X,\Tn(r)))_{\ol v} \cong
      \begin{cases}
       \H^m(X_{\ol K},\mu_{p^n}^{\otimes r})^{I_v}
       &  \hbox{ if } q=0\\
       0 &  \hbox{ if } q \not= 0,1\hbox{ or } 2
      \end{cases} \]
and an exact sequence
\begin{align*}
0 \lra \cH^1(\fH^m(X,\Tn(r)))_{\ol v}
 &\lra \H^1(I_v,\H^m(X_{\ol K},\mu_{p^n}^{\otimes r}))\\
 \os{\delta^+}{\lra}& \; \H^{m+2}_{Y_{\ol v}}(X_{\ol v}, \Tn(r))
    \lra \cH^2(\fH^m(X,\Tn(r)))_{\ol v} \ra 0.
\end{align*}
Here $\cH^q(-)$ denotes the $q$-th cohomology sheaf, and $\delta^+$ denotes the composite map $\delta_X \circ \vare$ in the diagram of Proposition \ref{prop-res}.
\end{enumerate}
\end{prop}

\begin{pf*}{Proof of Proposition \ref{prop:local-str}}\!
(1)\;
By the definition of $\fH^m(X,\Tn(r))$ in Definition \ref{def:mot} and the adjunction in \cite{sga4}\,\XVIII.3.1.12.2, we have
\begin{align}
\notag
R\iota_v^!\fH^m(X,\Tn(r)) &= R\iota_v^!R\sHom_{B,\Ln}(R^{2(d-1)-m}\pi_{X\!/\!B!}\Tn(d-r)_X,\Tn(1)_B)
 \phantom{\big|_{|}} \\
 \label{eq2-0}
 & \cong R\sHom_{v,\Ln}(\iota_v^*R^{2(d-1)-m}\pi_{X\!/\!B!}\Tn(d-r)_X,R\iota_v^!\Tn(1)_B)
\phantom{\big|_{|}} \\
 & \cong R\sHom_{v,\Ln}(R^{2(d-1)-m}\pi_{Y_v\nsp/\nsp v!}(\iota_{Y_v}^*\Tn(d-r)_X),\Ln)[-2],
\notag
\end{align}
where $\iota_{Y_v}$ denotes the closed immersion $Y_v \ra X$, and we have used the proper base change theorem for $R\pi_{X\!/\!B!}$ and the purity in Proposition \ref{thm-purity}\,(2) for $\Tn(1)_B$ in the last isomorphism.
In particular if $\ch(v)=p$ and $r>d$, then $\iota_{Y_v}^*\Tn(d-r)_X$ is zero by definition and we have $R\iota_v^!\fH^m(X,\Tn(r))=0$, which shows the third assertion of (1).
If $\ch(v) \ne p$ or $r=d$, then $R\iota_v^!\fH^m(X,\Tn(r))$ is acyclic outside of degree $2$ by \eqref{eq2-0} and Lemma \ref{lem1-1} for $Z=v$.
Moreover, if $r=d$, then we have
\begin{align*}
R^2\iota_v^!\fH^m(X,\Tn(d))& \!\os{\eqref{eq2-0}}\cong \!\! \sHom_{v,\Ln}(R^{2(d-1)-m}\pi_{Y_v\nsp/\nsp v!}\Ln,\Ln) \phantom{\big|_{|}} \\
& \;\sp \cong \cH^{m-2(d-1)}(R\sHom_{v,\Ln}(R\pi_{Y_v\nsp/\nsp v!}\Ln,\Ln))
\phantom{\big|_{|}} \\
& \;\sp \cong \cH^{m-2(d-1)}(R\pi_{Y_v\nsp/\nsp v*}R\pi_{Y_v\nsp/\nsp v}^!\Ln)
\end{align*}
again by Lemma \ref{lem1-1} for $Z=v$ and adjunction, and we have
\begin{equation}\label{eq2-3}
 R\pi_{Y_v\nsp/\nsp v}^!\Ln \cong R\pi_{Y_v\nsp/\nsp v}^!R\iota_v^!\Tn(1)_B[2] \cong R\iota_{Y_v}^!\Tn(d)_X[2d]
\end{equation}
by the purity in Proposition \ref{thm-purity}\,(2) for $v \hra B$ and Proposition \ref{thm-trace}\,(2). Hence we have
\[ (R^2\iota_v^!\fH^m(X,\Tn(d)))_{\ol v} \cong \H^{m+2}_{Y_{\ol v}}(X_{\ol v}, \Tn(d)). \]
The isomorphism in the case that $r>d$ and $\ch(v) \ne p$ is similar and left to the reader.
\par
(2)\;
We may assume that $B$ is local with closed point $v$, without loss of generality.
Put $\eta:=B \ssm v$, which is the generic point of $B$.
The sheaf $j_v^*R^{2(d-1)-m}\pi_{X\!/\!B!}\Tn(d-r)$ is locally constant on $\eta_\et$, and the object
\[ j_v^*\fH^m(X,\Tn(r)) = R\sHom_{\eta,\Ln}(j_v^*R^{2(d-1)-m}\pi_{X\!/\!B!}\Tn(d-r),\mu_{p^n}) \]
is isomorphic to the sheaf (on $\eta_\et$) associated with $\H^m(X_{\ol K},\mu_{p^n}^{\otimes r})$ placed in degree $0$ by Lemma \ref{lem1-1} for $Z=\eta$ and the Poincar\'e duality.
The assertion follows from this fact.
\par
(3)\; The assertion follows from Proposition \ref{prop:local-str}\,(1), (2) and the fact that the stalk at $\ol v$ of the connecting homomorphism
\[ \delta_{B,B\ssm v} : R^1j_{v*}j_v^* \fH^m(X,\Tn(r)) \lra \iota_{v*}R^2\iota_v^!\fH^m(X,\Tn(r)) \]
agrees with $\delta^+$ up to a sign.
\end{pf*}
\par\smallskip

The following corollary follows from Proposition \ref{prop:local-str}\,(1) and (3).
\begin{cor}\label{cor2-1}
\begin{enumerate}
\item[{\rm(1)}]
If $\ch(v)=p$ and $r > d$, then $\fH^m(X,\Tn(r)) \cong Rj_{v*}j_v^*\fH^m(X,\Tn(r))$.
\item[{\rm(2)}]
$\fH^m(X,\Tn(r))$ is concentrated in $[0,2]$, and $R \pi_{X\!/\!B*}\Tn(r)$ is concentrated in $[0,2d]$.
\end{enumerate}
\end{cor}

\subsection{Rigidity}
In this subsection, we assume further that $\fO$ is henselian local with finite residue field. Let $\fO'$ be the completion of $\fO$ at its maximal ideal, and put
\[ B':=\Spec(\fO)\quad \hbox{and} \quad X':=X \times_B B'. \]
Let $v$ be the closed point of $B'$, which we identify with the closed point of $B$. Let $Y'$ be the special fiber of $\pi_{X'\!/\!B'}: X' \to B'$, and let $Y$ be the special fiber of $\pi_{X\!/\!B}: X \to B$.
We have cartesian squares
\begin{equation}\label{cartesian}
\xymatrix{
Y' \,\ar@<-1pt>@{^{(}->}[r] \ar@{=}[d] \ar@{}[rd]|{\square}
& X' \ar[d]_{\alpha} \ar[r]^{\pi_{X'\!/\!B'}} \ar@{}[rd]|{\square}
 & B' \ar[d]_{\beta} \ar@{}[rd]|{\square} & \ar@<1pt>@{_{(}->}[l]_-{\iota_v} \, v \ar@<1pt>@{=}[d] \\
Y \, \ar@<-1pt>@{^{(}->}[r]&
 X \ar[r]^-{\pi_{X\!/\!B}} & B & \ar@<1pt>@{_{(}->}[l]_-{i_v} \, v.\!}
\end{equation}
We prove here the following preliminary result, where we do {\it not} assume that $\pi_{X\!/\!B}$ is proper:
\begin{prop}[rigidity]\label{thm:rigidity}
For any $r \geqq d$, there exist canonical isomorphisms
\begin{align*}
&\psi_1 :  R\pi_{X\!/\!B*}\Tn(r)_X \os{\simeq}\lra R\beta_*R\pi_{X'\!/\!B'*}\Tn(r)_{X'} \phantom{|_{\big|}} \\
&\psi_2^m : \fH^{\leqq m}(X,\Tn(r)) \os{\simeq}\lra R\beta_*\fH^{\leqq m}(X',\Tn(r))
 \qquad\quad \hspace{-.4pt} \hbox{{\rm(}\sp $^\forall m \in \bZ${\rm)}} \phantom{|_{\big|}} \\
&\psi_3^m : \fH^m(X,\Tn(r)) \os{\simeq}\lra R\beta_*\fH^m(X',\Tn(r))
 \qquad\qquad \hbox{{\rm(}\sp $^\forall m \in \bZ${\rm)}} \phantom{|_{\big|}} \\
&\psi_4^m : i_{v*}Ri_v^!\fH^m(X,\Tn(r)) \os{\simeq}\lra i_{v*}R\iota_v^!\fH^m(X',\Tn(r))
 \qquad\;\; \hspace{1.7pt} \hbox{{\rm(}\sp $^\forall m \in \bZ${\rm)}}
\end{align*}
in $D^b(B_{\et},\Ln)$, where $i_v : v \hra B$ and $\iota_v : v \hra B'$ are canonical closed immersions.
\end{prop}
\begin{cor}\label{cor:rigidity}
We have canonical isomorphisms for any $q, m \in \bZ$ and any $r \geqq d$
{\allowdisplaybreaks
\begin{align*}
 & \H^q(X,\Tn(r)_X) \cong \H^q(X',\Tn(r)_{X'}),\phantom{|_{\big|}}\\
 & \H^q_Y(X,\Tn(r)_X) \cong \H^q_{Y'}(X',\Tn(r)_{X'}),\phantom{|_{\big|}}\\
 & \H^q(B,\fH^m(X,\Tn(r))) \cong \H^q(B',\fH^m(X',\Tn(r))),\phantom{|_{\big|}} \\
 & \H^q_v(B,\fH^m(X,\Tn(r))) \cong \H^q_v(B',\fH^m(X',\Tn(r))).
\end{align*}
}\end{cor}
\begin{pf*}{\it Proof of Proposition \ref{thm:rigidity}}
Let $\res_X$ and $\res_B$ be the pull-back morphisms
\[ \res_X : \alpha^*\Tn(r)_X \lra \Tn(r)_{X'} \quad \hbox{ and } \quad
 \res_B :\beta^*\Tn(1)_B \ra \Tn(1)_{B'}\,. \]
We define $\psi_1$ as the composite
\[ \psi_1 :  R\pi_{X\!/\!B*}\Tn(r)_X \lra R\pi_{X\!/\!B*}R\alpha_*\Tn(r)_{X'} = R\beta_*R\pi_{X'\!/\!B'*}\Tn(r)_{X'}, \]
where the first arrow is the adjunction map of $\res_X$. We define $\psi_2^m$ as the composite
{\allowdisplaybreaks
\begin{align*}
\psi_2^m : \; & \fH^{\leqq m}(X,\Tn(r)) =
  R\sHom_{B,\Ln} (\tau_{\geqq 2(d-1)-m}R\pi_{X\!/\!B!}\Tn(d-r)_X,\Tn(1)_B)[2-2d]
 \phantom{|_{\big|}} \\
  &\;\; \lra R\beta_*R\sHom_{B',\Ln} (\tau_{\geqq 2(d-1)-m}
 \beta^*R\pi_{X\!/\!B!}\Tn(d-r)_X, \beta^*\Tn(1)_B)[2-2d]
 \phantom{|_{\big|}} \\
  &\;\; \lra  R\beta_*R\sHom_{B',\Ln}(\tau_{\geqq 2(d-1)-m}R\pi_{X'\!/\!B'!}\Tn(d-r)_{X'},\Tn(1)_{B'})[2-2d]
 \phantom{|_{\big|}} \\
  &\;\; = R\beta_*\fH^{\leqq m}(X',\Tn(r)),
\end{align*}
}where the second arrow is induced by $\res_B$ and the isomorphisms
{\allowdisplaybreaks
\begin{align*}
 \beta^*R\pi_{X\!/\!B!}\Tn(d-r)_X & \cong R\pi_{X'\!/\!B'!}\alpha^*\Tn(d-r)_X 
 \phantom{\big|_|} & \hbox{(proper base change)} \\
 & \cong R\pi_{X'\!/\!B'!}\Tn(d-r)_{X'} & \hbox{($r \geqq d$).}
\end{align*}
}We define $\psi_3^m$ in a similar way.
Note that the following square commutes by Corollary \ref{cor-trace}:
\begin{equation}\label{CD:filter}
\xymatrix{
R\pi_{X\!/\!B*}\Tn(r)_X \ar[d]_{\psi_1} \ar[r]_-{\eqref{filter}}^-\simeq &
\fH^{\leqq 2(d-1)}(X,\Tn(r)) \ar[d]^{\psi_2^{2(d-1)}} \\
R\beta_*R\pi_{X'\!/\!B'*}\Tn(r)_{X'} \ar[r]_-{\eqref{filter}}^-\simeq & R\beta_*\fH^{\leqq 2(d-1)}(X',\Tn(r)).}
\end{equation}
We define $\psi_4^m$ as the composite
\[\xymatrix{ \psi_4^m : i_{v*}Ri_v^!\fH^m(X,\Tn(r)) \ar[rr]^-{\text{base change}} &&
 i_{v*}R\iota_v^!\beta^*\fH^m(X,\Tn(r)) \os{\psi_3^m}\lra  i_{v*}R\iota_v^!\fH^m(X',\Tn(r)). }\]
See \cite{sga4}\,\XVIII.3.1.14.2 for the base change morphism.
This $\psi_4^m$ is an isomorphism, because both $Ri_v^!\fH^m(X,\Tn(r))$ and $R\iota_v^!\fH^m(X',\Tn(r))$ are isomorphic to
\begin{align*}
\begin{cases}
R\sHom_{v,\Ln}(R^{2(d-1)-m}\pi_{Y_v\nsp/\nsp v!}\Ln(d-r),\Ln)[-2] & \hbox{(if $\ch(v) \ne p$ or $r=d$)} \\
0 & \hbox{(if $\ch(v) = p$ and $r>d$)}
\end{cases}
\end{align*}
by \eqref{eq2-0} and Proposition \ref{prop:local-str}\,(1).
We prove that $\psi_1$, $\psi_2^m$ and $\psi_3^m$ are isomorphisms.
By the triangle \eqref{eq2-1} and the commutative diagram \eqref{CD:filter}, we are reduced to showing that $\psi_3^m$ is an isomorphism for any $m \in \bZ$.
Put $K':=\Frac(\fO')$, and let us note the following facts:
\begin{itemize}
\item[(i)]
$\H^m(X_{\ol K},\mu_{p^n}^{\otimes r}) \cong \H^m(X'_{\ol K{}'},\mu_{p^n}^{\otimes r})$,
 see \cite{milne}\,VI.4.3.
\item[(ii)]
$G_{\!K} \cong G_{\!K'}$, see \cite{milne:adual} p.\,160, Case 2 (i) and (ii).
\item[(iii)]
$\psi_4^m$ is an isomorphism.
\end{itemize}
By these facts and Proposition \ref{prop:local-str}\,(2),
 we see that $\psi_3^m$ is an isomorphism, which completes the proof of Proposition \ref{thm:rigidity}.
\end{pf*}

\section{Projective and inductive limits}\label{sect-limits}
Let $\pi_{X\!/\!B} : X \to B=\Spec(\fO)$ be as in \S\ref{sect1}.
We do {\it not} assume that $\pi_{X\!/\!B}$ is proper in this section,
but assume that $\fO$ and $K=\Frac(\fO)$ satisfy either of the following conditions:
\begin{enumerate}
\item[(L)]
$K$ is a non-archimedean local field of characterictic $0$, i.e., a finite field extension of $\bQ_\ell$ for some prime number $\ell$, and $\fO$ is the valuation ring of $K$.
\item[(G)]
$K$ is an algebraic number field, i.e., a finite field extension of $\bQ$, and $B=\Spec(\fO)$ is an open subset of $\Spec(\OK)$, where $\OK$ denotes the integer ring of $K$.
\end{enumerate}
The main aims of this section are to prove some standard finiteness results and to construct spectral sequences \eqref{ss:leray2}\sp--\sp\eqref{ss:leray4} below, under these assumptions.

\begin{prop}\label{prop1-1}
There is a canonical isomorphism
\begin{equation}\label{eq1-1}
\H^q(B,\fH^m(X,\Tn(r))) \cong \Ext_{B}^q(R^{2(d-1)-m}\pi_{X\!/\!B!}\Tn(d-r),\Gm)
\end{equation}
for any $q, m \geqq 0$, $n \geqq 1$ and $r \geqq d$.
Moreover, $\H^q(X,\Tn(r))$ and $\H^q(B,\fH^m(X,\Tn(r)))$ are finite for the same $(q, m, n, r)$.
\end{prop}
\begin{pf}
The isomorphism \eqref{eq1-1} follows from the definition of $\fH^m(X,\Tn(r))$ (see Definition \ref{def:mot}) and the canonical isomorphism
\[ R\sHom_B(\Ln,\Gm) \cong \Tn(1) \]
(a variant of \cite{Sa0} Proposition 4.5.1). See also \cite{JSS} (2.3.4).
The finiteness of the groups in \eqref{eq1-1} follows from the finiteness of Ext-groups in the Artin-Verdier duality (\cite{mazur}\,(2.4)) and the constructibility of $R^{2(d-1)-m}\pi_{X\!/\!B!}\Tn(d-r)$.
The finiteness of $\H^q(X,\Tn(r))$ follows from the spectral sequence \eqref{ss:leray} and that of $E_2$-terms.
\end{pf}

\subsection{Spectral sequences}
For $r \geqq d$, we introduce the following groups:
{\allowdisplaybreaks
\begin{align*}
\H^q(X,\zp(r)) := \varprojlim_{n \geqq 1} \ \H^q(X,\Tn(r)), \quad&\quad
\H^q(X,\qp(r)) := \H^q(X,\zp(r)) \otimes_{\zp}\qp, \phantom{\varprojlim_{n \geqq 1}}\\
\H^q(X,\QpZp(r)) &:= \varinjlim_{n \geqq 1} \ \H^q(X,\Tn(r)),\\
\H^q(B,\fH^m(X,\zp(r))) & := \varprojlim_{n \geqq 1} \ \H^q(B,\fH^m(X,\Tn(r))),\\
\H^q(B,\fH^m(X,\qp(r))) & := \H^q(B,\fH^m(X,\zp(r)))\otimes_{\zp} \qp,\phantom{\varprojlim_{n \geqq 1}}\\
\H^q(B,\fH^m(X,\QpZp(r))) & := \varinjlim_{n \geqq 1} \ \H^q(B,\fH^m(X,\Tn(r))).
\end{align*}
}Here the transition maps in the fourth group is defined by the commutative diagram
\[\xymatrix{
 \H^q(B,\fH^m(X,\fT_{n+1}(r))) \ar@{.>}[r] \ar[d]^{\cong}_{\eqref{eq1-1}} & \H^q(B,\fH^m(X,\Tn(r))) \ar[d]^{\cong}_{\eqref{eq1-1}} \\
 \Ext_B^q(R^{2(d-1)-m}\pi_{X\!/\!B!}\fT_{n+1}(d-r),\Gm) \ar[r] & \Ext_B^q(R^{2(d-1)-m}\pi_{X\!/\!B!}\Tn(d-r),\Gm)
}\]
with the bottom arrow induced by $\ul{p} : \Tn(d-r) \hra \fT_{n+1}(d-r)$ of Proposition \ref{prop:bock}.
The transition maps in the last group is defined by the commutative diagram
\[\xymatrix{
 \H^q(B,\fH^m(X,\Tn(r))) \ar@{.>}[r] \ar[d]^{\cong}_{\eqref{eq1-1}} & \H^q(B,\fH^m(X,\fT_{n+1}(r))) \ar[d]^{\cong}_{\eqref{eq1-1}} \\
 \Ext_B^q(R^{2(d-1)-m}\pi_{X\!/\!B!}\Tn(d-r),\Gm) \ar[r] & \Ext_B^q(R^{2(d-1)-m}\pi_{X\!/\!B!}\fT_{n+1}(d-r),\Gm)
}\]
with the bottom arrow induced by ${\mathscr R}^1 : \fT_{n+1}(d-r) \twoheadrightarrow \Tn(d-r)$
of Proposition \ref{prop:bock}.
Taking the projective limit of the spectral sequence \eqref{ss:leray} with respect to $n \geqq 1$,
 we obtain a convergent spectral sequence of $\zp$-modules
\begin{equation}\label{ss:leray2}
E^{a,b}_2=\H^a(B,\fH^b(X,\zp(r))) \Lra \H^{a+b}(X,\zp(r)).
\end{equation}
This spectral sequence yields a spectral sequence of $\qp$-vector spaces:
\begin{equation}\label{ss:leray3}
E^{a,b}_2=\H^a(B,\fH^b(X,\qp(r))) \Lra \H^{a+b}(X,\qp(r)).
\end{equation}
On the other hand, taking the inductive limit of \eqref{ss:leray} with respect to $n \geqq 1$,
 we obtain another convergent spectral sequence of $\zp$-modules
\begin{equation}\label{ss:leray4}
E^{a,b}_2=\H^a(B,\fH^b(X,\QpZp(r))) \Lra \H^{a+b}(X,\QpZp(r)).
\end{equation}

\subsection{Finite and cofinite generation}

The following preliminary results will be useful later:
\begin{thm}\label{lem:fg}
\begin{enumerate}
\item[{\rm(1)}]
 $\H^q(X,\zp(r))$ and $\H^q(B,\fH^m(X,\zp(r)))$ are finitely generated over $\zp$ for any $q, m \in \bZ$ and any $r \geqq d$.
\item[{\rm(2)}]
 $\H^q(X,\QpZp(r))$ and $\H^q(B,\fH^m(X,\QpZp(r)))$ are cofinitely generated over $\zp$ for any $q, m \in \bZ$ and any $r \geqq d$.
\item[{\rm(3)}]
We have\;
 $\rank_{\zp} \ \H^q(B,\fH^m(X,\zp(r))) =\corank_{\zp} \ \H^q(B,\fH^m(X,\QpZp(r)))$
 for any $q, m \in \bZ$ and any $r \geqq d$.
\end{enumerate}
\end{thm}
\begin{pf}
The assertions for $\H^q(X,\zp(r))$ and $\H^q(X,\QpZp(r))$ follow from a standard argument using Propositions \ref{prop1-1} and \ref{prop:bock}.
We prove the assertions for $\H^q(B,\fH^m(X,\zp(r)))$ and $\H^q(B,\fH^m(X,\QpZp(r)))$ in the case (G); the case (L) is similar and left to the reader.

We first show that $\H^q(B,\fH^m(X,\zp(r)))$ is finitely generated over $\zp$.
By the Artin-Verdier duality, it is enough to show that its Pontryagin dual
\[ \H^{3-q}_c(B,R^{m'}\!\pi_{X\!/\!B!}\qp/\zp(d-r)) :=
  \varinjlim_{n \geqq 1} \ \H^{3-q}_c(B,R^{m'}\pi_{X\!/\!B!}\fT_n(d-r)) \]
is cofinitely generated over $\zp$, where $m':=2(d-1)-m$.
Let $M^s_\Div$ to be the maximal $p$-divisible subsheaf of $M^s:=R^s\pi_{X\!/\!B!}\qp/\zp(d-r)$, i.e.,
\[ M^s_\Div:=\Image \Big( \sHom_B(\ul{\qp},M^s) \ra M^s \Big), \]
where $\ul{\qp}$ denotes the constant sheaf on $B_{\et}$ with values in $\qp$.
For each $n \geqq 1$, put ${}_{p^n}(M^s_\Div):=\Ker(\times p^n : M^s_\Div \to M^s_\Div)$,
which is a subquotient of $R^s\pi_{X\!/\!B!}\fT_n(d-r)$, hence constructible.
Moreover, there is a short exact sequence
\[ 0  \lra {}_{p^n}(M^s_\Div) \lra {}_{p^{n+n'}}(M^s_\Div) \lra {}_{p^{n'}}(M^s_\Div) \lra 0 \]
for any $n, n' \geqq 1$, and $\H^i_c(B,M^s_\Div)$ is cofinitely generated over $\zp$ for any $i$ by a standard argument. On the other hand, the quotient sheaf $M^s_\cotor:=M^s/M^s_\Div$ is the torsion part of $R^{s+1}\pi_{X\!/\!B!} \zp$, hence constructible (\cite{sga5}\,VI.2.2.2), and $\H^i_c(B,M^s_\cotor)$ is finite for any $i$.
Therefore by the long exact sequence
\begin{equation}\label{eq3-1}
 \dotsb \to \H^i_c(B,M^s_\Div) \to \H^i_c(B,M^s) \to \H^i_c(B,M^s_\cotor) \to \H^{i+1}_c(B,M^s_\Div) \to \dotsb,
\end{equation}
$\H^i_c(B,M^s)$ is cofinitely generated over $\zp$ for any $i$ and $s$, and $\H^q(B,\fH^m(X,\zp(r)))$ is finitely generated over $\zp$ for any $q$ and $m$.
\par
We next show that $\H^q(B,\fH^m(X,\QpZp(r)))$ is cofinitely generated over $\zp$.
By similar arguments as before, it is enough to show that the group
\[ \H^{3-q}_c(B,R^s \pi_{X\!/\!B!}\zp(d-r)) :=
  \varprojlim_{n \geqq 1} \ \H^{3-q}_c(B,R^s\pi_{X\!/\!B!}\fT_n(d-r)) \]
is finitely generated over $\zp$ for any $i$ and $s$.
Let $M^s$ and $M^s_\Div$ be as before, and put
\[ T^s_n:={}_{p^n}(M^s_\Div) \quad (n \geqq 1) \quad \hbox{and} \quad T^s:=(T^s_n)_{n \geqq 1} \]
Note that $T^s:=(T^s_n)_{n \geqq 1}$ is a constructible $\zp$-sheaf. Put
\[ \H^i_c(B,T^s) := \varprojlim_{n \geqq 1} \ H^i_c(B,T^s_n), \]
which is finitely generated over $\zp$ for any $i$ by a standard argument.
Noting that there is a short exact sequence of constructible $\zp$-sheaves
\[ 0\lra M^{s-1}_\cotor \lra  R^s\pi_{X\!/\!B!} \zp(d-r) \lra T^s \lra 0 \]
we obtain a long exact sequence
\begin{align}
\notag
 \dotsb \to \H^i_c(B,M^{s-1}_\cotor) \to \H^i_c(B,R^s\pi_{X\!/\!B!}\zp(d-r)) & \to \H^i_c(B,T^s)
\phantom{|_{\big|}} \\
& \to \H^{i+1}_c(B,M^{s-1}_\cotor) \to \dotsb,
 \label{eq3-2}
\end{align}
which shows that $\H^i_c(B,R^s\pi_{X\!/\!B!}\zp(d-r))$ is finitely generated over $\zp$ for any $i$.
Finally from the long exact sequence of $\zp$-modules
\begin{align}
\notag
 \dotsb \to \H^i_c(B,T^s) \to \H^i_c(B,R^s\pi_{X\!/\!B!}\qp(d-r))& \to \H^i_c(B,M^s_\Div)
\phantom{|_{\big|}} \\
& \to \H^{i+1}_c(B,T^s) \to \dotsb,
 \label{eq3-3}
\end{align}
we obtain
{\allowdisplaybreaks
\begin{align*}
 & \rank_{\sp \zp} \ \H^q(B,\fH^m(X,\zp(r)))
 \os{\text{(duality)}}= \corank_{\sp \zp} \ \H^{3-q}_c(B,M^{m'}) \quad \hbox{($m':=2(d-1)-m$)}
 \phantom{|_{\big|}} \\
 & \os{\eqref{eq3-1}}= \corank_{\sp \zp} \ \H^{3-q}_c(B,M^{m'}_\Div)
 \os{\eqref{eq3-3}}= \rank_{\sp \zp} \ \H^{3-q}_c(B,T^{m'}) \phantom{|_{\big|}} \\
 & \os{\eqref{eq3-2}}= \rank_{\sp \zp} \ \H^{3-q}_c(B,R^{m'}\!\pi_{X\!/\!B!}\zp(d-r))
  \os{\text{(duality)}}= \corank_{\sp \zp} \ \H^q(B,\fH^m(X,\QpZp(r))),
\end{align*}
}which shows the assertion (3).
\end{pf}

\section{Comparison with Selmer groups,\, local case}\label{sect3}
Let $\pi_{X\!/\!B} : X \to B=\Spec(\fO)$ be as in \S\ref{sect1}.
In this section, we always assume the following:
\begin{itemize}
\item
$\pi_{X\!/\!B}$ is {\it proper}, and the generic fiber $X_{\nsp K}$ is {\it geometrically connected over $K$}.
\item
$K$ is a non-archimedean local field of characteristic $0$, and $\fO$ is the valuation ring of $K$, i.e., the case (L) of \S\ref{sect-limits}.
\end{itemize}
Let $k$ be the residue field of $\fO$ and put $\ell:=\ch(k)$.
We will often write $Y$ (resp.\ $\ol Y$) for $X \otimes_\fO k$ (resp.\ $X \otimes_\fO \ol k$).
We put $V^m:=\H^m(X_{\ol K},\qp)$.
A main aim of this section is to compare $\H^m(X,\qp(r))$ with $\H^1_{\!f}(K,V^{m-1}(r))$, the local Selmer group of Bloch-Kato \cite{BK} \S3.
We will often write $\H^1_{\ovf}(K,-)$ for the quotient of $\H^1(K,-)$ by $\H^1_{\!f}(K,-)$.
The following fundamental fact (loc.\ cit., Proposition 3.8) will be useful:
\begin{lem}[{{\bf Bloch-Kato}}]\label{thm:fundbk}
Let $V$ be a finite-dimensional $\qp$-vector space on which the Galois group $G_{\nsp\nsp K}$ acts continuously.
If $\ell=p$, then assume further that $V$ is a de Rham representation.
Put $V^*:=\Hom_{\qp}(V,\qp)$.
Then under the perfect pairing of local Tate duality
\[ \H^1(K,V) \times \H^1(K,V^*(1)) \ra \H^2(K,\qp(1)) \cong \qp, \]
the subspaces $\H^1_{\!f}(K,V)$ and $\H^1_{\!f}(K,V^*(1))$ are the exact annihilators of each other.
\end{lem}

The following standard fact will be useful later in \S\S\ref{sect5}\sp--\ref{sect7} below.
\begin{lem}\label{lem4-1}
Assume that $\ell \ne p$, and that $\pi_{X\!/\!B} :  X \to B$ is smooth and proper.
Then we have $\H^a(B,\fH^m(X,\Tn(r))) = 0$ for any $a \geqq 2$, $m \geqq 0$, $n \geqq 1$ and $r \geqq d$.
\end{lem}
\begin{pf}
Under the assumptions,
$\fH^m(X,\Tn(r))$ is a locally constant sheaf on $B_\et$ placed in degree $0$, whose stalk at $\ol v$ is $\H^m(\ol Y,\mu_{p^n}^{\otimes r})$ by Lemma \ref{lem1-0}\,(1), Proposition \ref{ex:trace}\,(1) and the proper smooth base change theorem. Hence we have
\[ \H^a(B,\fH^m(X,\Tn(r))) \cong \H^a(v,\H^m(\ol Y,\mu_{p^n}^{\otimes r})) =0 \]
for any $a \geqq 2$, as claimed.
\end{pf}

\subsection{Comparison results}
The main result of this section is the following:
\begin{thm}\label{thm:local-cond}
For any $m \geqq 0$ and $r \geqq d$, we have canonical isomorphisms
\[ \H^q(B,\fH^m(X,\qp(r))) \cong
\begin{cases}
    \H^1_{\!f}(K,V^m(r)) &  \hbox{{\rm(}$q=1${\rm)}} \phantom{|_{\big|}}\\
      0  &  \hbox{{\rm(}otherwise{\rm)}}
\end{cases}
\]
Moreover, if $\ell \ne p$, then we have $\H^1(B,\fH^m(X,\qp(r)))=0$ for any $m \geqq 0$ and $r \geqq d$.
\end{thm}
\begin{rem}\label{rem:con}
{\rm If $\ell \not=p$, then we have $\H^m(X,\qp(r))=0$ for any $m \in \bZ$ and $r \geqq d$ by the proper base change theorem
\[ \H^m(X,\qp(r)) \cong \H^m(Y,\qp(r)) \]
and a theorem of Deligne {\rm\cite{d2}}\,3.3.4 on weights of $\H^*(\ol{Y},\qp)$ {\rm(}note that $\dim(Y)=d-1${\rm)}. Theorem \ref{thm:local-cond} for $\ell \ne p$ refines this fact.
}
\end{rem}
\noindent
We first state a few consequences of Theorem \ref{thm:local-cond}.
By the theorem and the spectral sequence \eqref{ss:leray3}, we obtain the following corollary:

\begin{cor}\label{cor:local-cond}
The spectral sequence \eqref{ss:leray3} degenerates at $E_2$, and we have
\[ \H^m(X,\qp(r)) \cong \H^1_{\!f}(K,V^{m-1}(r)) \]
for any $m \geqq 0$ and any $r \geqq d$.
\end{cor}

\noindent
The following corollary will be useful later:
\begin{cor}\label{cor4-fin}
\begin{enumerate}
\item[{\rm (1)}]
There exists a natural map
\[ H^1_{\!f}(K,V^m(r)) \lra \H^1(B,\fH^m(X,\QpZp(r))) \]
which fits into a commutative diagram
\[\xymatrix{
H^1_{\!f}(K,V^m(r)) \ar@{.>}[d]\ar[rd]^{\text{{\rm(}natural map{\rm)}}} \\
\H^1(B,\fH^m(X,\QpZp(r))) \; \ar@{^{(}->}[r] & \H^1(K,\H^m(X_{\ol K},\QpZp(r))).
}\]
See Proposition \ref{prop:local-str}\,{\rm(}1{\rm)} for the injectivity of the bottom arrow.
\item[{\rm (2)}]
$\H^a(B,\fH^m(X,\zp(r)))$ and $\H^a(B,\fH^m(X,\QpZp(r)))$ are finite for any $a \ne 1$, $m \geqq 0$ and $r \geqq d$.
\end{enumerate}
\end{cor}
\begin{pf}
The claim (1) immediately follows from Theorem \ref{thm:local-cond},
and the claim (2) follows from Theorems \ref{thm:local-cond} and \ref{lem:fg}.
\end{pf}

\medskip

We start the proof of Theorem \ref{thm:local-cond}.
A key step is to show Theorem \ref{thm:local-cond'} below.
Fix integers $m \geqq 0$ and $r \geqq d$, and put
\[ \H^q(B,R^m\pi_{X\!/\!B*}\qp(d-r)) := \qp \otimes_{\zp} \varprojlim_n \ \H^q(B,R^m\pi_{X\!/\!B*}\Tn(d-r)), \]
Under this notation, we will prove
\begin{thm}\label{thm:local-cond'}
We have
\begin{equation}\label{isom:local-cond2}
\H^q(B,R^m\pi_{X\!/\!B*}\qp(d-r)) \cong
  \begin{cases}
    V^m(d-r)^{G_{\nsp\nsp K}} \phantom{|_{|}}& \hbox{{\rm(}$q=0${\rm)}}\\
    \H^1_{\!f}(K,V^m(d-r)) \quad \phantom{|_{|}}& \hbox{{\rm(}$q=1${\rm)}}\\
    0 &  \hbox{{\rm(}$q \ne 0,1${\rm)}}\\
  \end{cases}
\end{equation}
and
\begin{equation}\label{isom:local-cond2'}
 V^m(r)^{G_{\nsp\nsp K}}=0.
\end{equation}
\end{thm}
\noindent
We have $\H^q(B,R^m\pi_{X\!/\!B*}\qp(d-r))=0$ if $r > d$.
Indeed, if $\ell = p$, then this vanishing follows from the definition of $\Tn(d-r)$ and the proper base change theorem; the case $\ell \ne p$ follows from similar arguments as in Remark \ref{rem:con}.
Consequently, the isomorphism \eqref{isom:local-cond2} asserts the vanishing of the right hand side for $r >d$.
We will prove Theorem \ref{thm:local-cond'} in \S\ref{sect2-2} and \S\ref{sect2-3} below.
\par\bigskip

\begin{pf*}{{\it Proof of ``Theorem \ref{thm:local-cond'} $\Lra$ Theorem \ref{thm:local-cond}''}}
Let $v$ be the closed point of $B$.
Put $s:=d-r (\leqq 0)$ for simplicity.
By the isomorphisms in \eqref{isom:local-cond2} and the localization long exact sequence
\begin{align*}
 \dotsb & \lra \H^{q-1}_v(B,R^m\pi_{X\!/\!B*}\qp(s)) \lra \H^{q-1}(B,R^m\pi_{X\!/\!B*}\qp(s)) \lra \H^{q-1}(K,V^m(s)) \\ & \lra \H^q_v(B,R^m\pi_{X\!/\!B*}\qp(s)) \lra \dotsb \phantom{\big|^{\big|}}
\end{align*}
we have
\begin{equation}\label{local1}
\H^q_v(B,R^m\pi_{X\!/\!B*}\qp(s)) \cong
  \begin{cases}
     0  \phantom{\big|_|}&  \hbox{($q \ne 2,3$)} \\
    \H^1_{\ovf}(K,V^m(s)) \quad \phantom{\big|_|}&  \hbox{($q=2$)}\\
    \H^2(K,V^m(s)) &  \hbox{($q=3$).}
  \end{cases}
\end{equation}
Theorem \ref{thm:local-cond} for $q \ne 0$ follows from \eqref{local1} with $2(d-1)-m$ in place of $m$, Lemma \ref{thm:fundbk} and the Tate duality for cohomology of $B$ (see \cite{mazur}\,(2.4)):
\[ \H^q(B,\fH^m(X,\qp(r))) \times \H^{3-q}_v(B,R^{2(d-1)-m}\pi_{X\!/\!B*}\qp(s)) \lra \H^{3}_v(B,\qp(1)) \cong \qp. \]
The assertion for $q = 0$ of Theorem \ref{thm:local-cond} is a consequence of the isomorphism
\[ \H^0(B,\fH^m(X,\qp(r))) \cong V^m(r)^{G_{\nsp\nsp K}} \]
(see Proposition \ref{prop:local-str}\,(1)) and the vanishing \eqref{isom:local-cond2'}. Finally if $\ell \ne p$, then $\H^1_{\!f}(K,V^m(r)) = \H^1(k,V^m(r)^{I_K}) =0$ again by \eqref{isom:local-cond2'} and the equality of dimensions
\begin{equation}\label{dim:equal}
\dim_{\qp} V^m(r)^{G_{\nsp\nsp K}} = \dim_{\qp} \H^1(k,V^m(r)^{I_K}),
\end{equation}
which is a consequence of the duality of Galois cohomology of $G_k$.
\end{pf*}

\subsection{Proof of Theorem \ref{thm:local-cond'}\, (the case $\bs{\ell \ne p}$)}\label{sect2-2}
Let $K^\ur$ be the maximal unramified extension of $K$, and let $I_K=\Gal(\ol K/K^\ur)$ be the inertia group of $K$.
Let $\fO^\ur$ be the valuation ring of $K^\ur$, and let $\cosp^m_X$ be the cospecialization map
\begin{equation}\label{def:cosp}
\cosp^m_X : \H^m(\ol Y,\qp) \cong \H^m(X^\ur,\qp) \lra \H^m(X_{\ol K},\qp)^{I_K} = (V^m)^{I_K}
\end{equation}
for $m \geqq 0$, where $X^\ur$ (resp.\ $\ol Y$) denotes $X \otimes_\fO \fO^\ur$ (resp.\ $Y \otimes_k \ol k$).
We first reduce Theorem \ref{thm:local-cond'} for $\ell \ne p$ to the following proposition:

\begin{prop}\label{prop:local-cond}
Assume that $\ell \not=p$, and let $m \geqq 0$ be an integer. Then{\rm:}
\begin{enumerate}
\item[{\rm(1)}]
We have $V^m(r)^{G_{\nsp\nsp K}}=0$\, for any $r \geqq d$.
\item[{\rm(2)}]
For any $s \leqq 0$ and $q=0,1$, the map $\cosp^m_X$ induces an isomorphism
\[ \H^q(k,\H^m(\ol Y,\qp(s))) \cong \H^q(k,V^m(s)^{I_K}). \]
\end{enumerate}
\end{prop}
\noindent
Proposition \ref{prop:local-cond}\,(1) is the same as \eqref{isom:local-cond2'} of Theorem \ref{thm:local-cond'}.

\bigskip

\begin{pf*}{{\it Proof of ``Proposition \ref{prop:local-cond} $\Lra$ Theorem \ref{thm:local-cond'}''}}
We have
\[ \H^q(B,R^m\pi_{X\!/\!B*}\qp(s)) \cong \H^q(k,\H^m(\ol Y,\qp(s))) \]
and the last group is zero unless $q=0$ or $1$, because $\cd(G_k)=1$.
The isomorphisms for $q=0,1$ of \eqref{isom:local-cond2} follow from Proposition \ref{prop:local-cond}\,(1) and the fact that
\[ \H^1_{\!f}(K,V^m(s)) = \H^1(k,V^m(s)^{I_K}) \]
by definition. Thus we obtain Theorem \ref{thm:local-cond'}, admitting Proposition \ref{prop:local-cond}.
\end{pf*}
\medskip
\begin{pf*}{Proof of Proposition \ref{prop:local-cond}}
If $X$ is smooth over $B$, then the assertions are clear by the proper smooth base change theorem and Deligne's proof of the Weil conjecture \cite{d2}\,3.3.9.
We are concerned with the case that $\pi_{X\!/\!B}:X \ra B$ is not smooth, in what follows.
\par
\bigskip
(I)\, {\it Strict semi-stable reduction case.}\;
We first prove Proposition \ref{prop:local-cond} assuming that $X$ has strict semi-stable reduction.
We introduce some notation. Let $\ol j$ be the canonical map $X_{\ol K} \ra X^\ur=X \otimes_\fO \fO^\ur$, and let $\ol {\iota}$ be the closed immersion $\ol Y \ra X^\ur$.
By the properness of $X/B$, we have the following Leray spectral sequence for any $n \geqq 1$:
\begin{equation}\label{ss:rz}
E_2^{a,b}=\H^a(\ol Y,\ol {\iota}{}^*R^b{\ol j}_*\Ln) \Lra \H^{a+b}(X_{\ol K},\Ln).
\end{equation}
By a theorem of Rapoport and Zink \cite{rz}\,2.23, there is an exact sequence on $(\ol Y)_\et$
\begin{align}
\notag
 0 \lra  \ol {\iota}{}^*R^b{\ol j}_*\Ln \lra u^{b+1}_*\Ln(-b)_{Z^{(b+1)}}
 & \lra u^{b+2}_*\Ln(-b)_{Z^{(b+2)}} \lra \\
\label{exact:rz}
 & \dotsb \lra u^{d}_*\Ln(-b)_{Z^{(d)}} \lra 0,
\end{align}
where for each $m > 0$, $Z^{(m)}$ denotes the disjoint union of $m$-fold intersections distinct irreducible components of $\ol Y$ and $u^{m}$ denotes the canonical (finite) map $Z^{(m)} \ra \ol Y$;
see \eqref{eqdef1-1} for $\Ln(-b)$. Hence the $E_2$-terms of the spectral sequence of \eqref{ss:rz} are finite and we obtain a spectral sequence
\begin{equation}\label{ss:rz2}
 E_2^{a,b}=\H^a(\ol Y,\ol {\iota}{}^*R^b{\ol j}_*\qp) \Lra \H^{a+b}(X_{\ol K},\qp)=V^{a+b}
\end{equation}
by taking the projective limit with respect to $n \geqq 1$ and the tensor product with $\qp$ over $\zp$.
Note that the canonical map $E_2^{m,0}=\H^m(\ol Y,\qp) \ra E^m=V^m$ agrees with the cospecialization map $\cosp_X^m$ of \eqref{def:cosp}, and that the inertia group $I_K$ acts trivially on the $E_2$-terms of \eqref{ss:rz2}.
We will prove the following:
\begin{lem}\label{lem-claim}
In the spectral sequence \eqref{ss:rz2}, we have $E_2^{a,b}=0$ unless $0 \leqq a \leqq 2(d-b-1)$ and $0 \leqq b \leqq d-1$.
Furthermore, for a pair $(a,b)$ with $0 \leqq a \leqq 2(d-b-1)$ and $0 \leqq b \leqq d-1$,
 the weights of $E_2^{a,b}$ are at least $\max\{2b,2(a+2b+1-d)\}$ and at most $a+2b$.
\end{lem}
By this lemma, the kernel and the cokernel of the map $\cosp_X^m$ in \eqref{def:cosp} have only positive weights and hence we obtain the assertion of Proposition \ref{prop:local-cond}\,(2).
Similarly, one can easily derive Proposition \ref{prop:local-cond}\,(1) from this lemma.
\par\bigskip

\begin{pf*}{\it Proof of Lemma \ref{lem-claim}}
By \eqref{exact:rz}, the sheaf $\ol {\iota}{}^*R^b{\ol j}_*\Ln$ (hence $E_2^{a,b}$ of \eqref{ss:rz2}) is zero unless $0 \leqq b \leqq d-1$.
Fix a $b \geqq 0$ in what follows. By the exact sequence \eqref{exact:rz}, we have a spectral sequence of finite-dimensional $G_k$-$\qp$-vector spaces:
\begin{equation}\label{ss:rz4}
'\!E_1^{s,t}=\H^t(Z^{(s+b+1)},\qp(-b))
   \Lra \H^{s+t}(\ol Y,\ol {\iota}{}^*R^b \ol j _*\qp),
\end{equation}
Here $'\!E_1^{s,t}$ is zero unless
\begin{equation}\label{uneq1}
0 \leqq t \leqq 2(d-s-b-1) \quad \hbox{ and } \quad 0 \leqq s \leqq d-b-1,
\end{equation}
because $\dim(Z^{(s+b+1)})=d-s-b-1$ and $Z^{(s+b+1)}=\emptyset$ if $s+b \geqq d$.
Using this spectral sequence, one can easily check that $E_2^{a,b}$ of \eqref{ss:rz2}) is zero unless $0 \leqq a \leqq 2(d-b-1)$. Moreover, $'\!E_1^{s,t}$ has weight $t+2b$ by \cite{d2}\,3.3.9.
Therefore one obtains the lemma by computing the span of $t+2b$ under the conditions \eqref{uneq1} and $a=s+t$.
\end{pf*}
\par\smallskip\noindent
This completes the proof Proposition \ref{prop:local-cond} in the strict semi-stable reduction case.
\par
\bigskip
(\II)\, {\it General case.}\;
We prove Proposition \ref{prop:local-cond} in the general case.
By the alteration theorem of de Jong \cite{dj}\,6.5, there exists a proper generically \'etale morphism $f : X' \ra X$ such that $X'$ is regular and flat over $B$ and has strict semi-stable reduction over the normalization $B'$ of $B$ in $X'$.
Let $L$ (resp.\ $k'$) be the function field of $B'$ (resp.\ the residue field of $L'$), $Y'$ for the special fiber of $\pi_{X'\!/\!B'}:X' \ra B'$.
Then Proposition \ref{prop:local-cond}\,(2) immediately follows from those for $X'$, proved in Step (I), and the fact that $V^m=\H^m(X_{\ol K},\qp)$ is a direct summand of $\H^m(X'_{\ol L},\qp)$ as $G_L$-$\qp$-vector spaces.
To prove Proposition \ref{prop:local-cond}\,(1), we consider the following commutative diagram:
\[
\xymatrix{
\H^q(k,\H^m(\ol Y,\qp(s))) \ar[r]^-{f^\sharp} \ar[d]_{\cosp_X^m}
    & \H^q(k',\H^m(\ol{Y'},\qp(s))) \ar[r]^-{\tr_f} \ar[d]_{\cosp_{X'}^m}
    & \H^q(k,\H^m(\ol Y,\qp(s))) \ar[d]_{\cosp_X^m} \\
\H^q(k,V^m(s)^{I_K}) \ar[r]^-{f^\sharp}
    & \H^q(k',\H^m(X'_{\ol L},\qp(s))^{I_L}) \ar[r]^-{\tr_f}
    & \H^q(k,V^m(s)^{I_K}),}
\]
where the right horizontal arrows are induced by the following homomorphism of \'etale sheaves on $B$:
\begin{multline*}
\tr_f : \pi_{B'\!/\!B*}R^m\pi_{X'\!/\!B'*}\Ln(s)_{X'} \cong R^m\pi_{X'\!/\!B*}\Ln(s)_{X'} \os{(*)}\cong
   R^m\pi_{X'/\!B*}(Rf^!\Ln(s)_{X}) \\ 
\xymatrix{= R^m\pi_{X/\!B*}(Rf_*Rf^!\Ln(s)_{X}) \ar[rr]^-{\text{adjunction}} && \, R^m\pi_{X/\!B*}\Ln(s)_{X}}
\end{multline*}
and we have used the absolute purity \cite{FG} to obtain the isomorphism $(*)$.
Since the middle vertical arrow in the above diagram is bijective by Step (I), the assertion of Proposition \ref{prop:local-cond}\,(1) for $X$ follows from the fact that the composite map
\[ R^m\pi_{X/\!B*}\Ln(s)_{X} \os{f^\sharp}{\lra} \pi_{B'\!/\!B*}R^m\pi_{X'\!/\!B'*}\Ln(s)_{X'}
 \os{\tr_f}{\lra} R^m\pi_{X/\!B*}\Ln(s)_{X} \]
on $B_\et$ agrees with the multiplication by the extension degree of function fields of $f : X' \to X$. This completes the proof.
\end{pf*}

\subsection{Proof of Theorem \ref{thm:local-cond'}\, (the case $\bs{\ell = p}$)}\label{sect2-3}
By the same arguments as in the proof of ``Proposition \ref{prop:local-cond} $\Ra$ Theorem \ref{thm:local-cond'}'' in \S\ref{sect2-2}, the assertions of Theorem \ref{thm:local-cond'} with $\ell=p$ is reduced to the following:
\begin{prop}\label{prop:local-cond2}
Assume that $\ell=p$. Let $m \geqq 0$ be an integer, and put $V:=V^m$. Then{\rm:}
\begin{enumerate}
\item[{\rm(1)}]
We have $V(r)^{I_K} = 0$ for any $r \geqq d$.
\item[{\rm(2)}]
For any $s < 0$, we have $V(s)^{I_K}=0$. For $s=0$, the cospecialization map
\[ \cosp^m_X : \H^m(\ol Y,\qp) \lra \H^m(X_{\ol K},\qp)^{I_K}=V^{I_K} \]
is bijective.
\item[{\rm(3)}]
For any $s \leqq 0$, we have $\H^1(k,V(s)^{I_K}) = \H^1_{\!f}(K,V(s))$ in $\H^1(K,V(s))$.
In particular, we have $\H^1_{\!f}(K,V(s))=0$ if $s < 0$.
\end{enumerate}
\end{prop}
We will first prove Proposition \ref{prop:local-cond2} assuming that $X$ has semi-stable reduction, and then prove the log smooth reduction case.
\par\bigskip
\begin{pf}
(I)\, {\it Semi-stable reduction case.}\;
See \cite{fo}\,1.5.5 for $B_\crys, B_\st, B_\dR^+$ and $B_{\dR}$.
Put $D:=\H^m_{\logcrys}(Y/\Ws(k))$.
By the Fontaine-Jannsen conjecture (\cite{hk}, \cite{tsuji}\,0.2), we have a $p$-adic period isomorphism
\begin{equation}\label{p-Hodge}
 V \otimes_{\qp} B_{\st} \cong D \otimes_{\Ws(k)} B_{\st},
\end{equation}
which preserves the Frobenius operator $\phi$, the monodromy operator $N$, the action of $G_{\nsp\nsp K}$, and the Hodge filtration $\rF^{\bullet}_H$ after taking $\otimes_{B_{\st}}B_{\dR}$.
By the isomorphism \eqref{p-Hodge}, we have
\[ V(r) \cong \left(D \otimes_{\Ws(k)}B_{\st} \right)^{N=0,\,\phi=p^r}
 \cap \rF_H^r\left(D \otimes_{\Ws(k)}B_{\dR} \right) \]
and
\begin{equation}\label{p-Hodge1}
V(r)^{I_K} \subset (\H^{m}_{\logcrys}(\ol Y/\Ws(\ol k))_{\qp})^{\varphi=p^r},
\end{equation}
for any $r \in \bZ$.
Here $\varphi$ denotes the Frobenius operator acting on $\H^{m}_{\logcrys}(\ol Y/\Ws(\ol k))$, and we have used the fact that $(B_{\st})^{I_K} = \Frac(\Ws(\ol k))$ (\cite{fo} 5.1.2, 5.1.3). Proposition \ref{prop:local-cond2}\,(1) and the case $s < 0$ of Proposition \ref{prop:local-cond2}\,(2) follow from \eqref{p-Hodge1} and the fact that  
\[ (\H^{m}_{\logcrys}(\ol Y/\Ws(\ol k))_{\qp})^{\varphi=p^r}=0 \quad\; \hbox{if \sp $r \geqq d$ \sp or \sp $r<0$.} \]
As for the case $s=0$ of Proposition \ref{prop:local-cond2}\,(2), the map $\cosp^m$ is bijective by \cite{Wu} Theorem 1. To prove Proposition \ref{prop:local-cond2}\,(3), it is enough to show the following two claims:
{\it
\begin{enumerate}
\item[{\rm(i)}]
The restriction map
\[ \H^1(K,V \otimes_{\qp} B_{\crys}) \lra \H^1(K^\ur,V \otimes_{\qp} B_{\crys}) \]
is injective. Consequently, the image of the inflation map
\[ \xymatrix{ \H^1(k,V(s)^{I_K}) \, \ar@<-1pt>@{^{(}->}[r] & \H^1(K,V(s)) } \]
is contained in $\H^1_{\!f}(K,V(s))$ for any $s \in \bZ$.
\item[{\rm(ii)}]
We have $\dim_{\qp}\ \H^1(k,V(s)^{I_K}) = \dim_{\qp}\ \H^1_{\!f}(K,V(s))$ for any $s \leqq 0$.
\end{enumerate}}
\par\bigskip\noindent
{\it Proof of the claim} (i).\;
By the inflation-restriction exact sequence
\[ (0 \to)\, \H^1(k, (V \otimes_{\qp} B_{\crys})^{I_K}) \lra
 \H^1(K,V \otimes_{\qp} B_{\crys}) \lra \H^1(K^{\ur},V \otimes_{\qp} B_{\crys})^{G_k}, \]
it is enough to show that the first term is zero. We have
\[ (V \otimes_{\qp} B_{\crys})^{I_K} \cong \H_{\logcrys}^m(\ol Y/\Ws(\ol k))_{\qp}^{N=0} \]
by the exact sequence (\cite{fo}\,3.2.3)
\[ 0 \lra B_{\crys} \lra B_{\st} \os{N}{\lra} B_{\st} \lra 0 \]
and the period isomorphism \eqref{p-Hodge}. Hence we have
\[ H^1(k,(V \otimes_{\qp} B_{\crys})^{I_K}) \cong
 \qp \otimes_{\zp} \varprojlim_{n \geqq 1} \ H^1(k,\H_{\logcrys}^m(\ol Y/\Wn(\ol k))^{N=0}). \]
Finally, the group on the right hand side is zero, because $\H_{\logcrys}^m(\ol Y/\Wn(\ol k))^{N=0}$ is a finite successive extension of $G_k$-modules which are isomorphic to the additive group of $\ol k$.
\par
\bigskip
\noindent
{\it Proof of the claim} (ii).\;
Since $V$ is a de Rham representation \cite{fal}, there is an exact sequence of finite-dimensional $\qp$-vector spaces (\cite{BK}\,Corollary 3.8.4):
\begin{multline}
0 \lra V(s)^{G_{\nsp\nsp K}} \lra \Cris(V) \oplus \DR(V(s))^0 \\
  \lra \Cris(V) \oplus \DR(V) \lra \H^1_{\!f}(K,V(s)) \lra 0,
\label{exact:bk}
\end{multline}
where $\Cris(V)$, $\DR(V(s))^0$ and $\DR(V)$ denote $(V \otimes_{\qp} B_{\crys})^{G_{\nsp\nsp K}}$, $(V(s) \otimes_{\qp} B_{\dR}^+)^{G_{\nsp\nsp K}}$ and $(V \otimes_{\qp} B_{\dR})^{G_{\nsp\nsp K}}$, respectively.
Moreover we have
\begin{equation}\label{same:faltings}
 \DR(V) \cong \H^m_{\dR}(X_{\nsp K}/K) = \rF_H^s \H^m_{\dR}(X_{\nsp K}/K) \cong \DR(V(s))^0
\end{equation}
for any $s \leqq 0$. Hence we obtain the claim (ii) from the equalities
\[ \dim_{\qp}\ \H^1_{\!f}(K,V(s)) \os{\eqref{exact:bk}}= \dim_{\qp}\ V(s)^{G_{\nsp\nsp K}} = \dim_{\qp}\ \H^1(k,V(s)^{I_K}), \]
where the right equality is similar to the equality of \eqref{dim:equal}.
This completes the proof of Proposition \ref{prop:local-cond2} in the semi-stable reduction case.
\pbn

(\II)\, {\it Log smooth reduction case.}\;
Let $f : X' \to X$, $B$ and $L$ be as in Step (\II) in the proof of Proposition \ref{prop:local-cond}.
The assertions in Proposition \ref{prop:local-cond2} other than the bijectivity of $\cosp^m$ are reduced to the semi-stable reduction case directly by a standard norm argument for $f_{\nsp K} : X'_{\nsp L} \to X_{\nsp K}$.
We derive the bijectivity of $\cosp^m$ for $X$ from that for $X'$.
Indeed, there exists a homomorphism $\tr_f: \pi_{B'\!/\!B*}R^m\pi_{X'\!/\!B'*}\Ln \to R^m\pi_{X/\!B*}\Ln$ of sheaves on $B_\et$ for each $n \geqq 1$ given by the following left commutative square, which is by definition the Pontryagin dual of the right commutative square ($m':=2d'-m$):
\[\xymatrix{ \H^m(\ol{Y'},\Ln) \ar[r] \ar[d]_{\cosp_{X'}^m}
    & \H^m(\ol Y,\Ln) \ar[d]_{\cosp_{X}^m} \\
\H^m(X'_{\ol L},\Ln)^{I_L} \ar[r]
    & \H^m(X_{\ol K},\Ln)^{I_K} }
\quad
\xymatrix{ \H^{m'+2}_{\ol{Y'}}((X')^\ur,\fT_n(d)) 
    & \ar[l]_-{\;\, f^\sharp}  \H^{m'+2}_{\ol Y}(X^\ur,\fT_n(d))  \\
\H^{m'}(X'_{\ol L},\mu_{p^n}^{\otimes d-1})_{I_L} \ar[u]^{\res_{X'}}
    & \ar[l]_-{\;\, f^\sharp} \H^{m'}(X_{\ol K},\mu_{p^n}^{\otimes d-1})_{I_K} \ar[u]^{\res_X}
}\]
where we put $(X')^\ur:=X' \times_{B'} (B')^{\ur}$ and $X^\ur:=X \times_B B^{\ur}$,
and the right square is the commutative diagram in Corollary \ref{cor-residue}.
Thus we see that $\cosp_X^m$ is bijective by a similar norm argument as in Step (\II) in the proof of Proposition \ref{prop:local-cond}.
This completes the proof of Proposition \ref{prop:local-cond2} and Theorem \ref{thm:local-cond'}.
\end{pf}

\par\medskip

By Proposition \ref{prop:local-cond2}\,(1) and \cite{BK} Corollary 3.8.4 for $V^m(r)=\H^m(X_{\ol K},\qp(r))$, we obtain the following corollary:
\begin{cor}\label{cor:exp}
The exponential map of Bloch-Kato induces an isomorphism
\[ \exp :  \H^m_{\dR}(X_{\nsp K}/K) \os{\simeq}\lra \H^1_{\!f}(K,V^m(r)) \]
for any $m \geqq 0$ and $r \geqq d$.
\end{cor}

\section{Comparison with Selmer groups,\, global case}\label{sect5}
Let $\pi_{X\!/\!B} : X \to B=\Spec(\fO)$ be as in \S\ref{sect1}.
In the rest of this paper, we always assume:
\begin{itemize}
\item
$\pi_{X\!/\!B}$ is {\it proper}, and the generic fiber $X_{\nsp K}$ is {\it geometrically connected over $K$}.
\item
$K$ is an algebraic number field, and $\fO$ is the integer ring of $K$, i.e., the case (G) of \S\ref{sect-limits}.
\end{itemize}
Put $V^m:=\H^{m}(X_{\ol K},\qp)$.
In this section, we compare $\H^m(X,\qp(r))$ with the Selmer group $\H^1_{\!f}(K,V^{m-1}(r))$, using the results of the previous section. See \cite{BK}\,\S5 for the definition of $\H^1_{\!f}(K,-)=\H^1_{\!f,B}(K,-)$.
For a place $v$ of $K$, we often write $K_v$ for the completion of $K$ at $v$.
For a finite place $v$ of $K$, we put $B_v:=\Spec(O_v)$ and $X_v:=X \times_B B_v$, where $O_v$ denotes the valuation ring of $K_v$.

\subsection{Fast computations}

\begin{prop}\label{thm:selmer}
Assume $r \geqq d$. Then{\rm:}
\begin{enumerate}
\item[{\rm(1)}]
$\H^q(B,\fH^m(X,\zp(r)))$ is finite in each of the following cases{\rm:}
\begin{itemize}
\item[{\rm(i)}] $m<0$  \qquad {\rm(ii)}\, $m>2(d-1)$ \qquad {\rm(iii)}\, $q \leqq 0$
 \qquad {\rm(iv)}\, $q>3$
\item[{\rm(v)}] $q=3$ \, and \,$(m,r) \ne (2(d-1),d)$
\end{itemize}
Consequently, the spectral sequence $\eqref{ss:leray2}$ degenerates at $E_2$-terms up to finite $p$-primary torsion.
\item[{\rm(2)}]
For any $m \geqq 0$, we have
\[ \H^1(B,\fH^m(X,\qp(r))) \cong \H^1_{\!f}(K,V^m(r)). \]
\end{enumerate}
\end{prop}

\begin{pf*}{Proof of Proposition \ref{thm:selmer}}
(1)\; We put
\[ \H^{q,m,r}:=\H^q(B,\fH^m(X,\zp(r))) \]
for simplicity.
The cases (i) and (ii) are clear by the definition of $\fH^m(X,\Tn(r))$ (see Definition \ref{def:mot}).
The case (iii) with $q<0$ follows from the fact that $\fH^m(X,\Tn(r))$ is concentrated in degrees $\geqq 0$
 (see Proposition \ref{prop:local-str}\,(3)). When $q=0$, the restriction map
\[ \H^{0,m,r} \lra \H^m(X_{\ol K},\zp(r))^{G_{\nsp\nsp K}} \]
is injective by Proposition \ref{prop:local-str}\,(1) and the last group is finite by \cite{d2}\,3.3.9. Hence $\H^{0,m,r}$ is finite. The case (iv) follows from the Artin-Verdier duality \cite{mazur}\,(2.4). Indeed, we have
\[ \H^q(B,\fH^m(X,\Tn(r))) \cong \Ext_B^q(R^{2(d-1)-m}\pi_{X/\!B*}\Tn(d-r),\Gm) \]
by \eqref{eq1-1}, and its dual
\[  \H^{3-q}_c(B,R^{2(d-1)-m}\pi_{X/\!B*}\Tn(d-r)) \]
is finite 2-torsion for any $n \geqq 1$ and $q>3$.
Finally we prove the case (v).
Fix a dense open subset $U \subset B[p^{-1}]$ such that $X_U \ra U$ is smooth (and proper).
Let $j$ be the open immersion $U \hra B$, and for each $v \in B$ let $\iota_v : v \hra B$ be the canonical map. There is an exact sequence
\[ \H^3(B,j_!j^*\fH^m(X,\zp(r))) \lra \H^{3,m,r} \lra \bigoplus_{v \in B \ssm U} \ \H^3(B_v,\fH^m(X_v,\zp(r))), \]
where we identified $\H^3(v,\iota_v^*\fH^m(X,\zp(r)))$ with $\H^3(B_v,\fH^m(X_v,\zp(r)))$ for each $v \in B \ssm U$ by Corollary \ref{cor:rigidity}.
The first term in this sequence is finite unless $(m,r)=(2(d-1),d)$ by the Artin-Verdier duality and a weight argument which is similar as for the case $q=0$. The last term is finite as well by Corollary \ref{cor4-fin}\,(2).
Thus $\H^{3,m,r}$ is finite in the case (v), which completes the proof of Proposition \ref{thm:selmer}\,(1).
\par
\medskip
(2)\; Let $S$ be a finite set of places of $K$ containing all places dividing $p$ or $\infty$, and all finite places where $X$ has bad reduction. Let $K_S$ be the maximal $S$-ramified extension of $K$ (i.e., the maximal Galois extension of $K$ which is unramified at every finite place of $K$ outside of $S$), and let $G_{\!S}$ be the Galois group $\Gal(K_S/K)$. To prove Proposition \ref{thm:selmer}\,(2), it is enough to check the following:
\begin{lem}\label{lem5-1}
There is an exact sequence of $\qp$-vector spaces
\begin{equation}\notag
0 \lra \H^{1,m,r} \otimes_{\zp} \qp \lra \H^1(G_{\!S},V^m(r)) \os{\res}{\lra} \bigoplus_{v \in S \cap B_0} \ \H^1_{\ovf}(K_v,V^m(r)).
\end{equation}
where $B_0$ denotes the set of closed points of $B$.
See \S\ref{sect3} for the definition of $\H^1_{\ovf}(K_v,-)$.
\end{lem}

\begin{pf}
Consider the localization long exact sequence of cohomology groups for each $n \geqq 1$
{\small
\begin{align*}
\dotsb & \to \H^q(B,\fH^m(X,\Tn(r))) \to \H^q(G_{\!S}, \H^m(X_{\ol K},\mu_{p^n}^{\otimes r})) \to \bigoplus_{v \in S \cap B_0} \ \H^{q+1}_v(B_v,\fH^m(X_v,\Tn(r))) \\
 &  \to \H^{q+1}(B,\fH^m(X,\Tn(r))) \to \dotsb,
\end{align*}
}where we have used the fact that $\fH^m(X,\Tn(r))|_{B \ssm S}$ is a locally constant sheaf on $B \ssm S$ associated with the $G_{\!S}$-module $\H^m(X_{\ol K},\mu_{p^n}^{\otimes r})$ (see Proposition \ref{ex:trace}\,(1)).
We have also used the isomorphisms
\[  \H^*_v(B,\fH^m(X,\Tn(r))) \cong \H^*_v(B_v,\fH^m(X_v,\Tn(r))) \qquad \hbox{($v \in S \ssm P_\infty$)} \]
obtained from \'etale excision and the rigidity of Corollary \ref{cor:rigidity}.
The groups in this long exact sequence are finite by Proposition \ref{prop1-1}.
Therefore we obtain the following long exact sequence by taking the projective limit with respect to $n \geqq 1$ and then $\otimes_{\zp}\qp$:
\begin{multline*}
 \dotsb \lra \bigoplus_{v \in S} \ \H^{q}_v(B_v,\fH^m(X_v,\qp(r))) \lra \H^{q,m,r} \otimes_{\zp} \qp \lra \H^q(G_{\!S},V^m(r)) \\
 \lra \bigoplus_{v \in S} \ \H^{q+1}_v(B_v,\fH^m(X_v,\qp(r))) \lra \dotsb.
\end{multline*}
Moreover we have
\[ \H^q_v(B_v,\fH^m(X_v,\qp(r))) \cong
 \begin{cases}
   0  & \quad \hbox{($q=1$)}\\
   \H^1_{\ovf}(K_v,V^m(r)) & \quad \hbox{($q=2$)} \\
   \H^2(K_v,V^m(r)) & \quad \hbox{($q=3$)}
 \end{cases} \]
by Theorem \ref{thm:local-cond}; the case $q=3$ will be useful later in the proof of Corollary \ref{cor5-4}\,(2) below. The assertion follows from these facts.
\end{pf}
\medskip\noindent
This completes the proof of Proposition \ref{thm:selmer}.
\end{pf*}

\begin{cor}\label{cor4-1}
For any $r \geqq d$, the spectral sequence \eqref{ss:leray3} degenerates at $E_2$, and we have
\[ \H^m(X,\qp(r)) \cong \begin{cases} \H^1_{\!f}(K,V^{m-1}(r)) \\
       \hspace{50pt} \oplus \, \H^2(B,\fH^{m-2}(X,\qp(r))) \quad & \hbox{{\rm(}$1 \leqq m \leqq 2d-1${\rm)}} \\
 \qp & \hbox{{\rm(}$(m,r) = (2d+1,d)${\rm)}}\\
 0 & \hbox{{\rm(}otherwise{\rm).}}
          \end{cases} \]
\end{cor}
\noindent
See Corollary \ref{cor-duality}\,(2) for the isomorphism $\H^{2d+1}(X,\qp(d)) \cong \qp$.
We will prove that $\H^2(B,\fH^m(X,\qp(r)))=0$ for any $(m,r)$ with $r \geqq d$, in Theorem \ref{thm5-1} below.
The following consequence of Proposition \ref{thm:selmer}\,(2) is a global analogue of Corollary \ref{cor4-fin}\,(1), which will be useful later.
\begin{cor}\label{cor5-3}
For any $r \geqq d$, there exists a natural map
\[ H^1_{\!f}(K,V^m(r)) \lra \H^1(B,\fH^m(X,\QpZp(r))) \]
which fits into a commutative diagram
\[\xymatrix{
H^1_{\!f}(K,V^m(r)) \ar@{.>}[d]\ar[rd]^{\text{{\rm(}natural map{\rm)}}} \\
\H^1(B,\fH^m(X,\QpZp(r))) \; \ar@{^{(}->}[r] & \H^1(K,\H^m(X_{\ol K},\QpZp(r))).
}\]
See Proposition \ref{prop:local-str}\,{\rm(}1{\rm)} for the injectivity of the bottom arrow.
\end{cor}

\begin{rem}\label{rem4-1}
{\rm
For any $s \leqq 0$, one can easily check the following canonical isomorphism by \eqref{isom:local-cond2}, \eqref{local1}
and similar arguments as for the proof of Proposition \ref{thm:selmer}:
\[ \H^1(B,R^m\pi_{X/\!B*}\qp(s)) \cong \H^1_{\!f}(K,V^m(s)). \]
}
\end{rem}

\subsection{A global finiteness of \'etale cohomology}\label{sect5-2}
In this subsection, we prove the following vanishing and finiteness result:
\begin{thm}\label{thm5-1}
For any $m \geqq 0$ and $r \geqq d$, we have
\[ \H^2(B,\fH^m(X,\qp(r))) = 0, \]
and the groups $\H^2(B,\fH^m(X,\zp(r)))$ and $\H^2(B,\fH^m(X,\QpZp(r)))$ are finite.
\end{thm}

\noindent
As a direct consequence of this theorem and Corollary \ref{cor4-1}, we obtain:

\begin{cor}\label{cor5-1}
For any $m \geqq 0$ and $r \geqq d$ with $(m,r) \ne (2d+1,d)$, we have
\[ \H^m(X,\qp(r)) \cong \H^1_{\!f}(K,V^{m-1}(r)). \]
\end{cor}

\noindent
On the other hand, Theorem \ref{thm5-1} and Remark \ref{rem4-1} imply the following vanishing result by the Artin-Verdier duality:
\begin{cor}\label{cor5-2}
For any $m \geqq 0$ and $s \leqq 0$, we have $\H^1_{\!f}(K,\H^m(X_{\ol K},\qp(s))) = 0$.
\end{cor}

\medskip

\begin{pf*}{Proof of Theorem \ref{thm5-1}}
By Theorem \ref{lem:fg}, it is enough to show that $\H^2(B,\fH^m(X,\QpZp(r)))$ is finite.
When $(m,r)=(2(d-1),d)$, we have
\[ \H^2(B,\fH^{2(d-1)}(X,\QpZp(d))) \os{\eqref{isom:trace}}\cong \H^2(B,\QpZp(1)) \cong \Br(\OK)\{p\}, \]
by the finiteness of $\Pic(\OK)$, and $\Br(\OK)$ is finite 2-torsion by the classical Hasse principle for Brauer groups, which implies the finiteness in question.
\par
In what follows, we assume $(m,r) \ne (2(d-1),d)$ and consider the following commutative diagram with exact rows, where both rows are obtained from localization sequences of \'etale cohomology, and the coefficients $\fH^m(X,\QpZp(r))$ (resp.\ $\fH^m(X_v,\QpZp(r))$) in the upper row (resp.\ the lower row) are omitted:
{\small
\begin{equation}\notag
\xymatrix{
\H^1(K) \ar[r]\ar[d]_{\alpha} &
 \displaystyle \bigoplus_{v \in B_0} \H^2_v(B) \ar[r] \ar[d]_{\delta}  &
 \H^2(B) \ar[r] \ar[d]_{\beta} &
 \H^2(K) \ar[r] \ar[d]_{\gamma} &
 \displaystyle \bigoplus_{v \in B_0} \H^3_v(B) \ar[d]_{\delta} \\
\displaystyle \bigoplus_{v \in B_0} \H^1_{\ovf}(K_v) \ar[r]^{(*)} &
 \displaystyle \bigoplus_{v \in B_0} \H^2_v(B_v) \ar[r] &
\displaystyle \bigoplus_{v \in B_0} \H^2(B_v) \ar[r] &
\displaystyle \bigoplus_{v \in B_0} \H^2(K_v) \ar[r] &
 \displaystyle \bigoplus_{v \in B_0} \H^3_v(B_v).
}\end{equation}
}Here we put
\[ \H^1_{\ovf}(K_v):=\Coker \big(\H^1_{\!f}(K_v,V^m(r)) \to \H^1(K_v,\H^m(X_{\ol {K_v}},\QpZp(r)))\big) \]
for each $v \in B_0$ (note also Proposition \ref{ex:trace}\,(1)), and
used Corollary \ref{cor4-fin}\,(1) to verify the existence of the bottom left arrow $(*)$.
The arrows $\delta$ are bijective by \'etale excision and the rigidity (Corollary \ref{cor:rigidity}).
The arrow $\gamma$ has finite kernel and cokernel by the Hasse principle of Jannsen \cite{J} p.\ 337, Theorem 3\,(c).
The arrow $\alpha$ has finite cokernel by \cite{BK} Proposition 5.14\,(ii).
Hence $\beta$ is bijective up to finite groups.
Finally, $\H^2(B_v,\fH^m(X_v,\QpZp(r)))$ is finite for all $v \in B_0$ by Corollary \ref{cor4-fin}\,(2),
and zero for any $v \in (B[p^{-1}])_0$ at which $X$ has good reduction by Lemma \ref{lem4-1}.
Thus $\H^2(B,\fH^m(X,\QpZp(r)))$ is finite.
\end{pf*}

\begin{rem}
{\rm
\begin{enumerate}
\item[(1)]
By Theorem \ref{thm5-1} for $r=d=2$ and $m=1$ and Lemma \ref{lem6-2}\,(3) below, Bloch's conjecture (\cite{B1} Remark 1.24) for a projective smooth curve $C$ over $K$ is reduced to a variant of Bass' conjecture (cf.\ \cite{Ba}) that the motivic cohomology $\H^3_{\cM}(X,\bZ(2))$ is finitely generated for a proper regular model $X\nsp /\nsp B$ of $C$.
\item[(2)]
Corollary \ref{cor5-2} removes an assumption of a result of Morin \cite{Mo} Theorem 1.5\,(3).
\end{enumerate}
}
\end{rem}

The following corollary of Theorem \ref{thm5-1} follows from a similar argument as for the proof of Lemma \ref{lem5-1} (see also \cite{J} p.\ 349, Question 2):
\begin{cor}\label{cor5-4}
Assume $r \geqq d$, and let $S$ and $G_S$ be as in the proof of Theorem \ref{thm:selmer}\,{\rm(}2{\rm)}.
Then{\rm:}
\begin{enumerate}
\item[{\rm(1)}]
For any $m$, the following map is surjective{\rm:}
\[ \H^1(G_{\!S},V^m(r)) \lra \bigoplus_{v \in S} \ \H^1_{\ovf}(K_v,V^m(r)). \]
\item[{\rm(2)}]
The restriction map
\[ \H^2(G_{\!S},V^m(r)) \lra  \bigoplus_{v \in S} \ \H^2(K_v,V^m(r)) \]
is bijective for any $(m,r) \ne (2(d-1),d)$ and injective for $(m,r)=(2(d-1),d)$.
In particular, if $r>d$ or $X_{\nsp K}$ has potentially good reduction at all finite places of $K$, then
\[ \H^2(G_{\!S},V^m(r))=0 \quad \hbox{ for any $(m,r) \ne (2(d-1),d)$}. \]
\end{enumerate}
\end{cor}
See Remark \ref{rem-log} for a remark on our log-smoothness assumption.

\section{$\bs{p}$-adic Abel-Jacobi mappings \quad ($\bs{d=2}$)}\label{sect6}
The setting remains as in \S\ref{sect5}.
From this section on, we assume further that $d=2$.

\subsection{Cycle class maps}\label{sect6-1}
See \S\ref{sect1-2} for the definition of the motivic complex $\bZ(r)$ on $(\Et/X)_\Zar$.
We regard $\bZ(r)$ as a complex on $X_\Zar$ by restriction of topology.
We define the motivic cohomology of $X$ as
\[\H^m_{\cM}(X,\bZ(r)) := \H^m_{\Zar}(X,\bZ(r)), \]
and define the motivic cohomology with $\Ln(=\bZ/p^n\bZ)$-coefficients as
\[ \H^m_{\cM}(X,\Ln(r)) := \H^m_{\Zar}(X,\bZ(r) \otimes \Ln) \qquad \hbox{($n \geqq 1$)}. \]
In this paper we do not consider the motivic complex $\bZ(r)$ on $X_\et$, mainly because it is not necessarily compared with $\fT_n(r)$ directly for the lack of the Gersten resolution for $\bZ(r) \otimes \Ln$ unless $X$ is smooth over $B$, cf.\ \cite{Sa0} Conjecture 1.4.1, \cite{Sa1} Remark 7.2, \cite{Z} Conjecture 2.2, Theorem 4.8, \cite{Ge0} Theorem 1.2\,(5).
\begin{lem}\label{lem6-2}
\begin{enumerate}
\item[{\rm(1)}]
We have
\[ \H^m_{\cM}(X,\bZ(r)) \cong \begin{cases}
 \H^m_{\cM}(K(X),\bZ(2)) \quad \phantom{\big|_|} & \hbox{\rm ($m \leqq 1$,\, $r=2$)} \\
 0 & \hbox{\rm ($m > r+2$)}
\end{cases} \]
where $K(X)$ denotes the function field of $X$.
\item[{\rm(2)}]
$\H^m_{\cM}(X,\bZ(2))$ is isomorphic to the cohomology at deree $m-2$ of the Gersten complex of Milnor \tK-groups
\begin{align*}
&\tK^M_2(K(X)) \lra \bigoplus_{x \in X^1} \ \kappa(x)^\times \lra \bigoplus_{x \in X^2} \ \bZ \\
&\quad \text{\rm (deg 0)} \phantom{\lra} \qquad \text{\rm (deg 1)} \phantom{\lra} \qquad  \text{\rm (deg 2)}
\end{align*}
for any $m \geqq 2$.
In particular, we have $\H^4_{\cM}(X,\bZ(2)) \cong \CH_0(X)$, the Chow group of $0$-cycles modulo rational equivalence.
\item[{\rm(3)}]
Assume that $r \geqq 2$, and that $p \geqq 3$ or $B(\bR)=\emptyset$.
Then the cycle class map {\rm(}see \S\ref{sect1-2}{\rm)}
\[ \cl_{\Ln}^{m,r} : \H^m_{\cM}(X,\Ln(r)) \lra \H^m(X,\Tn(r)) \]
is bijective for any $m \in \bZ$ with $(m,r) \ne (5,2)$ and any $n \geqq 1$.
Consequently, there exists a short exact sequence
\[ 0 \lra \H^m_{\cM}(X,\bZ(r))/p^n \lra \H^m(X,\Tn(r)) \lra  {}_{p^n}\H^{m+1}_{\cM}(X,\bZ(r)) \lra 0 \]
for the same $(m,n)$, where for an abelian group $M$, $_{p^n}M$ {\rm(}resp.\ $M/p^n${\rm)} denotes the kernel {\rm(}resp.\ cokernel{\rm)} of the map $M \os{\times p^n}\lra M$.
\end{enumerate}
\end{lem}

\begin{pf}
There exists a coniveau spectral sequence
\begin{equation}\label{eq6+0}
 E_1^{a,m}=\bigoplus_{x \in X^a} \H^{m-a}_{\cM}(x,\bZ(r-a)) \Lra \H^{a+m}_{\cM}(X,\bZ(r))
\end{equation}
by \cite{Ge} Proposition 2.1, whose $E_1^{a,m}$-terms are zero in each of the following cases for the reason of the dimension of cycles and the codimension of points:
\begin{enumerate}
\item[$\circ$]
$m > r$ \qquad $\circ$ \, $a < 0$ \qquad $\circ$\, $a > 2$
\qquad $\circ$\, $m < a=r$ \qquad $\circ$\, $m \leqq a=r-1$
\end{enumerate}
See \cite{B2} Theorem 6.1 for the vanishing in the last case.
The assertions (1) and (2) follow from these facts and the Nesterenko-Suslin-Totaro theorem
\[ \H^q_{\cM}(\Spec(F),\bZ(q)) \cong \tK^M_q(F) \]
for any field $F$ and any $q \geqq 0$, see \cite{NS}, \cite{To}.
\par
To prove the assertion (3), we consider a coniveau spectral sequence analogous to \eqref{eq6+0}
\begin{equation}\label{eq6+1}
E_1^{a,m}=\bigoplus_{x \in X^a} \H^{m-a}_{\cM}(x,\Ln(r-a)) \Lra \H^{a+m}_{\cM}(X,\Ln(r)),
\end{equation}
whose $E_1^{a,m}$-terms are zero in each of the following cases:
\begin{enumerate}
\item[$\circ$]
$m > r$ \qquad $\circ$\, $a < 0$ \qquad $\circ$\, $a > 2$
\end{enumerate}
On the other hand, since $r \geqq 2$, there is a coniveau spectral sequence of \'etale cohomology (see \cite{JSS} (5.10.1))
\begin{equation}\label{eq6+2}
E_1^{a,m}=\bigoplus_{x \in X^a} \H^{m-a}(x,\Ln(r-a)) \Lra \H^{a+m}(X,\Tn(r)),
\end{equation}
where the coefficients $\Ln(s)$ ($s \in \bZ$) on the points are those in \eqref{eq-RD}.
The $E_1^{a,m}$-terms of \eqref{eq6+2} are zero in each of the following cases:
\begin{enumerate}
\item[$\circ$]
$m > 3$ \qquad $\circ$\, $m < a$ \qquad $\circ$\, $a < 0$ \qquad $\circ$\, $a > 2$.
\end{enumerate}
Here we have used the well-known fact that the $\cd_p(\kappa(x)) = 3-a$ with for any $a \geqq 0$ and $x \in X^a$ (see e.g., \cite{T} Theorem 3.1, \cite{se} Chapter \II, \S4.2 Proposition 11).
There is a map of spectral sequences from \eqref{eq6+1} to \eqref{eq6+2} induced by cycle class maps of motivic cohomology groups by the commutative diagrams \eqref{eq-cd-cycle} and \eqref{eq-cd-cycle2} in \S\ref{sect1-2}. The cycle class map
\[ \H^{m-a}_{\cM}(x,\Ln(r-a)) \lra  \H^{m-a}(x,\Ln(r-a)) \]
is bijective for any $a \geqq 0$, any point $x \in X^a$ and any $m \leqq r$ by Rost-Voevodsky \cite{V1}, \cite{V2} Theorem 6.16 and Geisser-Levine \cite{GL2} Theorem 7.5 (resp.\ Bloch-Gabber-Kato \cite{BK1} Theorem 2.1 and Geisser-Levine \cite{GL1} Theorem 1.1), when $\ch(x) \ne p$ (resp.\ when $\ch(x) = p$).
If $r \geqq 3$, then the map $\cl_{\Ln}^{m,r}$ in question is bijective by these facts.
As for the case $r=2$, it remains to check that the $E_\infty^{a,3}$-terms of \eqref{eq6+2} are zero for $a=0$ and $1$,
which is a consequence of Kato's Hasse principle \cite{KCT} p.\ 145, Corollary.
\end{pf}
\begin{rem}
{\rm
If we assume the Beilinson-Soul\'e vanishing conjecture (\cite{Sou} p. 501, Conjecture) for points of $X$,
 then we would have
\[ \H^m_{\cM}(X,\bZ(r)) \cong
 \begin {cases} \H^1_{\cM}(K(X),\bZ(r)) \;\; \phantom{\big|_|} & \hbox{($m=1$)} \\
 0 & \hbox{($m \leqq 0$)}
\end{cases} \]
up to small torsion for any $r \geqq 2$, by the same arguments as in the proof of Lemma \ref{lem6-2}\,(1).
}
\end{rem}

\subsection{$\bs{p}$-adic Abel-Jacobi mappings and finiteness results}

Let $r$ be an integer with $r \geqq 2$. We define a $p$-adic cycle class map
\[ \cl_p^{m,r} : \H^m_{\cM}(X,\bZ(r)) \,\wh{\otimes}\, \zp \lra \H^m(X,\zp(r)) \]
as the projective limit with respect to $n \geqq 1$ of the cycle class map
\[ \cl_{/p^n}^{m,r} : \H^m_{\cM}(X,\bZ(r))/p^n \lra
 \H^m_{\cM}(X,\Ln(r)) \us{\simeq}{\os{\cl_{\Ln}^m}\lra} \H^m(X,\Tn(r)). \]
See Lemma \ref{lem6-2}\,(3) for the isomorphism $\cl_{\Ln}^{m,r}$.
Since $X_{\ol K}$ is a curve, $\H^m(X_{\ol{K}},\zp(r))$ is torsion-free, and
\begin{equation}\label{eq6-0}
 \H^0(B,\fH^m(X,\zp(r))) \subset \H^m(X_{\ol{K}},\zp(r))^{G_{\nsp\nsp K}} = 0
\end{equation}
by Proposition \ref{prop:local-str}\,(1) and for the reason of weights.
We define a $p$-adic Abel-Jacobi mapping
\[ \aj_p^{m,r} : \H^m_{\cM}(X,\bZ(r)) \,\wh{\otimes}\, \zp \lra \H^1(B,\fH^{m-1}(X,\zp(r))) \]
as the map induced by $\cl_p^{m,r}$ and an edge map of the spectral sequence \eqref{ss:leray2}:
\begin{equation}\label{eq6-1}
E^{a,b}_2=\H^a(B,\fH^b(X,\zp(r))) \Lra \H^{a+b}(X,\zp(r)).
\end{equation}
We first observe the following straight-forward remarks:

\begin{prop}\label{lem6-1}
Let $m$ and $r$ be integers with $r \geqq 2$ and $(m,r) \ne (5,2)$.
Assume that $p \geqq 3$ or $B(\bR)=\emptyset$.
Then the following five conditions are equivalent to one another{\rm :}
\begin{enumerate}
\item[{\rm (i)}]
$\aj_p^{m,r}$ has finite cokernel.
\item[{\rm (ii)}]
$\cl_p^{m,r}$ has finite cokernel.
\item[{\rm (iii)}]
$\cl_p^{m,r}$ is surjective.
\item[{\rm (iv)}]
$\H^{m+1}_{\cM}(X,\bZ(r))\{p\}$ is finite.
\item[{\rm (v)}]
$\H^{m+1}_{\cM}(X,\bZ(r))_{p\text{-}\Div}$ is uniquely $p$-divisible.
\end{enumerate}
Moreover if $m \leqq 1$, these conditions are equivalent to
\begin{enumerate}
\item[{\rm (i$'$)}]
$\aj_p^{m,r}$ is surjective.
\end{enumerate}
\end{prop}
\begin{pf}
The term $E_2^{a,m}$ of \eqref{eq6-1} is finite for any $a \geqq 2$ by Theorems \ref{thm:selmer}\,(1) and \ref{thm5-1}, which shows (iii)\,$\Ra$\,(i). The assertion (i)\,$\Ra$\,(ii) is a consequence of the following fact (a), and the assertion (ii)\,$\Ra$\,(iii) is a consequence of the fact (b) below, where $T_p$ denotes the $p$-Tate module:
\begin{enumerate}
\item[(a)] {\it The canonical map
\[ \H^m(X,\zp(r)) \lra \H^1(B,\fH^{m-1}(X,\zp(r))) \]
has finite kernel by Theorem \ref{thm5-1}.}
\item[(b)] {\it 
By taking the projective limit with respect to $n \geqq 1$ of the short exact sequence of Lemma \ref{lem6-2}\,{\rm(}3{\rm)}, we have $\Coker(\cl_p^{m,r}) \cong T_p(\H^{m+1}_{\cM}(X,\bZ(r)))$,
which are torsion-free.}
\end{enumerate}
We next prove (iii) $\Leftrightarrow$ (iv).
Indeed, by taking the inductive limit with respect to $n \geqq 1$ of the short exact sequence of Lemma \ref{lem6-2}\,(3),
we get an exact sequence
\begin{equation}\label{eq6-1+}
 0 \to \H^m_{\cM}(X,\bZ(r))\otimes \QpZp \to \H^m(X,\QpZp(r)) \to  \H^{m+1}_{\cM}(X,\bZ(r))\{p\} \to 0,
\end{equation}
which imply that $\H^{m+1}_{\cM}(X,\bZ(r))\{p\}$ is cofinitely generated over $\zp$,
 see Theorem \ref{lem:fg}\,(2). Hence
\[ \hbox{(iii)} \; \os{\text{(b)}}\Llra \; T_p(\H^{m+1}_{\cM}(X,\bZ(r))) = 0 \; \Llra \; \hbox{(iv)}. \]
The assertion (iv)\,$\Ra$\,(v) is obvious, and the assertion (v)\,$\Ra$\,(iv) also follows from the fact that $\H^{m+1}_{\cM}(X,\bZ(r))\{p\}$ is cofinitely generated over $\zp$.
Finally, if $m \leqq 1$, the canonical map in (a) is bijective by \eqref{eq6-0}, which shows that (iii) is equivalent to (i$'$).
\end{pf}
\par\medskip

The following lemma will be useful in what follows.

\begin{lem}\label{lem6+0}
Assume that $p \geqq 3$ or $B(\bR)=\emptyset$. Then{\rm:}
\begin{enumerate}
\item[{\rm (1)}]
$\cl_p^{m,r}$ is injective for any $m \in \bZ$ and $r \geqq 2$.
\item[{\rm (2)}]
We have $\H^5(X,\Tn(2)) \cong \Ln$ for any $n \geqq 1$, and $\H^m(X,\Tn(r))=0$ for any $m \geqq 5$, $r \geqq 2$ and $n \geqq 1$ with $(m,r) \ne (5,2)$.
\end{enumerate}
\end{lem}
\begin{pf}
The assertion (1) follows from Lemma \ref{lem6-2}\,(1) and (3).
The assertions in (2) follow from the duality (see Corollary \ref{cor-duality}\,(2), \eqref{eq-RHom})
\[ \H^m(X,\Tn(r)) \cong \H^{5-m}(X,\Tn(2-r))^*. \]
The details are straight-forward and left to the reader.
\end{pf}
\par\medskip

The following result gives an extension of the vanishing assertion in Lemma \ref{lem6-2}\,(1):

\begin{prop}\label{cor6-1}
Assume that $p \geqq 3$ or $B(\bR)=\emptyset$. Then
\[ \H^m_{\cM}(X,\bZ(r))\{p\}, \quad \H^m_{\cM}(X,\bZ(r))\,\wh{\otimes}\,\zp \quad \hbox{and} \quad \H^m(X,\zp(r)) \]
are zero for any $m \geqq 5$ and $r \geqq 3$.
In particular, $\H^m_{\cM}(X,\bZ(r))$ is uniquely $p$-divisible for the same $(m,r)$.
\end{prop}
\begin{pf}
We have $\H^m(X,\zp(r))=0$ by Lemma \ref{lem6+0}\,(2),
so $\H^m_{\cM}(X,\bZ(r))\,\wh{\otimes}\,\zp=0$ by Lemma \ref{lem6+0}\,(1).
To show that $\H^m_{\cM}(X,\bZ(r))\{p\}=0$, we use the surjectivity of the boundary map
\[ \H^{m-1}(X,\QpZp(r)) \twoheadrightarrow  \H^m_{\cM}(X,\bZ(r))\{p\} \]
of \eqref{eq6-1+}.
By Lemma \ref{lem6+0}\,(2), we have $\H^{m-1}(X,\QpZp(r))=0$ for any $m \geqq 6$, which implies that $\H^m_{\cM}(X,\bZ(r))\{p\}$ is zero for any $m \geqq 6$ by \eqref{eq6-1+}.
As for the case $m=5$, we have $\H^4(X,\QpZp(r))=0$.
Indeed, it is finite by Corollary \ref{cor5-1}, and $p$-divisible by the exact sequence (see Proposition \ref{prop:bock})
\[ \dotsb \lra \H^4(X,\QpZp(r)) \os{\times p}\lra \H^4(X,\QpZp(r)) \lra \H^5(X,\fT_1(r)) \lra \dotsb \]
and Lemma \ref{lem6+0}\,(2). Thus $\H^5_{\cM}(X,\bZ(r))\{p\}$ is zero.
\end{pf}

\begin{prop}\label{cor6-2}
Assume that $p \geqq 3$ or $B(\bR)=\emptyset$. Then for any $r \geqq 3$, we have
\[ \H^4_{\cM}(X,\bZ(r))\{p\} \cong \H^4_{\cM}(X,\bZ(r))\,\wh{\otimes}\,\zp \os{\simeq}{\us{\cl_p^{4,r}}\lra} \H^4(X,\zp(r)),\]
which are all finite.
\end{prop}
\begin{pf}
The cycle class map $\cl_p^{4,r}$ is injective by Lemma \ref{lem6+0}\,(1), and surjective by Proposition \ref{lem6-1}\,(iv)\,$\Ra$\,(iii) and the vanishing of $\H^5_{\cM}(X,\bZ(r))\{p\}$ in Proposition \ref{cor6-1}.
The finiteness of $\H^4(X,\zp(r))$ follows from Corollary \ref{cor5-1}.
\par
We next prove that $\H^4_{\cM}(X,\bZ(r))\{p\}$ is finite.
By Proposition \ref{lem6-1}\,(i)\,$\Ra$\,(iv), it is enough to check that the map
\[ \aj_p^{3,r} : \H^3_{\cM}(X,\bZ(r))\, \wh{\otimes} \, \zp \lra \H^1(B,\fH^2(X,\zp(r)))
 \cong \H^1(B[p^{-1}],\zp(r-1)) \]
has finite cokernel, where the last isomorphism follows from Proposition \ref{ex:trace}\,(2) and Lemma \ref{lem1-0}\,(2) for $B$. Take a finite morphism $f : B' \to X$ such that $B'$ is regular $1$-dimensional, and such that the composite $g : B' \to X \to B$ is finite flat.
To check the finiteness of $\Coker(\aj_p^{3,r})$, we construct a Chern character \eqref{eq-Chern} below, using the fact due to Quillen (\cite{Q1} Theorem 8) that the algebraic \tK-group $\tK_i(k)$ of a finite field $k$ is finite for any $i \geqq 1$.
By this fact and the localization sequence of algebraic \tK-groups (\cite{Q2} p.\ 113, Corollary of Theorem 5), we have
\begin{equation}\label{eq-Quillen}
 \tK_i(B') \otimes \bQ \cong \tK_i(L) \otimes \bQ \qquad \hbox{ for \, any \, } i \geqq 2,
\end{equation}
where $L$ denotes the function field of $B'$. On the other hand, the Chern character
\[  \ch_{F\nsp,\sp i}^{\mathscr {M}} : \tK_i(F) \otimes \bQ \lra \bigoplus_{j \geqq 0} \ \H^{2j-i}_{\cM}(F,\bZ(j))\otimes \bQ \]
is bijective for any field $F$ by Bloch \cite{B2} Theorem 9.1. Applying this fact to the closed points of $B'$ and Levine's localization \cite{Le} Theorem 1.7 to $B'$, we obtain
\begin{equation}\label{eq-Levine}
 \H^{2j-i}_{\cM}(B',\bZ(j)) \otimes \bQ \cong \H^{2j-i}_{\cM}(L,\bZ(j)) \otimes \bQ  \qquad \hbox{ for \, any \, } i \geqq 2.
\end{equation}
By \eqref{eq-Quillen} and \eqref{eq-Levine}, the Chern character $\ch_{L\nsp,\sp 2r-3}^{\mathscr {M}}$ defines a Chern character
\begin{equation}\label{eq-Chern}
 \ch_{B'\nsp,\sp 2r-3}^{\mathscr {M}\nsp,\sp r-1} : \tK_{2r-3}(B')\otimes \bQ \lra \H^1_{\cM}(B',\bZ(r-1))\otimes \bQ \qquad \hbox{($r \geqq 3$)},
\end{equation}
which fits into the following commutative diagram:
\[\xymatrix{
\tK_{2r-3}(B') \otimes \qp \ar[d]_{\ch_{B'\!,\sp 2r-3}^{\mathscr {M}\nsp,\sp r-1}\otimes \id_{\qp}} \ar[rd]^{\ch_{B'\!,\sp 2r-3}^{\et\nsp,\sp r-1}} \\
\H^1_{\cM}(B',\bZ(r-1)) \otimes \qp \ar[d]_{f_*} \ar[r] & \H^1(B'[p^{-1}],\qp(r-1)) \ar@{->>}[d]^{g_*} \\
\H^3_{\cM}(X,\bZ(r)) \otimes \qp \ar[r] & \H^1(B[p^{-1}],\qp(r-1)).
}\]
Here $\ch_{B'\!,\sp 2r-3}^{r-1\nsp,\sp\et}$ denotes the \'etale Chern character, and the middle and the bottom horizontal arrows are the $\qp$-linear extension of the following composite maps, respectively:
\begin{align*}
&\H^1_{\cM}(B',\bZ(r-1)) \lra \H^1_{\cM}(B',\bZ(r-1))\, \wh{\otimes} \, \zp
 \os{\cl_p^{1,r-1}}\lra \H^1(B'[p^{-1}],\zp(r-1)) \\
&\H^3_{\cM}(X,\bZ(r)) \lra \H^3_{\cM}(X,\bZ(r))\, \wh{\otimes} \, \zp \os{\aj_p^{3,r}}\lra \H^1(B[p^{-1}],\zp(r-1)),
\end{align*}
The arrow $g_*$ is surjective by a standard norm argument.
Now the finiteness of $\Coker(\aj_p^{3,r})$ in question follows from the surjectivity of $\ch_{B'\!,\sp 2r-3}^{\et\nsp,\sp r-1}$ (\cite{So2} Theorem 1, \cite{Ka} Theorem 5.3).
Thus $\H^4_{\cM}(X,\bZ(r))\{p\}$ is finite.
\par
Finally, the natural map $\H^4_{\cM}(X,\bZ(r))\{p\} \to \H^4_{\cM}(X,\bZ(r))\,\wh{\otimes}\,\zp$ is injective by the finiteness of $\H^4_{\cM}(X,\bZ(r))\{p\}$. To show the surjectivity of this map, consider the following commutative triangle:
\[\xymatrix{
\H^3(X,\QpZp(r)) \ar@{->>}[r]^{\delta} \ar[rd]_{\delta'} & \H^4_{\cM}(X,\bZ(r))\{p\} \ar[d]^{\cl_p^{4,r}|_\tors} \\
& \H^4(X,\zp(r)),
}\]
where the arrow $\delta$ denotes the boundary map of \eqref{eq6-1+}, and the arrow $\delta'$ denotes the boundary map of the long exact sequence obtained from Proposition \ref{prop:bock}
\[ \dotsb \to \H^3(X,\qp(r)) \to \H^3(X,\QpZp(r)) \os{\delta'}\lra \H^4(X,\zp(r)) \to \H^4(X,\qp(r)) \to \dotsb. \]
The arrow $\cl_p^{4,r}|_\tors$ means the restriction of $\cl_p^{4,r}$ to $\H^4_{\cM}(X,\bZ(r))\{p\}$.
Since $\delta'$ is surjective by the finiteness of $\H^4(X,\zp(r))$, $\cl_p^{4,r}|_\tors$ is surjective as well, which completes the proof.
\end{pf}

\pmn
The following corollary is a summary of known facts and our results on $\cl_p^{m,r}$ and $\aj_p^{m,r}$:

\begin{cor}\label{lem6-0}
Let $r$ be an integer with $r \geqq 2$, and assume that $p \geqq 3$ or $B(\bR)=\emptyset$.
Then{\rm:}
\begin{enumerate}
\item[{\rm (0)}]
$\H^m_{\cM}(X,\bZ(r))$ is uniquely $p$-divisible for any $m \leqq 0$ and any $m \geqq 5$, and zero for any $m > r+2$.
\item[{\rm (1)}]
$\cl_p^{1,r}$ and $\aj_p^{1,r}$ are injective.
\item[{\rm (2)}]
$\cl_p^{2,r}$ is injective, and $\aj_p^{2,r}$ has finite kernel.
\item[{\rm (3)}]
$\cl_p^{3,r}$ is bijective, and $\aj_p^{3,r}$ has finite kernel and cokernel.
\item[{\rm (4)}]
$\cl_p^{4,r}$ is bijective, and $\H^4_{\cM}(X,\bZ(r))\{p\}$ is finite.
Moreover, we have $\H^4_{\cM}(X,\bZ(r))\{p\} \cong \H^4_{\cM}(X,\bZ(r))\,\wh{\otimes}\,\zp$, and $\aj_p^{4,r}$ is zero.
\end{enumerate}
\end{cor}
\begin{pf}
The assertion (0) for $m \leqq 0$ follows from Lemmas \ref{lem6-2}\,(3) (for $m <0$) and \ref{lem6+0}\,(1)  (for $m=0$) and the vanishing of $\H^m(X,\fT_n(r))$ for $m < 0$ and $\H^0(X,\zp(r))$.
See Lemma \ref{lem6-2}\,(1) and Proposition \ref{cor6-1} for the other claims in (0).
The injectivity of $\cl_p^{m,r}$ in (1)--(4) is nothing other than Lemma \ref{lem6+0}\,(1),
and the finiteness of $\Ker(\aj_p^{m,r})$ in (2)--(4) follows from (a) in the proof of Proposition \ref{lem6-1}.
The injectivity of $\aj_p^{1,r}$ in (1) is that of $\cl_p^{1,r}$.
By Proposition \ref{lem6-1}, the surjectivity of $\cl_p^{m,r}$ and the finiteness of $\Coker(\aj_p^{m,r})$ are both equivalent to the finiteness of $\H^{m+1}_{\cM}(X,\bZ(r))\{p\}$.
This last finiteness for the case $m=4$ has been mentioned in (0);
the case $(m,r)=(3,2)$ is due to Bloch \cite{B1}, Kato-Saito \cite{KS}, see also Lemma \ref{lem6-2}\,(2);
the case $m=3$ and $r \geqq 3$ is a consequence of Proposition \ref{cor6-2}.
Finally, $\aj_p^{4,r}$ is zero for any $r \geqq 2$, because $\H^1(B,\fH^3(X,\zp(r)))=0$ by \eqref{filter2}.
\end{pf}

\subsection{$\bs{p}$-Tate-Shafarevich groups}
Let $r$ be an integer with $r \geqq 2$.
We put $T^m:=\H^m(X_{\ol K},\zp)$, $V^m:=T^m \otimes_{\zp} \qp$ and
\begin{align*}
\H^1_{\ovf}(K,T^m \otimes \QpZp(r))
 :=\frac{\,\H^1(K,T^m \otimes \QpZp(r))\,}{\,\text{Image of }\, \H^1_{\!f}(K,V^m(r))\,}.
\end{align*}
Note that $T^m \otimes \QpZp \cong \H^m(X_{\ol K},\QpZp)$, because $X_{\!K}$ is a curve by assumption.
Let $P$ (resp.\ $P_\infty$) be the set of all places of $K$ (resp.\ all infinite places of $K$).
We often identify a finite place of $K$ with a closed point of $B$. For each $v \in P$, we put
\begin{align*}
\H^1_{\ovf}(K_v,T^m \otimes \QpZp(r)) :=
\dfrac{\,\H^1(K_v,T^m \otimes \QpZp(r))\,}{\,\text{Image of }\, \H^1_{\!f}(K_v,V^m(r))\,}
\end{align*}
where $\H^1_{\!f}(K_v,V^m(r))$ means zero for any $v \in P_\infty$.
This group for $v \in B_0$ has been used in the proof of Theorem \ref{thm5-1}.
For $m \geqq 0$ and $r \geqq 2$ with $(m,r)\ne(2,2)$, the natural map
\begin{equation}\label{eq6+3}
 \alpha^{m,r} : \H^1_{\ovf}(K,T^m \otimes \QpZp(r))
 \lra \bigoplus_{v \in P} \ \H^1_{\ovf}(K_v,T^m \otimes \QpZp(r))
\end{equation}
has finite kernel and cokernel, and we have
\begin{equation}\label{eq6-3}
 \Coker(\alpha^{m,r}) \cong (T^{2-m}\otimes\QpZp(2-r))^{G_{\nsp\nsp K}})^*
\end{equation}
by \cite{BK}\,Proposition 5.14\,(i), (ii).
The $p$-Tate-Shafarevich group of the motive $\H^m(X_{\nsp K})(r)$ is defined as $\Ker(\alpha^{m,r})$ and often denoted by
$\sha^{(p)}(\H^m(X_{\nsp K})(r))$.
We fix a finite subset $S' \subset B_0$ containing all points of characteristic $p$ and all points where $X$ has bad reduction.

\begin{thm}\label{thm6-1}
Assume that $p \geqq 3$ or $B(\bR)=\emptyset$, and assume further that $\H^3_{\cM}(X,\bZ(r))\{p\}$ is finite. For each $v \in S'$ and $a = 2,3$, we put
\[ e_v^{a,m,r}:=\# \sp \H^a(B_v,\fH^m(X_v,\zp(r))),\]
which is finite by Corollary \ref{cor4-fin}\,{\rm(}2{\rm)}. Then we have
\begin{align*}
\frac{\chi(\alpha^{1,2})}{\,\chi(\alpha^{0,2})\,}
 & =
  \dfrac{\chi(\aj^{3,2}_p)}{\sp \chi(\aj^{2,2}_p) \sp}
  \ac \dfrac{\sp \#\sp \CH_0(X)\{ p \}\sp}{\sp \#\sp \Pic(\OK)\{ p \} \sp}
  \ac \prod_{v \in S'} \ \dfrac{\, e_v^{2,1,2} \ac e_v^{3,0,2}\,}{e_v^{2,0,2} \ac e_v^{3,1,2}}
  \quad  \phantom{\bigg|_{\big|}} & \hbox{\rm ($r=2$)\,}  \\
\frac{\chi(\alpha^{1,r})}{\sp \chi(\alpha^{0,r}) \ac \chi(\alpha^{2,r})\sp}
 & =
  \dfrac{\,\chi(\aj^{3,r}_p)\,}{\,\chi(\aj^{2,r}_p) \,} \ac \#\sp \H^4_{\cM}(X,\bZ(r))\{p\}
  \ac \prod_{v \in S'} \ \dfrac{\, e_v^{2,1,r} \ac e_v^{3,0,r} \ac e_v^{3,2,r}\,}{e_v^{2,0,r} \ac e_v^{2,2,r} \ac e_v^{3,1,r}}
  \quad & \hbox{\rm ($r \geqq 3$),}
\end{align*}
where we put $\chi(f):=\#\sp \Coker(f)/\#\sp \Ker(f)$ for a homomorphism $f : M \to N$ of abelian groups with finite kernel and cokernel.
\end{thm}

See Proposition \ref{cor6-2} for the finiteness of $\H^4_{\cM}(X,\bZ(r))\{p\}$.
The alternating products of local terms $e_v^{a,m,r}$ will be computed in \S\ref{sect7} below.
To prove Theorem \ref{thm6-1}, we first prove Lemma \ref{lem6-3} below as a preparation, which relies on the assumption that $d=2$.
We put
\begin{align*}
\H^1_{\ovf}(B,\fH^m(X,\QpZp(r)))
:=\frac{\,\H^1(B,\fH^m(X,\QpZp(r)))\,}{\,\text{Image of }\, \H^1_{\!f}(K,V^m(r))\,}
\end{align*}
using Corollary \ref{cor5-3}. For each $v \in B_0$, we put
\begin{align*}
\H^1_{\ovf}(B_v,\fH^m(X,\QpZp(r)))
 :=\frac{\,\H^1(B_v,\fH^m(X,\QpZp(r)))\,}
    {\,\text{Image of }\, \H^1_{\!f}(K_v,V^m(r))\,}
\end{align*}
using Corollary \ref{cor4-fin}\,(1).

\begin{lem}\label{lem6-3}
There are canonical isomorphisms of finite $p$-groups
\begin{align}
\label{eq6-4'}
\H^1_{\ovf}(B,\fH^m(X,\QpZp(r))) &\cong \H^2(B,\fH^m(X,\zp(r))), \phantom{\big|_|} \\
\label{eq6-5'}
\H^2(B,\fH^m(X,\QpZp(r))) &\cong \H^3(B,\fH^m(X,\zp(r))), \phantom{}
\intertext{for any $m \geqq 0$ and $r \geqq 2$.
Similarly, there are canonical isomorphisms of finite $p$-groups}
\label{eq6-4}
\H^1_{\ovf}(B_v,\fH^m(X_v,\QpZp(r))) &\cong \H^2(B_v,\fH^m(X_v,\zp(r))), \phantom{\big|_|} \\
\label{eq6-5}
\H^2(B_v,\fH^m(X_v,\QpZp(r))) &\cong \H^3(B_v,\fH^m(X_v,\zp(r))). \phantom{}
\end{align}
for any $m \geqq 0$, $r \geqq 2$ and $v \in B_0$.
Moreover, the groups in \eqref{eq6-4} and \eqref{eq6-5} are zero for any $v \in B_0 \ssm S'$. 
\end{lem}

\begin{pf}
We prove only \eqref{eq6-4'} and \eqref{eq6-5'}, and omit the proof of \eqref{eq6-4} and \eqref{eq6-5}.
We start with the following short exact sequence on $X_\et$, which is a simple case of Proposition \ref{prop:bock}:
\[ 0 \lra \fT_{n'}(2-r)_X \os{\ul {p^n}}{\lra} \fT_{n'+n}(2-r)_X \os{{\mathscr R}^n}\lra \Tn(2-r)_X \lra 0.\]
Concerning this exact sequence of \'etale sheaves, we first prove the following claim:
\begin{itemize}
\item[($\clubsuit$)]
{\it The associated long exact sequence of higher direct image sheaves breaks up into short exact sequences on $B_\et$
\[ \!\!\!\!\! 0 \to R^{2-m}\pi_{X\!/\!B*} \fT_{n'}(2-r)_X \to R^{2-m}\pi_{X\!/\!B*} \fT_{n'+n}(2-r)_X \to R^{2-m}\pi_{X\!/\!B*} \Tn(2-r)_X \to 0  \]
for $m=0,1,2$.}
\end{itemize}
{\it Proof of the claim {\rm($\clubsuit$)}}.\;
We write $\fS_m=\fS_{m,r,n',n}$ for the sequence on $B_\et$ in the display, and prove
that $\fS_m$ is exact for $m=0,1,2$.
We first note that $\fS_2$ is isomorphic to the short exact sequence (of sheaves) on $B_\et$
\begin{align*}
0 \lra \fT_{n'}(2-r)_B \os{\ul {p^n}}{\lra} \fT_{n'+n}(2-r)_B \os{{\mathscr R}^n}\lra \fT_{n}(2-r)_B \lra 0
\end{align*}
by the connectedness of geometric fibers of $\pi : X \to B$ and the assumption that $r \geqq 2$.
The stalks of the sheaves in $\fS_0$ at $v \in B_0$ are zero if $\ch(v)=p$ by \cite{sga4} X.5.2 (and those in $\fS_1$ at $v \in B_0$ are zero if $\ch(v)=p$ and $r >2$ by the proper base-change theorem). Thus it remains to check that the stalk of $\fS_0$ at $\ol x$ is exact for any point $x \in B$ with $\ch(x) \ne p$. Indeed, if $\ch(x) \ne p$, then one can check that the stalk $(R^2\pi_{X\!/\!B*} \fT_{n}(2-r)_X)_{\ol x}$
is isomorphic to the direct sum of copies of $\vL_n(1-r)$ over the set of the irreducible components of $X \times_B \ol x$, by taking a smooth dense open subset $U$ of $X \times_B \ol x$ and comparing the stalk in question with the cohomology of $U$ with compact support.
This completes the proof of ($\clubsuit$). \hfill \qed
\par\medskip
We return to the proof of Lemma \ref{lem6-3}.
From the short exact sequences in the claim ($\clubsuit$) for $n,n' \geqq 1$,
one obtains distinguished triangles in $D(B_\et)$
\[ \fH^m(X,\fT_n(r)) \lra \fH^m(X,\fT_{n+n'}(r)) \lra \fH^m(X,\fT_{n'}(r)) \lra \fH^m(X,\fT_n(r))[1], \]
which yield the following long exact sequence by Proposition \ref{prop1-1} and a standard argument:
\begin{multline*}
\dotsb \to \H^a(B,\fH^m(X,\qp(r))) \to \H^a(B,\fH^m(X,\QpZp(r))) \to \H^{a+1}(B,\fH^m(X,\zp(r)))
  \\ \phantom{\big{|}^{\big{|}}} \to \H^{a+1}(B,\fH^m(X,\qp(r))) \to \dotsb.
\end{multline*}
Now \eqref{eq6-5'} follows from the finiteness of $\H^2(B,\fH^m(X,\QpZp(r)))$ (Theorem \ref{thm5-1}) and the vanishing of $\H^3(B,\fH^m(X,\qp(r)))$ (Theorem \ref{thm:selmer}\,(1)).
Similarly, \eqref{eq6-4'} follows from Theorems \ref{lem:fg}\,(2) and \ref{thm:selmer}\,(2) and the vanishing of $\H^2(B,\fH^m(X,\qp(r)))$.
Finally, the groups on the right hand side of \eqref{eq6-4} and \eqref{eq6-5} are zero for any $v \in B_0 \ssm S'$ by Lemma \ref{lem4-1}.
\end{pf}
\pmn

\begin{pf*}{Proof of Theorem \ref{thm6-1}}
The map $\cl_p^{m,r}$ is bijective for $m=2$ by the finiteness assumption on $\H^3_{\cM}(X,\bZ(r))\{p\}$ (see Proposition \ref{lem6-1}\,(iv)\,$\Ra$\,(iii)), and bijective for $m=3,4$ by Corollary \ref{lem6-0}\,(3), (4).
In particular for $m=2,3$, the map $\aj_p^{m,r}$ is identified with the canonical map
\[ \H^m(X,\zp(r)) \lra \H^1(B,\fH^{m-1}(X,\zp(r))). \]
We put $e^{a,m,r}:=\# \sp \H^a(B,\fH^m(X,\zp(r)))$ for each $a \geqq 2$, $m \geqq 0$ and $r \geqq 2$ with $(a,m,r) \ne (3,2,2)$,
which is finite by Theorems \ref{thm:selmer}\,(1) and \ref{thm5-1}.
One can easily derive an equality
\[ \frac{\,\chi(\aj_p^{3,r}) \,}{\chi(\aj_p^{2,r})}
 = \dfrac{\, e^{2,0,r}\ac e^{2,2,r} \ac e^{3,1,r} \,}{ e^{2,1,r}\ac e^{3,0,r} \ac \#\sp \H^4(X,\zp(r)) } \]
for any $r \geqq 2$, from the spectral sequence \eqref{eq6-1} and the vanishing \eqref{eq6-0}.
Therefore by Corollary \ref{lem6-0}\,(4) and the isomorphisms
\begin{align*}
 \H^2(B,\fH^2(X,\zp(2))) & \os{\text{\eqref{isom:trace}}}\cong \H^2(B,\zp(1)) \cong \Pic(\OK) \otimes \zp \cong \Pic(\OK)\{p\}, \phantom{\big|_|} \\
 \H^3(B,\fH^2(X,\zp(r))) & \,\;\sp \cong \sp\;\, \H^3(B,\zp(r-1)) \cong \H^3(B[p^{-1}],\zp(r-1)) = 0 \qquad (r \geqq 3)
\end{align*}
we are reduced to showing that
\begin{equation}\label{eq6-2}
\chi(\alpha^{m,r})
 = \frac{\, e^{3,m,r}\,}{e^{2,m,r}} \times \prod_{v \in S'} \ \frac{\,e_v^{2,m,r}\,}{e_v^{3,m,r}}
 \qquad \hbox{for \, $^\forall (m,r) \ne (2,2)$, \, $r \geqq 2$.}
\end{equation}

To prove \eqref{eq6-2}, we use the same notation as in the proof of Theorem \ref{thm5-1}, and consider the following commutative diagram with exact rows for $(m,r) \ne (2,2)$ with $r \geqq 2$, where the coefficients $\fH^m(X,\QpZp(r))$ in the upper row and $\fH^m(X_v,\QpZp(r))$ in the lower row are omitted:
{
\begin{equation}\notag
\xymatrix{
\H^1_{\ovf}(B) \, \ar@{^{(}->}@<-1pt>[r]\ar[d] &
\H^1_{\ovf}(K) \ar[r]\ar[d]_{\alpha^{m,r}} &
 \displaystyle \bigoplus_{v \in B_0} \H^2_v(B) \ar@{-}[r]^-{(*)} \ar[d]_{\delta}^{\!\wr} & \qquad\quad \\
\displaystyle \bigoplus_{v \in B_0} \H^1_{\ovf}(B_v) \, \ar@{^{(}->}@<-1pt>[r] &
\displaystyle \bigoplus_{v \in B_0} \H^1_{\ovf}(K_v) \ar[r] &
 \displaystyle \bigoplus_{v \in B_0} \H^2_v(B_v) \ar@{-}[r]^-{(**)} & \qquad\quad \\
& \qquad\quad \ar[r]^-{(*)} & \H^2(B) \ar[r] \ar[d]_{\beta} &
 \H^2(K) \ar[r] \ar[d]_{\gamma}^{\!\wr} &
 \displaystyle \bigoplus_{v \in B_0} \H^3_v(B) \ar[d]_{\delta}^{\!\wr} \\
& \qquad\quad \ar[r]^-{(**)} & \displaystyle \bigoplus_{v \in B_0} \H^2(B_v) \ar[r] &
 \displaystyle \bigoplus_{v \in B_0} \H^2(K_v) \ar[r] &
 \displaystyle \bigoplus_{v \in B_0} \H^3_v(B_v) 
}\end{equation}
}In this diagram, the arrows $\delta$ are bijective as explained in the proof of Theorem \ref{thm5-1}.
The arrow $\gamma$ is bijective by the Hasse principle of Jannsen (\cite{J} p.\ 337, Theorem 3\,(d)) and the fact that $H^m(X_{\ol K},\QpZp(r))$ is divisible.
From the above commutative diagram, we obtain a six-term exact sequence
\begin{multline*}
0 \to \Ker(\alpha^{m,r}) \to \H^1_{\ovf}(B,\fH^m(X,\QpZp(r))) \to
\displaystyle \bigoplus_{v \in B_0} \H^1_{\ovf}(B_v,\fH^m(X_v,\QpZp(r))) \\
 \to \Coker(\alpha^{m,r}) \to \H^2(B,\fH^m(X,\QpZp(r))) \to \bigoplus_{v \in B_0} \H^2(B_v,\fH^m(X_v,\QpZp(r))) \to 0.
\end{multline*}
By Lemma \ref{lem6-3}, this sequence yields an exact sequence of the following from:
{
\begin{multline*}
0 \lra \Ker(\alpha^{m,r})  \lra \H^2(B,\fH^m(X,\zp(r)))  \lra
\displaystyle \bigoplus_{v \in S'} \ \H^2(B_v,\fH^m(X_v,\zp(r))) \\
 \lra \Coker(\alpha^{m,r})
 \lra \H^3(B,\fH^m(X,\zp(r))) \lra \bigoplus_{v \in S'} \ \H^3(B_v,\fH^m(X_v,\zp(r))) \lra 0,
\end{multline*}
}which implies the formula \eqref{eq6-2}.
\end{pf*}

\section{Local terms and zeta values \quad ($\bs{d=2}$)}\label{sect7}
In this section, we compute the local terms $e_v^{2,m,r}$ and $e_v^{3,m,r}$ that appear in Theorem \ref{thm6-1}.
The results in \S\S\ref{sect7-1}--\ref{sect7-2} below were obtained in discussions with Takao Yamazaki.
\par
The setting and the notation remain as in \S\ref{sect6}. In particular, we assume $d=2$.
Put $T^m:=\H^m(X_{\ol K},\zp)$ and $V^m:=T^m \otimes_{\zp} \qp$. 
We further fix the following notation. For a finite place $v$ of $K$,
we write $k_v$ (resp.\ $Y_v$, $Y_{\ol v}$) for the residue field at $v$ (resp.\ $X \otimes_{\OK} k_v$, $X \otimes_{\OK} \ol{k_v}$), and $X_v$ (resp.\ $X_{\ol v}$) for $X \otimes_{\OK} O_v$ (resp.\ $X \otimes_{\OK} O_{\ol v}^\sh$), where $O_v$ (resp.\ $O_{\ol v}^\sh$) denotes the completion of $\OK$ at $v$ (resp.\ the strict henselization of $O_v$ at its maximal ideal).
We put $q_v:=\# k_v$.

\subsection{Comparison with local points}\label{sect7-1}

We first show the following lemma, which refines the case of $q=1$ of Theorem \ref{thm:local-cond} under the assumption that $d=2$.

\begin{lem}\label{lem7-1}
We have
\[ \H^1(B_v,\fH^m(X_v,\zp(r)))=\H^1_{\!f}(K_v,T^m(r)) \]
as subgroups of $\H^1(K_v,T^m(r))$, for any finite place $v$ of $K$,  $m \geqq 0$ and $r \geqq 2$.
\end{lem}

\begin{pf}
Consider a commutative diagram
(see \S\ref{sect3} for the definition of $\H^1_{\ovf}(K_v,V^m(r))$)
\[\xymatrix{
\H^1(K_v,\fH^m(X_v,\zp(r))) \ar@{=}[d] \ar[rr]^-d && \H^2_v(B_v,\fH^m(X_v,\zp(r)))
 \ar[d]^b \\
\H^1(K_v,T^m(r)) \ar[r]^-a & \H^1_{\ovf}(K_v,V^m(r))\,
 \ar@<-1pt>@{^{(}->}[r]^-{d'} & \H^2_v(B_v,\fH^m(X_v,\qp(r))),
}\]
where the arrows $d$ and $d'$ are connecting maps of localization sequences of cohomology of $B_v$, and
the existence and the injectivity of $d'$ is a consequence of Theorem \ref{thm:local-cond} for $q=1$.
The arrow $a$ is the natural map, and we have $\Ker(a) = \H^1_{\!f}(K_v,T^m(r))$ by definition.
On the other hand, since $\H^1_v(B_v,\fH^m(X_v,\zp(r)))=0$ by Proposition \ref{prop:local-str}\,(1), we have
\[ \Ker(d) = \H^1(B_v,\fH^m(X_v,\zp(r))). \]
Thus it remains to check that the arrow $b$ is injective, which follows from the facts that
\[ \; \H^2_v(B_v,\fH^m(X_v,\zp(r))) = 0 \qquad \hbox{if \sp $v|p$ \sp and \sp $r \geqq 3$}
\qquad\qquad\;\; \hbox{(Corollary \ref{cor2-1}\,(1))}
 \]
and that otherwise
\[ \H^2_v(B_v,\fH^m(X_v,\zp(r))) \cong \H^1(k_v,\H^{2-m}(Y_{\ol v},\QpZp(2-r)))^* \qquad \hbox{(\cite{mazur}\,(2.4))} \]
is torsion-free because $\dim(Y_v)=1$ and $\cd(k_v)=1$.
\end{pf}
\par\medskip
The following corollary follows from Proposition \ref{prop:local-str}\,(1), Lemma \ref{lem7-1} and a similar argument as in the proof of Lemma \ref{lem5-1}:
\begin{cor}\label{cor7-0}
We have
\[ \H^1(B,\fH^m(X,\zp(r))) = \H^1_{\!f}(K,T^m(r)) \]
as subgroups of $\H^1(K,T^m(r))$, for any $m \geqq 0$ and $r \geqq 2$.
\end{cor}

\subsection{Comparison with zeta values of the fibers\, (the case $\bs{v\nd p}$)}\label{sect7-2}
In this subsection, we always assume that $v\nd p$ and $r \geqq 2$.
Note that $\H^a(B_v,\fH^m(X_v,\zp(r)))$ is finite for any $(a,m,r)$ by Theorems \ref{lem:fg}\,(1) and \ref{thm:local-cond}, and zero unless $a=0,1,2,3$ and $m=0,1,2$.
We put
\[ e_\vp^{a,m,r}:=\#\sp \H^a(B_v,\fH^m(X_v,\zp(r))) \]
for each $(a,m,r)$. Note that $\zeta(Y_v,r)$ is a non-zero rational number, since $\dim(Y_v)=1$.
Let $|\;|_p$ be the $p$-adic absolute value on $\qp$ such that $|p|_p=p^{-1}$.

\begin{lem}\label{lem7-2}
We have
\[ |\zeta(Y_v,r)|_p^{-1} = \prod_{(a,m)} \ (e_\vp^{a,m,r})^{(-1)^{a+m}}, \]
where $(a,m)$ on the right hand side runs through all pairs with $0 \leqq a \leqq 3$ and $0 \leqq m \leqq 2$.
\end{lem}

\begin{pf}
Let $G_v$ be the absolute Galois group of $k_v$, and
let $T_p$ be a free $\zp$-module of finite rank on which $G_v$ acts continuously and $\zp$-linearly.
Let $\varphi_v \in G_v$ be the arithmetic Frobenius element, and assume that $\varphi_v$ does not have eigenvalue $1$ on $T_p\otimes_{\zp} \qp$.
Then it is well-known that
\begin{equation}\label{eq7-1}
 \#\sp \H^1(k_v,T_p) = \big|\det{}_{\qp}(1-\varphi_v^{-1}\,|\,T_p \otimes_{\zp} \qp)\big|_p^{-1}.
\end{equation}
Now let $\Fr_v$ be the geometric Frobenius operator acting on $\H^i(Y_{\ol v},\qp)$.
We have $\varphi_v= q_v^r \ac \Fr_v^{-1}$ on $\H^i(Y_{\ol v},\qp(r))$, and
{\allowdisplaybreaks
\begin{align*}
|\zeta(Y_v,r)|_p^{-1}
 & = \prod_{i \geqq 0} \;
  \big|\det{}_{\qp}(1-q_v^{-r} \ac \Fr_v\,|\,\H^i(Y_{\ol v},\qp))\big|_p^{(-1)^i}
  & \hbox{(trace formula \cite{G}, \S2)} \\
 & = \prod_{i \geqq 0} \; (\#\sp \H^1(k_v,\H^i(Y_{\ol v},\zp(r))))^{(-1)^{i+1}}
  & \hbox{(by \eqref{eq7-1})} \\
 & \os{\text{($\star$)}}= \prod_{i \geqq 0} \; (\#\sp \H^i(Y_v,\zp(r)))^{(-1)^i}
  & \hbox{(see below)} \\
 & = \prod_{i \geqq 0} \; (\#\sp \H^i(X_v,\zp(r)))^{(-1)^i}
  & \hbox{(proper base change)} \\
 & = \hspace{-3pt}\prod_{(a,m)} \ (e_\vp^{a,m,r})^{(-1)^{a+m}}
  & \hbox{(spectral sequence \eqref{ss:leray2})}
\end{align*}
}as claimed, where the equality ($\star$) follows from the fact that
 $\H^i(Y_{\ol v},\zp(r))^{G_v}=0$ for any $i \geqq 0$ (because $\dim(Y_v)=1$ and $r \geqq 2$).
\end{pf}
\medskip
If $X_v$ is smooth over $O_v$ (and $v \nd p$), then one obtains easily from \eqref{eq7-1} that
\[ \#\sp \H^1_{\!f}(K_v,T^1(r)) = \big|\det{}_{\qp} (1-q_v^{-r} \ac \Fr_v \,|\, \H^1(Y_{\ol v},\qp))\big|_p^{-1}. \]
See also \cite{BK} Theorem 4.1\,(i). The following theorem extends this fact to the general $v\nd p$ case (see also Lemma \ref{lem4-1}):

\begin{thm}\label{thm7-1}
We have $e_\vp^{a,2,r}=1$ for $a=2,3$, and
\[ 
  \frac{\sp \# \sp \H^1_{\!f}(K_v,T^1(r)) \sp}
  {\sp \big|\zeta(Y_v,r)(1-q_v^{1-r})(1-q_v^{-r})\big|_p^{-1} \sp}
 = \frac{\,e_\vp^{2,1,r}\ac e_\vp^{3,0,r}\,}{e_\vp^{2,0,r}\ac e_\vp^{3,1,r}}. \]
\end{thm}

\begin{pf}
We first show that $e_\vp^{a,2,r}=1$ for $a=2,3$. Indeed, we have
\[ \H^a(B_v,\fH^2(X_v,\zp(r))) \os{\eqref{isom:trace}}\cong \H^a(B_v,\zp(r-1)) \cong \H^a(v,\zp(r-1)) = 0 \]
for any $a \geqq 2$. To prove the second assertion, we note the following facts:
\begin{itemize}
\item[(a)]
{\it $e_\vp^{0,m,r}=1$ for any $m \geqq 0$,
 by Proposition \ref{prop:local-str}\,{\rm(}1{\rm)}, Theorem \ref{thm:local-cond} and the fact that $T^m$ is torsion-free.}
\item[(b)]
{\it $e_\vp^{1,m,r} = \#\sp \H^1_{\!f}(K_v,T^m(r))$ by Lemma \ref{lem7-1}.}
\item[(c)]
{\it $e_\vp^{1,0,r}=|1-q_v^{-r}|_p^{-1}$ and $e_\vp^{1,2,r}=|1-q_v^{1-r}|_p^{-1}$,
 by {\rm(}b{\rm)} and {\rm\cite{BK}} Theorem 4.1\,{\rm(}i{\rm)}.}
\end{itemize}
Combining these facts with Lemma \ref{lem7-2}, we have
{\allowdisplaybreaks
\begin{align*}
 &\big|\zeta(Y_v,r) \sp (1-q_v^{1-r})(1-q_v^{-r})\big|_p^{-1} \phantom{|_{\big|}} & \\
 &\quad = \big|(1-q_v^{1-r})(1-q_v^{-r})\big|_p^{-1} \ac 
\frac{\,e_\vp^{1,1,r}\ac e_\vp^{2,2,r}\,}{\,e_\vp^{1,0,r}\ac e_\vp^{1,2,r} \ac e_\vp^{3,2,r}\,} \ac \frac{e_\vp^{2,0,r}\ac e_\vp^{3,1,r}}{\,e_\vp^{2,1,r}\ac e_\vp^{3,0,r}\,}
  \phantom{\Big|_{\big|}}  & \hbox{(Lemma \ref{lem7-2} and (a))} \\
 &\quad = \#\sp \H^1_{\!f}(K_v,\H^1(X_{\ol K},\zp(r))) \ac \frac{e_\vp^{2,0,r}\ac e_\vp^{3,1,r}}{\,e_\vp^{2,1,r}\ac e_\vp^{3,0,r}\,},
 & \hspace{-35pt} \hbox{((b),\, (c),\, $e_\vp^{2,2,r}=e_\vp^{3,2,r}=1$)}
\end{align*}
}which shows the assertion.
\end{pf}

\subsection{Comparison with zeta values of the fibers\, (the case $\bs{v|p}$)}
Let $v$ be a finite place of $K$ dividing $p$.
We assume here that $X_v$ is smooth over $O_v$.
For each $m=0,1,2$, we fix a Haar measure $\mu_\vp^m$ on $H^m_{\dR}(X_{K_v}/K_v)$ such that
\[ \mu_\vp^m(\H^m_{\dR}(X_v/O_v)) = 1, \]
where $\H^m_{\dR}(X_v/O_v)$ denotes the (usual) algebraic de Rham cohomology of $X_v/O_v$. 
Via the exponential isomorphism of Corollary \ref{cor:exp}:
\[ \exp: \H^m_{\dR}(X_{K_v}/K_v) \os{\simeq}\lra \H^1_{\!f}(K_v,V^m(r))\qquad
 \hbox{($r \geqq 2$)}, \]
we regard $\mu_\vp^m$ as a Haar measure on $\H^1_{\!f}(K_v,V^m(r))$.
Let $K_0=K_{v\nsp,0}$ be the fraction field of the Witt ring $W:=W\nsp(k_v)$, and let $\sigma$ be the Frobenius automorphism of $K_0$.
Let $|\;|_p$ be the $p$-adic absolute value on $\bQ$ such that $|p|_p=p^{-1}$.
We prove here a $p$-adic counterpart of Theorem \ref{thm7-1} under some assumptions.

\begin{thm}\label{thm7-2}
Assume that $p-2 \geqq r \geqq 2$ and that $K_v/\qp$ is unramified {\rm(}i.e., $X_v$ is smooth over $\zp${\rm)}.
Then we have $e_\vp^{a,2,r}=1$ for $a=2,3$, and
\[ \frac{\, \mu_\vp^1(\H^1_{\!f}(K_v,T^1(r)))\sp}
  {\sp \big|\zeta(Y_v,r)(1-q_v^{1-r})(1-q_v^{-r})\big|_p^{-1}\sp}
 = \frac{\,e_\vp^{2,1,r}\ac e_\vp^{3,0,r}\,}{e_\vp^{2,0,r}\ac e_\vp^{3,1,r}},
 \]
where we put $e_\vp^{a,m,r}:=\#\sp \H^a(B_v,\fH^m(X_v,\zp(r)))$ for $a \ne 1$.
\end{thm}

To prove this theorem, we first show Lemma \ref{lem7-3} below, which is a $p$-adic analogue of Lemma \ref{lem7-2}
(compare with \cite{FM} Proposition 5.10).
For a continuous homomorphism $\phi : M \to N$ of locally compact groups with finite kernel and with open image, and for a Haar measure $\nu$ on $N$, we define a Haar measure $\nu'$ on $M$ by
\[ \nu'(Z):= \sum_{i=1}^r \ \nu(\phi(Z_i)) \]
for any Borel subset $Z \subset M$, where $Z=Z_1 \amalg Z_2 \amalg \dotsb \amalg Z_r$ is a partition of $Z$ by Borel subsets $Z_1,Z_2,\dotsc,Z_r$ with each $\phi|_{Z_i}$ injective.
We call $\nu'$ the {\it measure induced by} $\nu$ and often denote it by $\nu$.

\begin{lem}\label{lem7-3}
Under the same assumptions as in Theorem \ref{thm7-2}, we have
\[ |\zeta(Y_v,r)|_p^{-1} = \prod_{(a,m)} \ (e_v^{a,m,r})^{(-1)^{a+m}} \]
where $(a,m)$ on the right hand side runs through all pairs with $0 \leqq a \leqq 3$ and $0 \leqq m \leqq 2${\rm;} we put
\[ e_v^{1,m,r}:=\mu_v^m(\H^1(B_v,\fH^m(X,\zp(r))))\]
with $\mu_v^m$ the measure on $\H^1(B_v,\fH^m(X,\zp(r)))$ induced by that on $\H^1_{\!f}(K_v,V^m(r))$.
\end{lem}

\begin{pf}
We first note that $e_\vp^{0,m,r}=1$ for any $m \geqq 0$, by Proposition \ref{prop:local-str}\,(1), Theorem \ref{thm:local-cond} and the fact that $\H^m(X_{\ol K},\zp)$ is torsion-free.
Hence there exists an edge map induced by the spectral sequence \eqref{ss:leray2}
\[ \H^{m+1}(X_v,\zp(r)) \lra \H^1(B_v,\fH^m(X,\zp(r))), \]
which has finite kernel and cokernel by Theorem \ref{thm:local-cond}.
Concerning the Haar measure $\mu_v^m$ on $\H^{m+1}(X_v,\zp(r))$ induced by that on $\H^1(B_v,\fH^m(X,\zp(r)))$, we have
\[ \prod_{i \geqq 0} \; \mu_v^{i-1}(\H^i(X_v,\zp(r)))^{(-1)^i}
 = \prod_{(a,m)} \ (e_\vp^{a,m,r})^{(-1)^{a+m}}\]
by the spectral sequence \eqref{ss:leray2}. It remains to show that
\begin{equation}\label{eq7-2}
 |\zeta(Y_v,r)|_p^{-1} = \prod_{i \geqq 0} \; \mu_v^{i-1}(\H^i(X_v,\zp(r)))^{(-1)^i}.
\end{equation}

By the assumption on $O_v$, it is isomorphic to $W:=W(k_v)$, the ring of Witt vectors in $k_v$.
For each $n \geqq 1$, we put $\Xn:=X_v \otimes_W \Wn$, and let $\cS_n(r)_{X_v}$ be the syntomic complex associated with the smooth scheme $X_v$ over $W=O_v$.
Let $p(r)\sp \Omega_{\Xn/\Wn}^\bullet$ (resp.\ $p(r)\sp \Omega_{X_v/W}^\bullet$) be the subcomplex
\[  p^r \ac \cO_{\Xn} \os{d}\lra p^{r-1} \ac \Omega_{\Xn/\Wn}^1 \qquad
 \big(\hbox{resp.\ $p^r \ac \cO_{X_v} \os{d}\lra p^{r-1} \ac \Omega_{X_v/W}^1$}\big) \]
of the de Rham complex $\Omega_{\Xn/\Wn}^\bullet$ (resp.\ $\Omega_{X_v/W}^\bullet$).
We note the following facts:
\begin{itemize}
\item[(a)]
{\it There exists an isomorphism
\[ \big(p(r)\sp \Omega_{\Xn/\Wn}^\bullet\big)_n[-1] \cong  (\cS_n(r)_{X_v})_n\]
for any $r$ with $2 \leqq r < p$ in
the derived category of complexes of pro-sheaves on $(Y_v)_\et$ by} \cite{BEK} {\it Theorem 5.4.}
\item[(b)]
{\it The Euler characteristic
\begin{align*}
 \chi(X_v,\Omega_{X_v/W}^\bullet/p(r)\Omega_{X_v/W}^\bullet)
&:= \prod_{i \geqq 0} \ (\#\sp \H^i(X_v,\Omega_{X_v/W}^\bullet/p(r)\sp \Omega_{X_v/W}^\bullet))^{(-1)^i} \\
&\,= \prod_{(a,b)} \ (\#\sp \H^a(Y_v,\Omega_{Y_v/k_v}^b))^{(-1)^{a+b}(r-b)}
\end{align*}
agrees with $|\zeta(Y_v,r)|_p^{-1}$ for any $r \geqq 2$} (\cite{milne:value} {\it Theorem 0.1}).
\item[(c)]
{\it We have $\cS_n(r)_{X_v} \cong i^*\Tn(r)$ in $D(Y_v,\Ln)$ for any $r$ with $r < p-1$ and any $n \geqq 1$ {\rm(\cite{Ku}}\,p.\ 275, Theorem{\rm)},
where $i$ denotes the closed immersion $Y_v \hra X_v$.}
\end{itemize}
By these facts, we have
{\allowdisplaybreaks
\begin{align*}
 |\zeta(Y_v,r)|_p^{-1}
 & = \chi(X_v,\Omega_{X_v/W}^\bullet/p(r)\Omega_{X_v/W}^\bullet) \phantom{\big|_|}
 & \hbox{(by (b))} \\
 & = \prod_{i \geqq 0} \
 \frac{\mu_v^i(\H^i_{\dR}(X_v/W))^{(-1)^i}}{\,\mu_v^i(\H^i(X_v,p(r)\sp \Omega_{X_v/W}^\bullet))^{(-1)^i}} \\
 & = \prod_{i \geqq 0} \ \mu_v^i(\H^i(X_v,p(r)\sp \Omega_{X_v/W}^\bullet))^{(-1)^{i+1}}
 & \hbox{($\mu_v^i(\H^i_{\dR}(X_v/W))=1$)} \\
 & = \prod_{i \geqq 0} \ \mu_v^i(\H^{i+1}(X_v,\zp(r)))^{(-1)^{i+1}}
 & \hbox{(by (a), (c)).}
\end{align*}
Thus we obtain \eqref{eq7-2} and Lemma \ref{lem7-3}.
}
\end{pf}
\par\bigskip

\begin{pf*}{Proof of Theorem \ref{thm7-2}}
We first show that $e_\vp^{a,2,r}=1$ for any $a \geqq 2$. Indeed, we have
\[ \H^a(B_v,\fH^2(X_v,\zp(r))) \os{\eqref{isom:trace}}\cong \H^a(B_v,\zp(r-1)). \]
If $r=2$, then the last group is zero for any $a \geqq 2$ because $\H^a(B_v,\Gm)=0$ for any $a \geqq 1$.
On the other hand, if $r \geqq 3$, then by the Tate duality, we have
\[ \H^a(B_v,\zp(r-1)) \cong \H^a(K_v,\zp(r-1)) \cong \H^{2-a}(K_v,\qp/\zp(2-r))^*, \]
which is zero for any $a \geqq 2$ by the assumptions on $K_v$ and $p$.
Noting that
\begin{itemize}
\item[(a$^+$)]
{\it $e_\vp^{1,0,r}=|1-q_v^{-r}|_p^{-1}$ and $e_\vp^{1,2,r}=|1-q_v^{1-r}|_p^{-1}$
 by} \cite{BK} {\it Theorem 4.2 for $V=\qp(r)$ and $\qp(r-1)$,
 and again by the assumptions on $K_v$ and $p$,}
\end{itemize}
one obtains the second assertion from the same computations as in Theorem \ref{thm7-1}.
\end{pf*}

\section{Global points and zeta values \quad ($\bs{d=2}$)}\label{sect8}
The setting and the notation remain as in \S\ref{sect6} (in particular, $d=2$).
Put $T^m:=\H^m(X_{\ol K},\zp)$ and $V^m:=T^m \otimes_{\zp} \qp$. 
In this section, we relate the formula in Theorem \ref{thm6-1} with zeta values
assuming Conjecture \ref{conj7-1} below for the motives $\H^m(X_{\nsp K})(r)$ with $m=0,1,2$, a weak version of $p$-Tamagawa number conjecture \cite{BK} \S5.
Let $S'$ be a finite set of closed points of $B$ containing all points of characteristic $p$, and all points where $X$ has bad reduction.
For $m=0,1,2$ and $r \geqq 2$ with $(m,r) \ne (2,2)$, we put
\[ \tL_{S'}(\H^m(X_{\nsp K}),r):=\!\! \prod_{v \in B_0 \ssm S'}
 \det(1-q_v^{-r}\ac \Fr_v\,|\,V^m)^{-1}. \]
This infinite product on the right hand side converges, because $m-2r \leqq -3$.
Let $\bZ_{(p)}$ be the localization of $\bZ$ at the prime ideal $(p)$.

\subsection{$\bs{p}$-Tamagawa number conjecture}

\begin{conj}[Bloch-Kato]\label{conj7-1}
For any $m = 0,1,2$ and $r \geqq 2$ with $(m,r) \ne (2,2)$, there exists a finite-dimensional $\bQ$-subspace $\Phi^{m,r}=\Phi^{m,r}_p$ of the $\bQ$-vector space
\[ \H^{m+1}_{\cM}(X_{\nsp K},\bQ(r))_\bZ:= \Image\big(\H^{m+1}_{\cM}(X,\bQ(r)) \to \H^{m+1}_{\cM}(X_{\nsp K},\bQ(r))\big) \]
satisfying the following conditions {\rm (i)} and {\rm (ii):}
\begin{itemize}
\item[{\rm(i)}]
The $p$-adic Abel-Jacobi map
\[ \H^{m+1}_{\cM}(X_{\nsp K},\bQ(r)) \lra \H^1(K,V^m(r)) \]
induces an isomorphism $\Phi^{m,r} \otimes \qp \cong \H^1_{\!f}(K,V^m(r))$,
and Beilinson's regulator map to the real Deligne cohomology
\[ \H^{m+1}_{\cM}(X_{\nsp K},\bQ(r)) \lra \H^{m+1}_{\cD}(X_{\nsp/\bR},\bR(r)) \]
induces an isomorphism
$\Phi^{m,r} \otimes \bR \cong \H^{m+1}_{\cD}(X_{/\bR},\bR(r))$.
\item[{\rm(ii)}]
We define $A^{m,r}_p(K)$, the group of $p$-global points as the pull-back of $\Phi^{m,r}$ under the natural map
\[ \H^1_{\!f}(K,T^m(r)) \lra \H^1_{\!f}(K,V^m(r)) \cong \Phi^{m,r} \otimes \qp, \]
which is a finitely generated $\bZ_{(p)}$-module. We further fix an $\OK$-lattice $L^m$ of the de Rham cohomology $\H^m_{\dR}(X_{\nsp K}/K)$, and define a number $R_\Phi^{m,r} \in \bR^\times/\bZ_{(p)}^\times$ to be the volume of the space
\[  \H^{m+1}_{\cD}(X_{\nsp/\bR},\bZ_{(p)}(r))/{\text{Image of }} A^{m,r}_p(K) \]
with respect to $L^m$. 
See Remark {\rm\ref{rem8-1}\,(1)} below for an explicit description of the Deligne cohomology $\H^{m+1}_{\cD}(X_{\nsp/\bR},\bZ_{(p)}(r))$. On the other hand, for each $v \in B_0$ we put
\[ A^{m,r}_p(K_v):=\H_{\!f}^1(K_v,T^m(r)), \]
which we call the group of $p$-local points at $v$. Then we have
\begin{align}\label{eq7-3}
\tL_{S'}(\H^m(X_{\nsp K}),r) \equiv \chi(\alpha^{m,r})^{-1} \ac R^{m,r}_\Phi \ac \prod_{v \in S'} \ \mu_v^m(A^{m,r}_p(K_v))
 \;\; \mod \bZ_{(p)}^\times,
\end{align}
where $\mu_v^m$ for $v \nd p$ means the cardinality, and $\mu_v^m$ for $v|p$ denotes the Haar measure on \newline
$A^{m,r}_p(K_v)$ constructed from that on $\H^m_{\dR}(X_{\nsp K_v}/K_v)$ such that $\mu_v^m(L^m \otimes_{\OK} O_v)=1${\rm;} see \eqref{eq6+3} for the map $\alpha^{m,r}$.
\end{itemize}
\end{conj}

\begin{rem}\label{rem8-1}
{\rm
\begin{enumerate}
\item[(1)]
The map $A^{m,r}_p(K) \to \H^{m+1}_{\cD}(X_{\nsp/\bR},\bZ_{(p)}(r))$ induced by the regulator map is injective, by the condition (i) for $\Phi^{m,r}$ and \cite{BK} Lemma 5.10. Here
\[ \H^{m+1}_{\cD}(X_{\nsp/\bR},\bZ_{(p)}(r))
 = \bigg(\frac{\,\H^{m}_{\dR}(X/\bZ)\otimes \bC\,}{\,\H^{m}_\sing(X \otimes_{\bZ} \bC,(2\pi i)^r \ac \bZ_{(p)})\,}\bigg)^{\!\!+} \]
for any $m = 0,1,2$ and $r \geqq 2$, by definition.
\item[(2)]
The product on the right hand side of \eqref{eq7-3} is independent of the choice of $L^m$.
\item[(3)]
Conjecture \ref{conj7-1} for $m=0$ (resp.\ $m=2$) implies that
\[ \zeta_K(r) \equiv \chi(\alpha^{0,r})^{-1} \ac R^{0,r}_\Phi \qquad
 \hbox{(resp.\ $\zeta_K(r-1) \equiv \chi(\alpha^{0,r-1})^{-1} \ac R^{0,r-1}_\Phi$)} \]
modulo $\bZ_{(p)}^\times$ if $r \geqq 2$ (resp.\ $r \geqq 3$) and $p$ is unramified in $K$.
Here we have used the fact (c) in the proof of Theorem \ref{thm7-1} for all $v \nd p$ belonging to $S'$, and the fact (a$^+$) in the proof of Theorem \ref{thm7-2}.
See also \cite{FM} \S5.8.3.
\item[(4)]
We have $R_\Phi^{m,r}=1$ for any $m \geqq 3$, because $\H^{m+1}_{\cD}(X_{\nsp/\bR},\bZ_{(p)}(r))$ is zero for such $m$'s.
\item[(5)]
If $(m,r)=(2,2)$, there exists a $\bQ$-subspace $\Phi^{2,2}$ of $\H^3_{\cM}(X_{\nsp K},\bQ(2))_\bZ$ which is isomorphic to $\H^1_{\cM}(B,\bQ(1))$ under the push-forward map
\[ \H^3_{\cM}(X_{\nsp K},\bQ(2)) \lra \H^1_{\cM}(\Spec(K),\bQ(1)) \cong K^\times \otimes \bQ. \]
Indeed, by a standard norm argument, the push-forward map
\[ \H^3_{\cM}(X,\bQ(2)) \lra \H^1_{\cM}(B,\bQ(1)) \cong \OK^\times \otimes \bQ \]
is surjective, and there is a $\bQ$-subspace $\wt{\Phi}{}^{2,2} \subset \H^3_{\cM}(X,\bQ(2))$ which maps bijectively onto $\H^1_{\cM}(B,\bQ(1))$. One can define a desired space $\Phi^{2,2}$ by
\[ \Phi^{2,2}:=\Image\big(\wt{\Phi}{}^{2,2}\to\H^3_{\cM}(X_{\nsp K},\bQ(2))\big). \]
By this construction of $\Phi^{2,2}$, we have
\[ \Phi^{2,2} \otimes \qp \cong \H^1_{\!f}(K,V^2(2)) \, (= \H^1_{\!f}(K,\qp(1))). \]
See also Corollary \ref{lem6-0}\,(3). For $(m,r)=(2,2)$, we will use the classical class number formula instead of \eqref{eq7-3}, later in Theorem \ref{thm7-3} below.
\end{enumerate}}
\end{rem}

\begin{prop}\label{prop7-1}
Let $r$ be an integer, and let $p$ be a prime number.
Assume all the following conditions{\rm:}
\begin{itemize}
\item[{\rm (i)}]
$p-2 \geqq r \geqq 2$.
\item[{\rm (ii)}]
For any $v \in B_0$ dividing $p$, $v$ is absolutely unramified and $X$ has good reduction at $v$.
\item[{\rm (iii)}]
Conjecture \ref{conj7-1} holds for $m=0,1$ {\rm(}resp.\ $m=0,1,2${\rm)}, if $r=2$ {\rm(}resp.\ $r \geqq 3${\rm)}.
\end{itemize}
Then the equivalent conditions {\rm(i)--(v)} of Proposition \ref{lem6-1} are satisfied for $m=1,2$ {\rm(}resp.\ $m=1,2,3${\rm)}, if $r=2$ {\rm(}resp.\ if $r \geqq 3${\rm)}. Moreover, we have
{\allowdisplaybreaks
\begin{align*}
& \displaystyle \us{s=2}{\Res} \ \zeta(X,s) \equiv
\us{s=1}{\Res} \ \zeta_K(s) \ac
\frac{\, \chi(\aj_p^{3,2}) \ac \#\sp \CH_0(X) \ac R^{0,2}_\Phi  \,}{\sp \chi(\aj_p^{2,2}) \ac \#\sp \Pic(\OK) \ac R^{1,2}_\Phi \sp}
\;\; \mod \bZ_{(p)}^\times \phantom{\big|_{\Big|}} & \hbox{{\rm (}$r = 2${\rm)}} \\
& \zeta(X,r) \equiv
\frac{\,\chi(\aj_p^{3,r}) \ac \#\sp \H^4_{\cM}(X,\bZ(r))\{p\} \ac R^{0,r}_\Phi \ac R^{2,r}_\Phi \,}{\sp \chi(\aj_p^{2,r}) \ac R^{1,r}_\Phi \sp}
\;\; \mod \bZ_{(p)}^\times & \hbox{{\rm (}$r \geqq 3${\rm)}}
\end{align*}
}
\end{prop}

\begin{pf}
The first assertion is obvious.
For any $r \geqq 2$, we have
{\allowdisplaybreaks
\begin{align*}
& \us{s \to r}{\text {lim}} \ \dfrac{\zeta(X,s)}{\sp \zeta_K(s)\zeta_K(s-1) \sp}\phantom{\bigg|_|} =
 \dfrac{1}{\tL_{S'}(\H^1(X_{\nsp K}),r)}
 \ac \prod_{v \in S'} \dfrac{\zeta(Y_v,r)}{\sp (1-q_v^{-r})^{-1}(1-q_v^{1-r})^{-1} \sp} \\
 &\quad \equiv
 \dfrac{\sp \chi(\alpha^{1,r}) \sp}{ R^{1,r}_\Phi } \ac
 \prod_{v \in S'} \ \frac{1}{\sp \mu_v^1(A^{1,r}_p(K_v)) \sp}
 \ac \prod_{v \in S'} \frac{\, e_\vp^{2,0,r}\ac e_\vp^{3,1,r}\ac \mu_v^1(A^{1,r}_p(K_v)) \,}{\,e_\vp^{2,1,r}\ac e_\vp^{3,0,r}\,} \;\;  \mod  \bZ_{(p)}^\times \\
 &\quad =
 \dfrac{\sp \chi(\alpha^{1,r}) \sp}{ R^{1,r}_\Phi }
 \ac \prod_{v \in S'} \frac{\sp e_\vp^{2,0,r}\ac e_\vp^{3,1,r}\sp}{\,e_\vp^{2,1,r}\ac e_\vp^{3,0,r}\,}
\end{align*}
}by the assumptions (i)\sp--\sp(iii) for $m=1$ and Theorems \ref{thm7-1} and \ref{thm7-2} (see Remark \ref{rem8-1}\,(2)).
Hence for $r=2$, we have
{\allowdisplaybreaks
\begin{align*}
\us{s=2}{\Res} \ \zeta(X,s) & \equiv
\us{s=1}{\Res} \ \zeta_K(s) \ac \frac{R^{0,2}_\Phi \ac \chi(\alpha^{1,2})}{\sp \chi(\alpha^{0,2}) \ac R^{1,2}_\Phi \sp} \ac \! \prod_{v \in S'} \frac{\sp e_\vp^{2,0,2}\ac e_\vp^{3,1,2}\sp}{\,e_\vp^{2,1,2}\ac e_\vp^{3,0,2}\,} \;\;  \mod  \bZ_{(p)}^\times \\
& = \us{s=1}{\Res} \ \zeta_K(s) \ac \frac{\, \chi(\aj_p^{3,2}) \ac \#\sp \CH_0(X) \ac R^{0,2}_\Phi  \,}{\sp \chi(\aj_p^{2,2}) \ac \#\sp \Pic(\OK) \ac R^{1,2}_\Phi \sp}
\end{align*}
}by the assumption (iii) for $m=0$ and Theorem \ref{thm6-1}.
See also Remark \ref{rem8-1}\,(3).
Similarly for any $r \geqq 3$, we have
{\allowdisplaybreaks
\begin{align*}
\zeta(X,r) & \equiv \frac{R^{0,r}_\Phi \ac \chi(\alpha^{1,r}) \ac R^{2,r}_\Phi}{\sp \chi(\alpha^{0,r}) \ac R^{1,r}_\Phi \sp \ac \chi(\alpha^{2,r}) } \ac \! \prod_{v \in S'} \frac{\sp e_\vp^{2,0,r}\ac e_\vp^{3,1,r}\sp}{\,e_\vp^{2,1,r}\ac e_\vp^{3,0,r}\,} \;\; \mod  \bZ_{(p)}^\times \\
& = \dfrac{\sp\chi(\aj^{3,r}_p)\ac \#\sp \H^4_{\cM}(X,\bZ(r))\{p\} \ac R^{0,r}_\Phi \ac R^{2,r}_\Phi\sp}{\,\chi(\aj^{2,r}_p)\ac R^{1,r}_\Phi \,}
\end{align*}
}as claimed.
\end{pf}

\subsection{Zeta value formula without \'etale cohomology}
Let $p$ be an arbitrary prime number.
Assuming Conjecture \ref{conj7-1} for $p$, we define a number $R^{m,r}_{\cM}=R^{m,r}_{\cM\nsp\nsp,\sp p} \in \bR^\times/\bZ_{(p)}^\times$ ($m \geqq 0$, $r \geqq 2$) as follows.
We first take the inverse image $\wt{A}{}_p^{m,r}$ of $A^{m,r}_p(K)$ under the composite map
\[ \H^{m+1}_{\cM}(X,\bZ(r)) \otimes \bZ_{(p)} \to
 \H^{m+1}_{\cM}(X_{\nsp K},\bZ(r)) \otimes \bZ_{(p)} \to
 \H^1(K,T^m(r)), \]
where for $(m,r)=(2,2)$, $A^{2,2}_p(K)$ is considered with respect to $\Phi^{2,2}$ constructed in Remark \ref{rem8-1}\,(5).
Since $A_p^{m,r}(K)$ is finitely generated over $\bZ_{(p)}$, the canonical map $\wt{A}{}_p^{m,r} \to A_p^{m,r}(K)$ induces a homomorphism
\[ c^{m,r} : \ol{A}{}_p^{m,r} :=\wt{A}{}_p^{m,r}\!\big/\nsp\big(\wt{A}{}_p^{m,r}\big)_\Div \lra A_p^{m,r}(K). \]
Here `$\Div$' means the maximal divisible subgroup.
This map fits into a commutative diagram
\begin{equation}\label{eq7-4}
\xymatrix{
 \ol{A}{}_p^{m,r}\otimes \zp
 \ar[r]^-{\gamma^{m,r}} \ar[d]_{c^{m,r}\otimes \id} &
 \H^{m+1}_{\cM}(X,\bZ(r))\,\wh{\otimes}\,\zp \ar[r]^-{\aj^{m+1,r}_p} &
 \H^1(B,\fH^m(X,\zp(r))) \ar@{=}[d] \\
A_p^{m,r}(K) \otimes \zp \, \ar@{=}[rr]^-{\sim} && \H^1_{\!f}(K,T^m(r)),
}\end{equation}
where $\gamma^{m,r}$ denotes the natural map.
See Corollary \ref{cor7-0} for the right vertical equality.
\begin{lem}\label{lem7-4}
Assume that $p \geqq 3$ or $B(\bR)=\emptyset$, and that Conjecture \ref{conj7-1} holds. Then
$\gamma^{m,r}$ and $c^{m,r}$ have finite cokernel for any $m \geqq 0$ and $r \geqq 2$.
\end{lem}
\begin{pf}
$\Coker(c^{m,r})$ is finite, because it is finitely generated over $\bZ_{(p)}$ and torsion by the definition of $\wt{A}{}_p^{m,r}$ (essentially by Conjecture \ref{conj7-1}).
The map $\gamma^{m,r}$ has finite cokernel as well, because $c^{m,r}\otimes \id_{\zp}$ has finite cokernel and $\aj^{m+1,r}_p$ has finite kernel by Corollary \ref{lem6-0}.
\end{pf}
\bigskip\noindent
By the finiteness of $\Coker(c^{m,r})$, we define $R_{\cM}^{m,r} \in \bR^\times/\bZ_{(p)}^\times$ to be the volume of the space
\[ 
\begin{cases}
 \H^{m+1}_{\cD}(X_{\nsp/\bR},\bZ_{(p)}(r))/\sp {\text{Image of }} \ol{A}{}_p^{m,r}
  \qquad \phantom{\big|_{\big|}} & \hbox{(for $(m,r) \ne (2,2)$)} \\
 \wt{\H}{}^3_{\cD}(X_{\nsp/\bR},\bZ_{(p)}(2))/\sp {\text{Image of }} \ol{A}{}_p^{2,2} &
 \hbox{(for $(m,r)=(2,2)$)}
\end{cases}
\]
with respect to $L^m$ that we fixed in Conjecture \ref{conj7-1}, where $\wt{\H}{}^3_{\cD}(X_{\nsp/\bR},\bZ_{(p)}(2))$ denotes the kernel of the canonical trace map
\[ \tr : \wt{\H}{}^3_{\cD}(X_{\nsp/\bR},\bZ_{(p)}(2)) \lra \bR. \]
We have $R_{\cM}^{m,r}=1$ for any $m \geqq 3$ by definition.

\begin{prop}\label{rem7-2}
{\rm
If $p \geqq 3$ or $B(\bR)=\emptyset$, then $\gamma^{3,r}$ is bijective for any $r \geqq 2$.
}
\end{prop}
\begin{pf}
Since $T^3=0$, we have $A_p^{3,r}(K)=0$ and
\[ \ol{A}{}_p^{3,r} = \H^4_{\cM}(X,\bZ(r))\otimes \bZ_{(p)}\big/ \nsp (\H^4_{\cM}(X,\bZ(r))\otimes \bZ_{(p)})_\Div \]
by definition. Since $\H^4_{\cM}(X,\bZ(r))\{p\}$ is finite by Corollary \ref{lem6-0}\,(4), the natural maps
\[ \H^4_{\cM}(X,\bZ(r))\{p\} \lra \ol{A}{}_p^{3,r} \lra \H^4_{\cM}(X,\bZ(r)) \, \wh{\otimes} \, \zp \]
are injective, and moreover bijective by Corollary \ref{lem6-0}\,(4), which shows the assertion.
\end{pf}

\begin{thm}\label{thm7-3}
Under the same assumptions as in Proposition \ref{prop7-1}, assume further that
\begin{itemize}
\item[{\rm (iv)}]
$\gamma^{m,r}$ of \eqref{eq7-4} is bijective for any $m=0,1,2$.
\end{itemize}
Then $c^{m,r}$ has finite kernel for any $m=0,1,2,3$, and we have
\begin{align*}
& \zeta^*(X,r) \equiv
\prod_{m=0}^3 \ \bigg( \frac{\, R^{m,r}_{\cM}\,}{\,\#\Ker(c^{m,r}) \,}\bigg)^{\!\!(-1)^m}
    \;\; \mod \bZ_{(p)}^\times,
\end{align*}
where $\zeta^*(X,r)$ denotes $\us{s=2}\Res \ \zeta(X,s)$ {\rm(}resp.\ $\zeta(X,r)${\rm)}
 if $r=2$ {\rm(}resp.\ $r \geqq 3${\rm)}.
\end{thm}

\begin{rem}\label{rem7-3}
{\rm
A stronger version of Conjecture \ref{conj7-1} asserts that
\begin{itemize}
\item[(h\nsp 1)]
{\it The $\bQ$-space $\Phi^{m,r}$ agrees with $\H^{m+1}_{\cM}(X_{\nsp K},\bQ(r))_\bZ$\, for any $m=0,1,2$ and $r \geqq 2$.}
\end{itemize}
The above condition (iv) holds true, under this stronger hypothesis and the following variant of Bass' conjecture (cf.\ \cite{Ba}):
\begin{itemize}
\item[(h2)]
 {\it $\H^{m+1}_{\cM}(X,\bZ(r))$ is finitely generated for any $m=0,1,2$ and $r \geqq 2$.}
\end{itemize}
Under the hypotheses (h\nsp 1) and (h2), $\Ker(c^{m,r})$ agrees with the $p$-primary torsion part of the kernel of the regulator map
\[ \reg_{\cD}^{m+1,r} : \H^{m+1}_{\cM}(X,\bZ(r)) \lra \H^{m+1}_{\cD}(X_{\nsp/\bR},\bZ(r)) \]
by Remark \ref{rem8-1}\,(1), and $R^{m,r}_{\cM}$ is exactly the volume of its cokernel (modulo $\bZ_{(p)}^\times$).
}
\end{rem}
\smallskip
\begin{pf*}{Proof of Theorem \ref{thm7-3}}
The hypothesis (iv) and Proposition \ref{rem7-2} imply that $\gamma^{m,r}$ is injective for $m = 0, 1, 2, 3$. Hence the finiteness of $\ker(\aj_p^{m+1,r})$ (see Corollary \ref{lem6-0}) implies that
$c^{m,r}\otimes \id$ in \eqref{eq7-4} has finite kernel.
Thus $c^{m,r}$ has finite kernel, because $\zp$ is faithfully flat over $\bZ_{(p)}$.
\par
We rewrite the number on the right hand side in the formulas in Proposition \ref{prop7-1}.
By the classical class number formula, we have
\[ \us{s=1}{\Res} \ \zeta_K(s) = \vol(\Coker(\varrho)) \ac \# \sp \Pic(\OK), \]
where $\varrho=\varrho_K$ denotes the regulator map to (the reduced part of) the integral Deligne cohomology
\[ \varrho : \OK^\times \lra \wt{\H}{}^1_\cD(B_{\nsp/\bR},\bZ(1))
 := \Ker\big(\tr: \H^1_{\cD}(B_{\nsp/\bR},\bZ(1)) \to \bR\big) \]
and the volume of $\Coker(\varrho)$ has been taken with respect to $\OK \subset K=\H^0_{\dR}(\Spec(K)/K)$
(note that $\varrho$ is injective).
To prove the formula in Theorem \ref{thm7-3}, it remains to check
\begin{equation}\label{eq7-5}
 \frac{\, R^{m,r}_{\cM}\,}{\,\#\Ker(c^{m,r})\,} =
\begin{cases}
R^{0,r}_\Phi \qquad \phantom{\Big|} & \hbox{($m=0$)} \\
\chi(\aj_p^{2,r}) \ac R^{1,r}_\Phi \qquad\phantom{\Big|} & \hbox{($m=1$)} \\
\chi(\aj_p^{3,2}) \ac \vol(\Coker(\varrho)) \qquad \phantom{\Big|} & \hbox{($(m,r)=(2,2)$)} \\
\chi(\aj_p^{3,r}) \ac R^{2,r}_\Phi \qquad \phantom{\Big|} & \hbox{($m=2$, $r \geqq 3$)} \\
(\#\sp \CH_0(X)\{p\})^{-1} \qquad \phantom{\Big|} & \hbox{($(m,r)=(3,2)$)} \\
(\#\sp \H^4_{\cM}(X,\bZ(r))\{p\})^{-1} \qquad \phantom{\Big|} & \hbox{($m=3$, $r \geqq 3$)}
\end{cases}
\end{equation}
We have
\begin{equation}\label{eq7-6}
 \Ker(\aj_p^{m,r})=\Ker(c^{m,r}) \quad \hbox{and} \quad \Coker(\aj_p^{m,r}) \cong \Coker(c^{m,r})
\end{equation}
for any $(m,r)$ by the diagram \eqref{eq7-4}, the hypothesis (iv) and Proposition \ref{rem7-2}.
This fact implies \eqref{eq7-5} for $m=0,1,2$ with $(m,r) \ne (2,2)$.
See also Proposition \ref{lem6-1} and Corollary \ref{lem6-0}\,(1) for the fact that $\chi(\aj_p^{1,r})=\#\sp\Ker(c^{1,r})=1$.
The formula \eqref{eq7-5} for $m=3$ follows from \eqref{eq7-6} and the fact that $R^{m,r}_{\cM}=1$ for $m \geqq 3$.
Finally, noting that $\gamma^{2,2}$ is bijective by assumption, consider the diagram \eqref{eq7-4} for $(m,r)=(2,2)$:
\begin{equation*}
\xymatrix{
 \ol{A}{}_p^{2,2}\otimes \zp \ar[r]^{c^{2,2}\otimes \id}\ar[d]_-{\aj^{3,2}_p}
 & A_p^{2,2}(K) \otimes \zp \ar@{=}[r]^-\sim \ar@{=}[d]^\wr & \OK^\times \otimes \zp \ar@{=}[d]^\wr \\
 \H^1(B,\fH^2(X,\zp(2))) \ar@{=}[r] & \H^1_{\!f}(K,T^2(2)) \ar@{=}[r]^-\sim & \H^1_{\!f}(K,\zp(1)),
}\end{equation*}
which shows \eqref{eq7-5} for $(m,r)=(2,2)$. This completes the proof.
\end{pf*}

\subsection*{Acknowledgments}
The key idea of the fundamental results of this paper is based on the joint works \cite{JSS}, \cite{SS} with Profs.\ Uwe Jannsen and Shuji Saito.
The first draft of \S\S3\sp--\sp6.1 of this paper was written while the author stayed at University of Southern California from October 2001 to September 2003 supported by JSPS Postdoctoral Fellowship for Research Abroad.
He would like to thank Profs.\ Wayne Raskind and Thomas Geisser for fruitful discussions on arithmetic duality again.
The author expresses his gratitude to Takao Yamazaki for stimulating discussions and comments to the research for this article. The results in \S\S\ref{sect7-1}--\sp\ref{sect7-2} were obtained in a joint work with him in 2019.
He is also grateful to Profs.\ Takeshi Saito and Yukiyoshi Nakkajima, who gave him valuable comments to his talks on earlier versions of this research.
Thanks are also due the referee, who read the manuscript of this paper very carefully and gave the author constructive comments.

\bigskip
\noindent
Department of Mathematics, Chuo University
\par\noindent
1-13-27 Kasuga, Bunkyo-ku,
Tokyo 112-8551,
JAPAN

\smallskip

\noindent
{\it Email} : \textbf{kanetomo@math.chuo-u.ac.jp}
\end{document}